\documentclass[preprint]{elsarticle}
\usepackage[utf8]{inputenc}
\usepackage{bm}
\usepackage{color}
\usepackage{graphicx}
\usepackage{amsmath,amssymb}
\usepackage[bottom=2cm,top=2cm,left=2cm,right=3cm]{geometry}
\usepackage{subfigure}
\usepackage{longtable,tabularx}
\usepackage[table]{xcolor}
\usepackage{booktabs}
\usepackage{sidenotes}
\usepackage{url}
\usepackage[normalem]{ulem}
\newcounter{lnote}

\usepackage[ruled]{algorithm2e}

\SetKwProg{Fn}{}{}{end} 

\newcommand{\aalpha}{\bm{\alpha}}
\newcommand{\bbeta}{\bm{\beta}}
\newcommand{\yy}{\mathbf{y}}
\newcommand{\ii}{\mathbf{i}}
\newcommand{\jj}{\mathbf{j}}
\newcommand{\ee}{\mathbf{e}}
\newcommand{\dd}{\, \mathrm{d}}
\newcommand{\MISCsurr}{\mathcal{S}}
\newcommand{\SRBFsurr}{\mathcal{S}}
\newcommand{\MISCquad}{\mathcal{R}}

\begin{document}

\begin{frontmatter}	
	
\title{Comparing Multi-Index Stochastic Collocation and Multi-Fidelity Stochastic Radial Basis Functions for Forward Uncertainty Quantification of Ship Resistance}

\author[imati]{Chiara Piazzola\corref{cor1}}
\ead{chiara.piazzola@imati.cnr.it}
\author[imati]{Lorenzo Tamellini} 
\ead{tamellini@imati.cnr.it}
\author[inm]{Riccardo Pellegrini}
\ead{riccardo.pellegrini@inm.cnr.it}
\author[inm]{Riccardo Broglia}
\ead{riccardo.broglia@cnr.it}
\author[inm]{Andrea Serani}
\ead{andrea.serani@cnr.it}
\author[inm]{Matteo Diez}
\ead{matteo.diez@cnr.it}

\cortext[cor1]{Corresponding author}
\address[imati]{Consiglio Nazionale delle Ricerche - Istituto di Matematica Applicata e Tecnologie Informatiche ``E. Magenes'' (CNR-IMATI), Via Ferrata 5/A, 27100 Pavia, Italy}
\address[inm]{Consiglio Nazionale delle Ricerche - Istituto di Ingegneria del Mare (CNR-INM), Via di Vallerano 139, 00128, Roma, Italy}

\begin{abstract} 
This paper presents a comparison of two multi-fidelity methods for the forward uncertainty quantification of a naval engineering problem. 
Specifically, we consider the problem of quantifying the uncertainty of the hydrodynamic resistance of a 
roll-on/roll-off passengers ferry advancing in calm water and subject to two operational uncertainties
(ship speed and payload). The first four statistical moments
(mean, variance, skewness, kurtosis), and the probability density function for such quantity
of interest (QoI) are computed with two multi-fidelity methods, i.e., the Multi-Index Stochastic Collocation (MISC) 
and {an} adaptive multi-fidelity Stochastic Radial Basis Functions (SRBF).
The QoI is evaluated via computational fluid dynamics simulations, which are performed with the in-house unsteady Reynolds-Averaged Navier-Stokes (RANS) multi-grid solver $\chi$navis. The different fidelities employed by both methods are obtained by stopping the RANS solver at different grid levels of the multi-grid cycle.
The performance of both methods are presented and discussed: in a nutshell, the findings suggest that, at least for the current implementations of both methods,
MISC could be preferred whenever a limited computational budget is available, whereas for a larger computational budget SRBFs seem to be preferable, thanks to its robustness to the numerical noise in the evaluations of the QoI.
\end{abstract}


\end{frontmatter}

\section{Introduction}
{Aerial, ground, and water-born vehicles must perform, in general, under a variety of environmental and operating conditions and, therefore, their design analysis and  optimization processes cannot avoid  taking into account the stochasticity associated to environmental and operational parameters.
An example is given by ships and their subsystems, which are required to operate under a variety of highly stochastic conditions, such as speed, payload, sea state, and wave heading \cite{serani2021hull}. 
In this context, the accurate prediction of relevant design metrics (i.e., resistance and powering requirements; seakeeping, maneuverability, and dynamic stability; structural response and failure) requires prime-principles-based high-fidelity computational tools (e.g., computational fluid/structural dynamics, CFD/CSD), especially for innovative configurations and off-design conditions. These tools are, however, generally computationally expensive, making the quantification of the relevant statistical indicators (through the use of many function evaluations) a technological challenge. 
As an example, an accurate hull-form optimization based on unsteady Reynolds Averaged Navier-Stokes (URANS) solvers under stochastic conditions may require up to 500k CPU hours on high performance computing (HPC) systems, even if computational cost reduction methods are used \cite{serani2021hull}. Similarly, a URANS-based statistically significant evaluation of ship maneuvering performance in irregular waves may require up to 1M CPU hours on HPC systems \cite{serani2021urans}.
In this context, the use of efficient uncertainty quantification (UQ) methods is essential to make the design analysis and optimization processes affordable. 
The development and application of UQ methods for sea-vehicle problems were discussed in \cite{stern2017development}. Moreover, the numerical UQ of a high-speed catamaran was performed and discussed for calm water \cite{diez2014uncertainty}, regular \cite{he2013urans} and irregular \cite{diez2018statistical} waves conditions. An experimental UQ was presented in \cite{durante2020accurate} for validation purposes of the same model. The efficiency of the UQ methods is, in general, problem dependent and has to be carefully assessed. As an example, several UQ methods were compared for an airfoil benchmark problem in \cite{quagliarella2019benchmarking}.}

In general, there is by now a large consensus in the UQ and computational sciences communities on the fact that large-scale UQ analyses can only be performed by leveraging on multi-fidelity methodologies, i.e., methodologies that explore the bulk of the variability of the quantities of interest (QoI) of the simulation over coarse grids (or more generally, computationally inexpensive models with e.g.\ simplified physics), and resort to querying high-fidelity models (e.g., refined grids or full-physics models) only sparingly, to correct the initial guess produced with the low-fidelity models, see e.g.\ \cite{beran2020}. Within this general framework, several approaches can be conceived, depending on the kind of fidelity models considered and on the strategy used to sample the parameter space (i.e., for what values of the uncertain parameters the different fidelity models should be queried/evaluated).

One large class of methods that has received increasing attention in this context is the family of multi-level/multi-index methods, due to its effectiveness and solid mathematical ground. The hierarchy of models considered by these methods is usually obtained by successive (most often -- but not necessarily -- dyadic) refinements of a computational grid. The multi-level/multi-index distinction arises from the number of {hyper-}parameters that are considered to control the overall discretization of the problem, i.e., how many hyper-parameters are used to determine the computational grids (e.g.\ one or multiple size parameters for the grid elements and/or time-stepping) and the number of samples from the parameters space to be solved on each grid (e.g.\ specified by a single number or by a tuple of different numbers along different directions in the parametric space).

Combining {the above} considerations with a specific sampling strategy over the parameter space results in different variations of the method.
{One approach is to use random/quasi random sampling methods: this leads to methods}
such as Multi-Level Monte Carlo \cite{giles:MLMC,scheichl.giles:MLMC}, Multi-Index Monte Carlo \cite{hajiali.eal:MultiIndexMC}, Multi-Level/Multi-Index Quasi-Monte Carlo \cite{kss12}.

{A different option is to resort to methods that build a polynomial approximation over the parameter space: methods such as }
Multi-Level Stochastic Collocation \cite{teckentrup.etal:MLSC},
Multi-Index Stochastic Collocation \cite{beck.eal:MISC-IGA,jakeman2019adaptive,hajiali.eal:MISC1,hajiali.eal:MISC2}, 
Multi-Level Least-Squares polynomial approximation \cite{hajiali2017multilevel}, etc. {fall in this category. Note that}
{the wording ``Stochastic Collocation'' is to be understood as a synonym of
``sampling in the parametric space'': it refers to the fact that the parameters of the problem can be seen as random (stochastic) variables, and sampling the parametric space can be seen as ``collocating the approximation problem at points of the stochastic domain''.}

Another widely studied class of multi-level methods employs  
kernel-based surrogates such as hierarchical kriging \cite{han2012-AIAA}, co-kriging \cite{debaar2015-CF}, Gaussian processes \cite{wackers2020-AIAA}, and radial-basis functions \cite{serani2019-IJCFD}. Additive, multiplicative, or hybrid correction methods, also known as ``bridge functions'' or ``scaling functions'' \cite{han2013-AST}, are used to build multi-fidelity surrogates. Further efficiency of multi-fidelity surrogates is gained using dynamic/adaptive sampling strategies, for which the multi-fidelity design of experiments for the surrogate training is not defined a-priori but dynamically updated, exploiting the information that becomes available during the training process. Training points are dynamically added with automatic selection of {both} their location and the desired fidelity level, {with} the aim of reducing the computational cost required to properly represent the function \cite{serani2019-IJCFD}.

{Moving away from the multi-level/multi-index paradigm, multi-fidelity methods that are based on different physical models rather than multiple discretizations have been proposed e.g. in \cite{peherstorfer:MFsurvey,Gorodetsky2020b,Pisaroni:MLMC-opt,pisaroni:CMLMC,geraci2019}.}

The objective of the present work is to assess and compare the use of two methods, one from each methodological family, for the forward UQ analysis 
of a naval engineering problem. Specifically, the performance of the 
Multi-Index Stochastic Collocation (MISC)  
and adaptive multi-fidelity Stochastic Radial Basis Functions (SRBF) methods 
is compared {on: i) an analytical test function and ii) the forward UQ analysis of a roll-on/roll-off passengers (RoPax) ferry sailing in calm water with two operational uncertainties, specifically ship speed and draught, the latter being directly linked to the payload. The estimation of the expected value, variance, skewness, kurtosis, and of the probability density function (PDF) of the function and the hydrodynamic resistance of the RoPax, is presented and discussed.
The test function considered in the analytical test is tailored to resemble the surrogate model of the naval engineering problem: the results of this preliminary test can then be considered as a baseline for the assessment of the relative performances of the two methods, 
and help in interpreting the results of the naval test case.}
{In the RoPax problem the} hydrodynamic resistance for each value of speed and draught requested by MISC and SRBF is computed by the URANS equation solver $\chi$navis \cite{dimascio2007,dimascio2009,broglia2018}, developed at CNR-INM. $\chi$navis {embeds a multi-grid approach for iterations acceleration, based on a sequence of grids obtained by} derefining an initial fine grid. More specifically, in this work four grids are used, and leveraged by both MISC and SRBF to vary the fidelity of the simulations.
Therefore, both MISC and SRBF are used as multi-index methods with only one component controlling the spatial discretization. 
Another relevant aspect is that $\chi$navis is an iterative solver, and, as such, it stops as soon as a suitable norm of the residual drops below a prescribed tolerance. The fact that the RANS equations are not solved at machine-precision introduces in practice some noise in the evaluation of the resistance, which needs to be dealt with during the computations of the UQ indicators (statistical moments, PDF, etc). 
 
A preliminary version of this work is available as proceedings of the AIAA Aviation 2020 Forum, see \cite{piazzola2021:ferry-AIAA}. With respect to that version, the manuscript was significantly improved in many ways. First, the discussion on MISC is now focused on the construction of the surrogate model rather than on computing statistical moments, and in particular we added some (we believe) interesting considerations about the fact that the MISC surrogate model is \emph{not} interpolatory, even when using nested points in the parametric space; to the best of the authors' knowledge, this fact was never mentioned in previous literature. {Second, for the SRBF method applied to the RoPax UQ analysis, a methodological advancement is used. In the previous work 
 interpolation was enforced for the early iterations. Then, when a certain number of training points was available, an optimization process was performed to automatically select the number and the position of centers of the SRBF, thus automatically selecting whether to perform interpolation or regression. Differently, in this work the optimization process is performed since the first iteration, thus making the methodology fully adaptive.} 
Finally, 
the numerical results section has been enriched by including an analytical test, by adding a reference solution for the naval problem (which is obtained by a sparse-grid sampling of the highest fidelity at our disposal), and by  discussing a possible strategy to mitigate the impact of the RANS noise on the MISC framework.

The remainder of this paper is organized as follows. Section \ref{sect:forward-UQ} introduces the general framework and notation for the UQ problem, and the two methodologies considered in this work; in particular, MISC is presented in Section \ref{sect:misc}, while SRBFs are presented in Section \ref{sect:SRBF}. 
{Section \ref{sect:numerical_tests} presents the numerical results: a preliminary analytical test (see Section \ref{sect:problem_description_analytic}) and then the naval problem (see Section \ref{sect:problem_description_ropax}). }
Finally, a summary of the findings of the numerical tests and an outlook on future work is presented in Section \ref{sect:conclusions}.

\section{Forward uncertainty quantification methods}\label{sect:forward-UQ}
Let us assume that we are interested in the outcome of a CFD simulation that depends on the value of $N$ random/uncertain parameters collected in the vector ${\bf y}=[y_1,y_2,\ldots,y_N]$; we denote by $\Gamma \subseteq \mathbb{R}^N$ the set of all possible values of ${\bf y}$, and by $\rho({\bf y})$ the PDF of ${\bf y}$ over $\Gamma$. The goal of a forward UQ analysis is to compute statistical indicators of the QoI, $G$, of such CFD simulation,   to quantify its variability due to the uncertainties on ${\bf y}$.  For instance, we might be interested in computing expected values and/or higher-order moments of $G$  (in the numerical tests we will report on mean, variance, skewness, and kurtosis, denoted by  $\mathbb{E}[G]$, $\text{Var}[G]$, $\text{Skew}[G]$, and $\text{Kurt}[G]$, respectively),  and the PDF of $G$, which completely describes its statistical variability.

This analysis is often performed by a sampling approach, i.e., the CFD simulation is run 
for several possible values of ${\bf y}$, and the corresponding results are post-processed to get the indicators of interest.  For instance, the statistical moments can be approximated by weighted averages of the values obtained,
while the PDF can be approximated by histograms or, e.g.,\ kernel density methods \cite{rosenblatt:kde,parzen:kde}. 
Clearly, these analyses require large datasets of evaluations of $G$: if computing a single instance of $G$  requires a significant amount of computational time, obtaining the dataset can become prohibitively expensive.  A possible workaround is then to replace the evaluations of $G$ with the evaluations of a surrogate model,  which is ideally a good approximation of the original $G$, cheap to evaluate and obtained by suitably combining together a relative small number of evaluations of $G$  (less than what would be needed to perform the UQ analysis of the full model).  The two methods that we consider in this work are both methods to construct such surrogate model, and in particular they leverage the fact that CFD simulations can be performed over multiple grid resolutions to further reduce the computational costs. 

Before describing in detail each method, we need to introduce some notation. To this end, let us assume that the computational domain of our CFD simulation can be discretized by a grid with non-cubic hexahedral elements of the same size\footnote{The assumption that all elements must be of the same size can be relaxed, but it is kept for simplicity of exposition.} and let us also assume for a moment that the level of refinement of the grid along each physical direction can be specified by prescribing some integer values $\alpha_1, \alpha_2,\alpha_3$; to fix ideas, one can think, e.g.,\ that the number of elements of the grid scales as $2^{\alpha_1} \times 2^{\alpha_2} \times 2^{\alpha_3}$, but this is not necessary. 
The three values of $\alpha_i$ are collected in a multi-index $\aalpha = [\alpha_1, \alpha_2, \alpha_3]$; prescribing the multi-index $\aalpha$ thus prescribes the computational grid to be generated. If this flexibility is not allowed by the grid-generator (or by the problem itself), it is possible to set $\alpha_1 = \alpha_2 = \alpha_3=\alpha$, i.e., controlling the grid-generation by a single integer value $\alpha$ {(this is actually the case for the RoPax ferry example considered in this work)}. 
The same  philosophy applies both to single- and multi-patch grids, where in principle there could be up to three values $\alpha_i$ for each patch.
In general, we assume that $\aalpha$ has $d$ components, $\aalpha \in \mathbb{N}^d_+$.
The QoI of the CFD simulation computed over the grid specified by $\aalpha$ is denoted by $G_{\aalpha}$; this could be, e.g.,\ the full velocity field or a scalar quantity associated to it.

\subsection{Multi-Index Stochastic Collocation (MISC)} \label{sect:misc}

In this section, the MISC method is introduced.
As already mentioned, the MISC method is a multi-fidelity method that falls under the umbrella of multi-index/multi-level methods:
in particular, the single-fidelity models upon which MISC is built are global Lagrangian interpolants over $\Gamma$.   

\subsubsection{Tensorized Lagrangian interpolant operators}

The first step to derive the MISC surrogate model is to select a sequence of collocation points for each uncertain parameter $y_n$,
i.e., for each direction of $\Gamma_n$ of $\Gamma$.
For computational efficiency, these points 
should be chosen according to $\rho(\yy)$, and they should be of nested type (i.e., collocation grids of increasing refinement
should be subset of one another). 
In the RoPax ferry example considered in this work, the uncertain parameters $\yy$ can be modeled as uniform and independent {random variables} (see Sect.~\ref{sect:problem_description_ropax}) 
for which we choose to employ Clenshaw--Curtis (CC) points, see e.g. \cite{trefethen:comparison}. A set of $K$ univariate CC points can be obtained as 
\begin{equation}\label{eq:CC_points}
t_{K}^{(j)} =\cos\left(\frac{(j-1) \pi}{K-1}\right), \quad 1\leq j \leq K,
\end{equation}
and two sets of CC points, with $K_1$ and $K_2$ points, are nested if $(K_2  - 1)/(K_1 - 1)= 2^\ell$ for some integer $\ell$, see also below.
Other nested alternatives for uniformly distributed {parameters} are Leja points \cite{narayan:Leja,nobile.etal:leja} and Gauss--Patterson points \cite{patterson:addition}.
Next, we introduce the function 
\begin{equation}\label{eq:lev_fun}
m(0)=0, \ m(1)=1, \ m(\beta_n) = 2^{\beta_n-1}+1 \ \text{for} \ \beta_n \geq 2, 
\end{equation} 
and denote by $\mathcal{T}_{n,\beta_n}$ the set of $m(\beta_n)$ CC points along $y_n$, i.e.
\[
  \mathcal{T}_{n,\beta_n} = \left\{y_{n,m(\beta_n)}^{(j_n)} \biggr\rvert j_n=1, \ldots, m(\beta_n)\right\} \quad \text{ for } n=1,\ldots,N.
\]
Note that this choice of $m$ guarantees nestedness of two sets of CC points, i.e.,
$\mathcal{T}_{n,\beta} \subseteq \mathcal{T}_{n,\gamma}$ if $\gamma \geq \beta$.

An $N$-dimensional interpolation grid can then be obtained by taking the Cartesian product of the $N$ univariate sets just introduced.
The number of collocation points in this grid is specified by a multi-index $\bbeta \in \mathbb{N}^N_+$: 
such multi-index plays thus a similar role for the parametric domain $\Gamma$ as the multi-index $\aalpha$ for the physical domain.
We denote such tensor interpolation grid by $\mathcal{T}_{\bbeta} = \bigotimes_{n=1}^{N} \mathcal{T}_{n,\beta_n}$  
and its number of points by $M_{\bbeta} = \prod_{n=1}^{N} m(\beta_n)$: using standard multi-index notation,
they can be written as
\[
  \mathcal{T}_{\bbeta} = \left\{\yy_{m(\bbeta)}^{(\jj)}\right\}_{\jj \leq m(\bbeta)},
  \quad  \text{ with } \quad
  \yy_{m(\bbeta)}^{(\jj)} = \left[y_{1,m(\beta_1)}^{(j_1)}, \ldots, y_{N,m(\beta_N)}^{(j_N)}\right],
\]
where $m(\bbeta) = \left[m(\beta_1),\,m(\beta_2),\ldots,m(\beta_N) \right]$ and $\jj \leq m(\bbeta)$ means that $j_n \leq m(\beta_n)$ for every $n = 1,\ldots,N$.
For fixed $\aalpha$, the approximation of $G_{\aalpha}(\yy)$ based on global Lagrangian polynomials collocated at these grid points
(single-fidelity approximation) has the following form  
\begin{equation} \label{eq:MISC-singlefidelity}
G_{\aalpha}(\yy) \approx \mathcal{U}_{\aalpha,\bbeta}(\yy) := \sum_{\jj \leq m(\bbeta)} G_{\aalpha}\left(\yy_{m(\bbeta)}^{(\jj)}\right) \mathcal{L}_{m(\bbeta)}^{(\jj)}(\yy),
\end{equation}
where $\left\{ \mathcal{L}_{m(\bbeta)}^{(\jj)}(\yy) \right\}_{\jj \leq m(\bbeta)}$ are $N$-variate Lagrange basis polynomials, 
defined as tensor products of univariate Lagrange polynomials, i.e.  
\begin{equation}\label{eq:lagrange_basis}
  \mathcal{L}_{m(\bbeta)}^{(\jj)}(\yy) = \prod_{n=1}^{N} l_{n,m(\beta_n)}^{(j_n)}(y_n)
  \quad \text{ with } \quad
  l_{n,m(\beta_n)}^{(j_n)}(y_n) = \prod_{k=1, k\neq j_n}^{m(\beta_n)} \frac{y_n-y_{n,m(\beta_n)}^{(k)}}{y_{n,m(\beta_n)}^{(k)}-y_{n,m(\beta_n)}^{(j_n)}}. 
\end{equation}

Naturally, the single-fidelity approximation $\mathcal{U}_{\aalpha,\bbeta}$ is more and more accurate the higher the number of collocations points in each direction.
Hence, ideally one would choose both multi-indices $\aalpha$ and $\bbeta$ with large components,
say $\aalpha = \aalpha^\star$ and $\bbeta = \bbeta^\star$, i.e., to consider many CFD simulations over a refined computational grid;
however, this is typically infeasible due to the computational cost of a single CFD simulation.

\subsubsection{MISC surrogate model}

The above discussion on the costs of $\mathcal{U}_{\aalpha^\star,\bbeta^\star}$ motivates the introduction of MISC.
MISC is a multi-fidelity approximation method that replaces $\mathcal{U}_{\aalpha^\star,\bbeta^\star}$
with a linear combination of multiple coarser $\mathcal{U}_{\aalpha,\bbeta}$:
as will be clearer later, the components of such linear combination are chosen
obeying to the idea that whenever the spatial discretization $\aalpha$ is refined,
the order of the interpolation $\bbeta$ is kept to a minimum and vice versa.

\begin{figure}[t]
	\centering
	\subfigure[]{\includegraphics[width=0.3\linewidth]{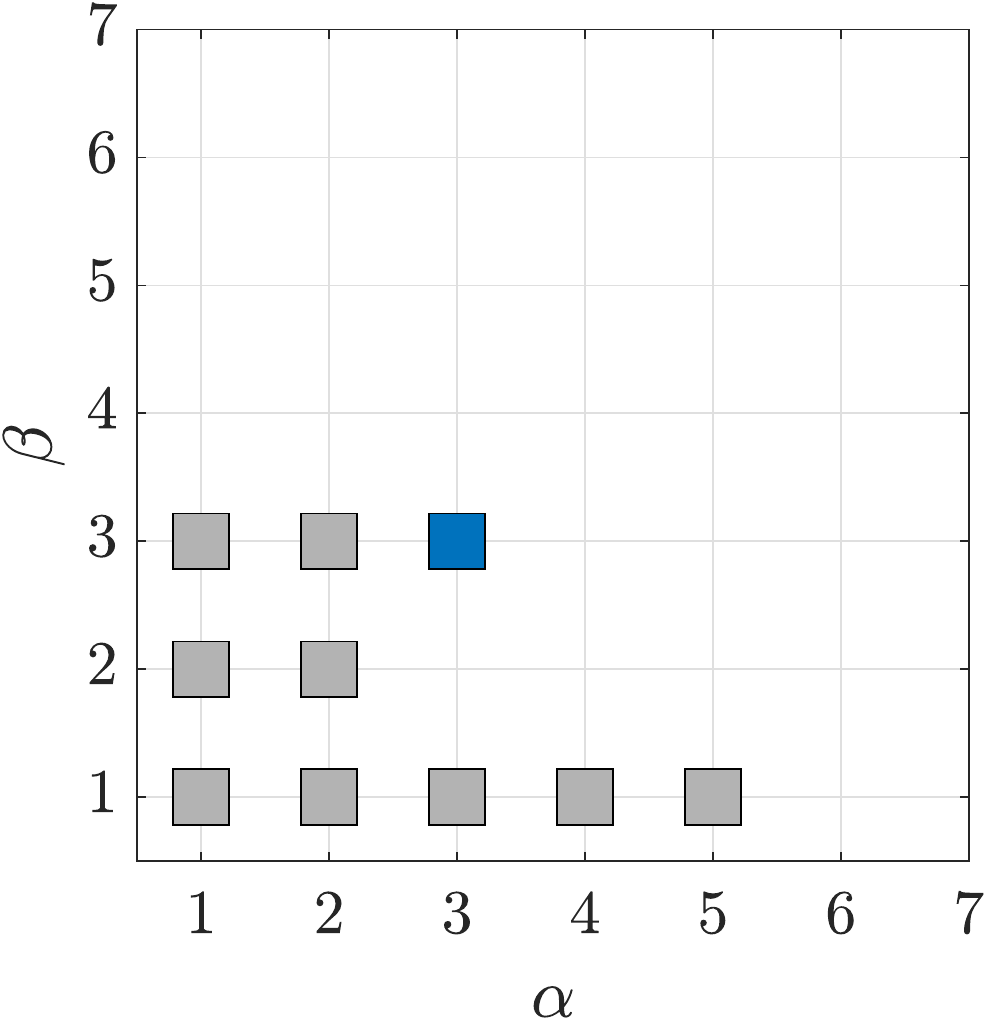}} \hspace{1em}
	\subfigure[]{\includegraphics[width=0.3\linewidth]{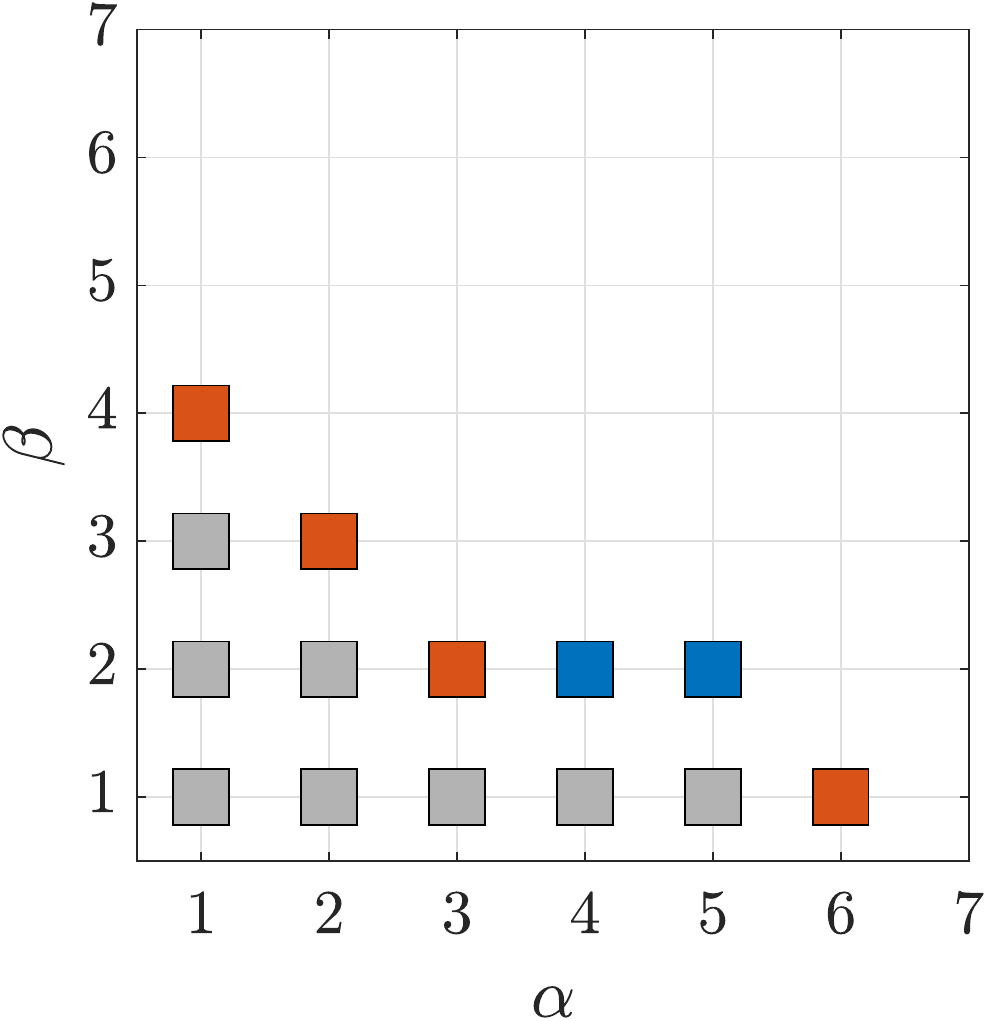}}
	\caption[index_sets]{Multi-index sets for the construction of the MISC approximation {(in the case $d=N=1$)}. (a): the gray set is downward closed, whereas adding the blue multi-index to it would result in a set not downward closed; (b): a downward closed set (in gray) and its margin (indices marked in red and blue). If Algorithm \ref{algo:misc_implementation} reaches the gray set, it will next explore all indices marked in red (their addition to the gray set keeps the downward closedness property) but not those marked in blue. The red set is also known as ``reduced margin''. }
	\label{fig:index_sets}
\end{figure}

To build a MISC approximation, the so-called ``detail operators'' (univariate and multivariate) on the physical and parametric domains have to be introduced.
They are defined as follows, with the understanding that $\mathcal{U}_{\aalpha,\bbeta} (\yy) =0$ when at least one component of $\aalpha$ or $\bbeta$ is zero.
In the following the dependence of the interpolation operator on the parameters $\yy$ is omitted for sake of compactness.  
Thus, we denote by $\ee_i$ the canonical multi-index, i.e.\ $(\ee_i)_k = 1$ if $i=k$ and 0 otherwise, and define
\begin{alignat*}{2}
&\text{\textbf{Univariate physical detail: }}
&& \Delta_i^{\text{phys}}[\mathcal{U}_{\aalpha,\bbeta}]=\mathcal{U}_{\aalpha,\bbeta}-\mathcal{U}_{\aalpha-\ee_i,\bbeta} \text{ with } 1 \leq i \leq d; \\
&\text{\textbf{Univariate parametric detail: }}
&& \Delta_i^{\text{param}}[\mathcal{U}_{\aalpha,\bbeta}]=\mathcal{U}_{\aalpha,\bbeta}-\mathcal{U}_{\aalpha,\bbeta-\ee_i} \text{ with } 1 \leq i \leq N; \\
&\text{\textbf{Multivariate physical detail: }}
&& \bm{\Delta}^{\text{phys}}[\mathcal{U}_{\aalpha,\bbeta}] = \bigotimes_{i=1}^d \Delta_i^{\text{phys}}[\mathcal{U}_{\aalpha,\bbeta}]; \\
&\text{\textbf{Multivariate parametric detail: }}
&& \bm{\Delta}^{\text{param}}[\mathcal{U}_{\aalpha,\bbeta}] = \bigotimes_{j=1}^N \Delta_j^{\text{param}}[\mathcal{U}_{\aalpha,\bbeta}]; \\
&\text{\textbf{Mixed multivariate detail: }}
&& \bm{\Delta}^{\text{mix}}[\mathcal{U}_{\aalpha,\bbeta}] = \bm{\Delta}^{\text{param}}\left[\bm{\Delta}^{\text{phys}}[\mathcal{U}_{\aalpha,\bbeta}] \right].  
\end{alignat*}
Observe that taking tensor products of univariate details amounts to composing their actions, i.e. 
\[
\bm{\Delta}^{\text{phys}}[\mathcal{U}_{\aalpha,\bbeta}]
= \bigotimes_{i=1}^d \Delta_i^{\text{phys}}[\mathcal{U}_{\aalpha,\bbeta}]
= \Delta_1^{\text{phys}}\left[ \, \cdots \left[ \Delta_d^{\text{phys}}\left[ \mathcal{U}_{\aalpha,\bbeta} \right] \, \right] \, \right],
\]
and analogously for the multivariate parametric detail operators $\bm{\Delta}^{\text{param}}[\mathcal{U}_{\aalpha,\bbeta}]$.
By replacing the univariate details with their definitions, we can then see that this implies
that the multivariate operators can be evaluated by evaluating certain full-tensor approximations  
$\mathcal{U}_{\aalpha,\bbeta}$ introduced in the previous subsection, and then taking linear combinations:
\begin{align*}
\bm{\Delta}^{\text{phys}}[\mathcal{U}_{\aalpha,\bbeta}]
& = \Delta_1^{\text{phys}}\left[ \, \cdots \left[ \Delta_d^{\text{phys}}\left[ \mathcal{U}_{\aalpha,\bbeta} \right] \, \right] \, \right] 
= \sum_{\jj \in \{0,1\}^d} (-1)^{\lVert\jj\rVert_1} \mathcal{U}_{\aalpha-\jj,\bbeta};\\
\bm{\Delta}^{\text{param}}[\mathcal{U}_{\aalpha,\bbeta}]
&  = \sum_{\jj \in \{0,1\}^N} (-1)^{\lVert\jj\rVert_1} \mathcal{U}_{\aalpha,\bbeta-\jj}.
\end{align*}
The latter expressions are known as ``combination-technique'' formulations, and can be very
useful for practical implementations. In particular, they allow to evaluate e.g.\ $\bm{\Delta}^{\text{phys}}[\mathcal{U}_{\aalpha,\bbeta}]$
by calling pre-existing softwares on different grids up to $2^d$ times in a ``black-box'' fashion.
Analogously, evaluating  $\bm{\Delta}^{\text{param}}[\mathcal{U}_{\aalpha,\bbeta}]$
requires evaluating up to $2^N$ operators $\mathcal{U}_{\aalpha,\bbeta}$ over different interpolation grids,
and evaluating $\bm{\Delta}^{\text{mix}}[\mathcal{U}_{\aalpha,\bbeta}]$ requires evaluating
up to $2^{d+N}$ operators $\mathcal{U}_{\aalpha,\bbeta}$ over different parametric grids and physical grids.
Observe that by introducing these detail operators a hierarchical decomposition of $\mathcal{U}_{\aalpha,\bbeta}$ can be obtained;
indeed, the following telescopic identity holds true:
\begin{equation}\label{eq:telescopic_sum}
\mathcal{U}_{\aalpha,\bbeta} = \sum_{[\ii,\jj] \leq [\aalpha,\bbeta]} \bm{\Delta}^{\text{mix}}[\mathcal{U}_{\ii,\jj}].
\end{equation}
As an example, the case of $d=N=1$ (i.e., one-dimensional physical and parametric spaces) can be considered.
Recalling that by definition $\mathcal{U}_{i,j} = 0$ when either $i=0$ or $j=0$, it can be seen that
\begin{align*}
\sum_{[i,j] \leq [2,2]} \bm{\Delta}^{\text{mix}}[\mathcal{U}_{i,j}] 
& = \bm{\Delta}^{\text{mix}}[\mathcal{U}_{1,1}]
+ \bm{\Delta}^{\text{mix}}[\mathcal{U}_{1,2}]
+ \bm{\Delta}^{\text{mix}}[\mathcal{U}_{2,1}]
+ \bm{\Delta}^{\text{mix}}[\mathcal{U}_{2,2}] \\
& = \mathcal{U}_{1,1}
+ ( \mathcal{U}_{1,2} - \mathcal{U}_{1,1} )
+ ( \mathcal{U}_{2,1} - \mathcal{U}_{1,1} )
+ (\mathcal{U}_{2,2} - \mathcal{U}_{2,1} - \mathcal{U}_{1,2} + \mathcal{U}_{1,1}) \nonumber \\
& = \mathcal{U}_{2,2}. \nonumber
\end{align*}
The crucial observation is that, under suitable regularity assumptions for $G(\yy)$, see e.g.\ \cite{hajiali.eal:MISC1,hajiali.eal:MISC2}),
not all of the details in the hierarchical decomposition in Eq.~\eqref{eq:telescopic_sum} contribute equally to the approximation,
i.e., some of them can be discarded and the resulting formula will retain good approximation properties at a fraction of the computational cost
(roughly, the multi-indices to be discarded are those corresponding to ``high-order'' details, i.e., those for which $\| \aalpha \|_1 + \| \bbeta \|_1$
is sufficiently large). Upon collecting the multi-indices $[\aalpha,\bbeta]$ to be retained in the sum in a multi-index set $\Lambda \subset \mathbb{N}_+^{d+N}$,
the MISC multi-fidelity approximation of $G$ can be introduced as

\begin{equation}\label{eq:misc_surr_delta}
G(\yy) \approx \MISCsurr_\Lambda (\yy) :=\sum_{[\aalpha,\bbeta] \in \Lambda} \bm{\Delta}^{\text{mix}}[\mathcal{U}_{\aalpha,\bbeta}(\yy)].
\end{equation}
To obtain a meaningful expression, $\Lambda$ should be chosen as downward closed, i.e.\ (see Fig.~\ref{fig:index_sets}a)
\[
\forall \mathbf{k} \in \Lambda, \quad \mathbf{k} - \ee_j \in \Lambda \text{ for every } j=1,\ldots,d+N \text{ such that } k_j > 1.  
\]

Clearly, the MISC formula in Eq.~\eqref{eq:misc_surr_delta} has a combination-technique expression as well, which can be written in compact form as
\begin{equation}\label{eq:misc}
\MISCsurr_\Lambda (\yy)
=\sum_{[\aalpha,\bbeta] \in \Lambda} \bm{\Delta}^{\text{mix}}[\mathcal{U}_{\aalpha,\bbeta}(\yy)]
=\sum_{[\aalpha,\bbeta] \in \Lambda} c_{\aalpha,\bbeta} \, \mathcal{U}_{\aalpha,\bbeta}(\yy), 
\end{equation}
where the coefficients $c_{\aalpha,\bbeta}$ are defined as 
\[
c_{\aalpha,\bbeta} = \sum_{\substack{[\ii,\jj] \in \{0,1\}^{d+N}\\ [\aalpha+\ii,\bbeta+\jj] \in \Lambda}}
(-1)^{\lVert [\ii,\jj] \rVert_1}.
\]

This is the approximation formula which is used in {our} practical implementation of the MISC method,
which shows our initial statement that the MISC evaluation is computed by evaluating full-tensor interpolation operators
$\mathcal{U}_{\aalpha,\bbeta}$ independently and combining them linearly, as specified by Eq.~\eqref{eq:misc}.
Before going further, we remark a few important points:
\begin{enumerate}
\item The effectiveness of the MISC approximation depends on the choice of the multi-index set $\Lambda$.
    The optimal choice of $\Lambda$ depends on the regularity assumptions on $G(\yy)$:
    in general, the result will be a method akin to classical multi-level schemes
    such as Multi-Level Monte Carlo, where most of the statistical variability of the QoI is explored by solving many CFD simulations with coarse grids (large $\|\bbeta\|_1$ with small $\|\aalpha\|_1$)
    and then the result is corrected with a few CFD simulations with refined grids (large {$\|\aalpha\|_1$} with small $\|\bbeta\|_1$).
    A practical adaptive algorithm to construct $\Lambda$ is presented in Sect.~\ref{sect:adaptive-MISC}, following the
    discussion in, e.g.,\ \cite{jakeman2019adaptive}.
    Another option is to design $\Lambda$ a-priori, by a careful analysis of the PDE at hand, see, e.g.,\ \cite{beck.eal:MISC-IGA,hajiali.eal:MISC1,hajiali.eal:MISC2}.

  \item MISC works well only if the levels are sufficiently separated,
    i.e.\ if the number of degrees of freedom of the computational grid (and the corresponding computational cost)
    grows significantly from one level to the next one: to fix ideas, one such case is if
    the number of elements in the grid scales, e.g.,\ as $2^{\alpha_1} \times 2^{\alpha_2} \times 2^{\alpha_3}$,
    but not if, e.g.,\ increasing $\alpha_1$ to $\alpha_1+1$ adds only one element to the grid.
    If this separation does not hold, the cost of computing all the components in Eq.~\eqref{eq:misc} would
    exceed the cost of the construction of a highly-refined single-fidelity surrogate model $\mathcal{U}_{\aalpha^\star,\bbeta^\star}$.

  \item The MISC surrogate model $\MISCsurr_\Lambda (\yy)$ is not interpolatory, even when nested nodes
    are used in the parametric space (as it is the case here).
    To illustrate this, let us consider as an example the case $d=1, N=2$, with $\Gamma = [-1,1]^2$.
    We construct the MISC approximation based on the multi-index set $\Lambda = \{[1, 1\, 1],\,  [1, 2 \,  1],\, [2, 1 \, 1]\}$:
    formula \eqref{eq:misc} results in 
    \[
      \MISCsurr_\Lambda (\yy) = -\mathcal{U}_{[1, 1\, 1]}(\yy)+ \mathcal{U}_{[2, 1 \,1 ]}(\yy) + \mathcal{U}_{[1, 2 \, 1]}(\yy),
    \]
    where the first and the second operators are constant interpolants ($\beta = [1 \, 1]$) whose value is equal to the first and second fidelity
    evaluated at the center of the parametric domain, $\yy^{(C)}=[0,0]$, respectively; the third operator is instead is an interpolant of degree two
    based on the value of the first fidelity evaluated at the following three CC bivariate points: $\yy^{(L)}=[-1,0]$, $\yy^{(C)}=[0,0]$,
    and $\yy^{(R)}=[1,0]$. Then, evaluating the MISC approximation at, e.g.,\ $\yy^{(R)}$ results in   
    \[
    \MISCsurr_\Lambda (\yy^{(R)}) =  - G_1(\yy^{(C)}) + G_2(\yy^{(C)}) +  G_1(\yy^{(R)})  \neq G_1(\yy^{(R)}), 
    \]
    i.e.\ the value of $\MISCsurr_{\Lambda}$ at $\yy^{(R)}$ is different from the only model evaluation available at such point (which is $G_1(\yy^{(R)})$).
	This is in contrast with the well-known property of single-fidelity sparse-grid surrogate models,
	for which the use of nested nodes over $\Gamma$ guarantees that $\MISCsurr_\Lambda(\yy_{SG}) = G(\yy_{SG})$ at the sparse-grid points $\yy_{SG}$.
  \end{enumerate}

\subsubsection{MISC quadrature}

By taking the integral of the MISC surrogate model defined in Eq.~\eqref{eq:misc} it is straightforward to obtain a quadrature formula $\MISCquad_\Lambda $
to approximate the expected value of $G$:
\[
\mathbb{E}[G] \approx \mathbb{E}[\MISCsurr_{\Lambda}] =: \MISCquad_\Lambda = \int_{\Gamma} \MISCsurr_\Lambda (\yy) \rho(\yy) \dd \yy
=\sum_{[\aalpha,\bbeta] \in \Lambda} c_{\aalpha,\bbeta} \int_{\Gamma} \mathcal{U}_{\aalpha,\bbeta}(\yy) \rho(\yy) \dd \yy. 
\]
By recalling the definition of $\mathcal{U}_{\aalpha,\bbeta}$ given in Eq.~\eqref{eq:MISC-singlefidelity}
and of the multivariate Lagrange polynomials in Eq.~\eqref{eq:lagrange_basis},
each of the integrals at the right-hand side of the previous formula can be rewritten in compact form as
tensor quadrature operators, i.e.  
\begin{align*}
  \mathcal{Q}_{\aalpha,\bbeta} := \int_{\Gamma} \mathcal{U}_{\aalpha,\bbeta}(\yy) \rho(\yy) \dd \yy 
  & = \sum_{\jj \leq m(\bbeta)} G_{\aalpha}\left(\yy_{m(\bbeta)}^{(\jj)}\right) \left( \prod_{n=1}^{N} \int_{\Gamma_n} l_{n,m(\beta_n)}^{(j_n)}(y_n) \rho(y_n) \dd y_n\right) \\
  & = \sum_{\jj \leq m(\bbeta)} G_{\aalpha}\left(\yy_{m(\bbeta)}^{(\jj)}\right)  \left( \prod_{i=1}^{N} \omega_{n,m(\beta_n)}^{(j_n)} \right)
    = \sum_{\jj \leq m(\bbeta)} G_{\aalpha}\left(\yy_{m(\bbeta)}^{(\jj)}\right) \omega_{m(\bbeta)}^{(\jj)},
\end{align*}
where $\omega_{n,m(\beta_n)}^{(j_n)}$, are the standard quadrature weights obtained by computing the integrals of the associated univariate
Lagrange polynomials (available as analytical or tabulated values for most families of collocation points),
and $\omega_{m(\bbeta)}^{(\jj)}$ are their multivariate counterparts. 
The quadrature formula $\MISCquad_\Lambda$ can then be understood as a linear combination of tensor quadrature operators, in complete analogy
with the MISC surrogate model construction:
\begin{equation}\label{eq:misc_quad_exp}
\MISCquad_\Lambda =  \sum_{[\aalpha,\bbeta] \in \Lambda} c_{\aalpha,\bbeta} \mathcal{Q}_{\aalpha,\bbeta}. 
\end{equation}
Equivalently, one can also write 
\[
  \MISCquad_\Lambda =  \sum_{[\aalpha,\bbeta] \in \Lambda} \bm{\Delta}^{\text{mix}}[\mathcal{Q}_{\aalpha,\bbeta}],  
\]
where the definition of $\bm{\Delta}^{\text{mix}}[\mathcal{Q}_{\aalpha,\bbeta}]$ can be easily deduced by replacing the interpolation operators with quadrature operators in the definition of the detail operators given earlier in this section.   
Clearly, formula \eqref{eq:misc_quad_exp} easily generalizes to the computation of higher-order moments:
\begin{equation} \label{eq:MISC_quad}
  \mathbb{E}[G^r] 
  \approx   \sum_{[\aalpha,\bbeta] \in \Lambda}  c_{\aalpha,\bbeta} \sum_{\jj \leq m(\bbeta)} G^r_{\aalpha}\left(\yy_{m(\bbeta)}^{(\jj)}\right) \omega_{m(\bbeta)}^{(\jj)}, \text{ with } r\geq 1.
\end{equation}
However, this formula might not be the most effective approach to approximate $\mathbb{E}[G^r]$, especially in case of noisy evaluations of $G$;
this aspect is discussed in more detail in Sect.~\ref{sect:RoPax_results}.

We close the discussion on the MISC quadrature by remarking that the collocation points have a twofold use, i.e., they are both interpolation and quadrature points.
This aspect significantly differentiates the MISC method from the approach based on radial basis functions presented in Sect.~\ref{sect:SRBF}, where two distinct sets of points are considered: one for constructing a surrogate models (``training points''), and one for obtaining sample values of the surrogate models and deriving an estimate of expected value and higher moments of $G$ (``quadrature points'').

\subsubsection{An adaptive algorithm for the multi-index set $\Lambda$} \label{sect:adaptive-MISC}
As already mentioned, the effectiveness of the MISC approximation depends on the choice of the multi-index set $\Lambda$:
in this work such set is built with an adaptive algorithm, see \cite{jakeman2019adaptive}. 
We begin by introducing the following decomposition of the quadrature error
\begin{align*}
\lvert \mathbb{E}[G] - \MISCquad_\Lambda \rvert
&  = \Big \lvert \mathbb{E}[G] - \sum_{[\aalpha,\bbeta] \in \Lambda} \bm{\Delta}^{\text{mix}}[\mathcal{Q}_{\aalpha,\bbeta}] \Big\rvert \nonumber \\
&  = \Big \lvert \sum_{[\aalpha,\bbeta] \not \in \Lambda} \bm{\Delta}^{\text{mix}}[\mathcal{Q}_{\aalpha,\bbeta}] \Big \rvert
\leq \sum_{[\aalpha,\bbeta] \not \in \Lambda} \big \lvert \bm{\Delta}^{\text{mix}}[\mathcal{Q}_{\aalpha,\bbeta}] \big \rvert
= \sum_{[\aalpha,\bbeta] \not \in \Lambda} \mathcal{E}^{\MISCquad}_{\aalpha,\bbeta},  
\end{align*}
where $\mathcal{E}^{\MISCquad}_{\aalpha,\bbeta} := \big \lvert \bm{\Delta}^{\text{mix}}[\mathcal{Q}_{\aalpha,\bbeta}] \big \rvert$.
$\mathcal{E}^{\MISCquad}_{\aalpha,\bbeta}$ thus represents the ``error contribution'' of $[\aalpha,\bbeta]$,
i.e., the reduction in the quadrature error due to having added $[\aalpha,\bbeta]$ to the current
index-set $\Lambda$. In practice, $\mathcal{E}^{\MISCquad}_{\aalpha,\bbeta}$ can be conveniently computed by 
\begin{equation}\label{eq:error_contr_quad}
\mathcal{E}^{\MISCquad}_{\aalpha,\bbeta} = \lvert \MISCquad_{\Lambda \cup [\aalpha,\bbeta]} -\MISCquad_\Lambda \rvert;
\end{equation}
for any $\Lambda$ downward-closed set such that $\Lambda \cup [\aalpha,\bbeta]$ is also downward-closed. 
A similar quadrature-based error contribution is considered in \cite{jakeman2019adaptive},
where a convex combination of the error in the computation of the mean and of the variance of the QoI is used. 
Another possibility is to introduce an error decomposition based on the point-wise accuracy of the surrogate model,
following the same arguments above. The ``error contribution'' of $[\aalpha,\bbeta]$ is then taken as 
\begin{equation}\label{eq:error_contr_Linf}
  \mathcal{E}^{\MISCsurr}_{\aalpha,\bbeta} = \lVert \MISCsurr_{\Lambda \cup [\aalpha,\bbeta]} -\MISCsurr_{\Lambda} \rVert_{L_{\infty}}
   \approx \max_{\yy \in \mathcal{H}}  \left\lvert \MISCsurr_{\Lambda \cup [\aalpha,\bbeta]} (\yy) -\MISCsurr_\Lambda(\yy) \right\rvert,
\end{equation}
where $\mathcal{H} \subset \Gamma$ is a suitable set of ``testing points''.
Note that a similar criterion has been proposed also in the context of sparse-grid methods:
different choices of $\mathcal{H}$ can be considered, depending whether nested or non-nested points are used
  (cf., e.g.,\ \cite{chkifa:adaptive-interp,schillings.schwab:inverse} and \cite{nobile.eal:adaptive-lognormal}, respectively).
  In this work we consider a set of 10000 random points (note that this operation is not expensive since it does not require evaluations of the full model).

Similarly to the ``error contribution'', the ``work contribution'' $\mathcal{W}_{\aalpha,\bbeta}$ of $[\aalpha,\bbeta]$ is defined as the work required
to add $[\aalpha,\bbeta]$ to the current index-set $\Lambda$.
It is the product of the computational cost associated to the spatial grid identified by the multi-index $\aalpha$,
denoted by $\text{cost}(\aalpha)$ (see details in Sect.~\ref{sect:RoPax_problem}, Eq.~\eqref{eq:cost}),
times the number of new evaluations of the PDE required by the multi-index $\bbeta$, i.e.  
\begin{equation}\label{eq:work_contr}
\mathcal{W}_{\aalpha,\bbeta} = \text{cost}(\aalpha)  \prod_{n=1}^N (m(\beta_n)-m(\beta_n-1)),
\end{equation}
with $m$ defined as in Eq.~\eqref{eq:lev_fun}. Note that the expression above is based on the fact that the
collocation points used here are nested.

{We then introduce the so-called ``profit'' associated to the multi-index $\aalpha,\bbeta$, which is defined in correspondence with the two choices of error contribution above as 
\begin{equation}\label{eq:profits}
	P^{\MISCsurr}_{\aalpha,\bbeta} = \frac{\mathcal{E}^{\MISCsurr}_{\aalpha,\bbeta}}{\mathcal{W}_{\aalpha,\bbeta}} \quad \text{ or } \quad P^{\MISCquad}_{\aalpha,\bbeta} = \frac{\mathcal{E}^{\MISCquad}_{\aalpha,\bbeta}}{\mathcal{W}_{\aalpha,\bbeta}}. 
\end{equation}
An effective strategy to build adaptively a MISC approximation can then be broadly described as follows: given the MISC approximation
associated to a multi-index set $\Lambda$, a new MISC approximation is obtained by adding to $\Lambda$ the multi-index $[\aalpha,\bbeta] \not \in \Lambda$
with the largest profit (either $P^{\MISCsurr}$ or $P^{\MISCquad}$, depending on the goal of the simulation), such that $\Lambda \cup \{[\aalpha,\bbeta]\}$ is downward closed.}
In practice the implementation reported in Algorithm \ref{algo:misc_implementation} is used:
it makes use of an auxiliary multi-index set, i.e.\ the margin of a multi-index set $\Lambda$, $\text{Mar}(\Lambda)$,
which is defined as the set of multi-indices that can be reached ``within one step'' from $\Lambda$ (see Fig.~\ref{fig:index_sets}b)
\begin{equation*}
\text{Mar}(\Lambda) = \{\ii \in \mathbb{N}^{d+N} \text{ s.t. } \ii = \jj + \ee_k \text{ for some } \jj \in \Lambda \text{ and some } k \in \{1,\ldots,d+N\} \}.    
\end{equation*}
This algorithm was first proposed in the context of sparse-grids quadrature in \cite{gerstner.griebel:adaptive}
and its MISC implementation was first proposed in \cite{jakeman2019adaptive}.
It is an a-posteriori algorithm and as such it determines the error contribution $\mathcal{E}^{\mathcal{F}}_{\aalpha,\bbeta}$ of $[\aalpha,\bbeta]$
{after} having added $[\aalpha,\bbeta]$ to the grid. Therefore, at the end of the algorithm we do not have $\MISCsurr_{\Lambda}$, but actually $\MISCsurr_{J}$,
where $J$ is the set of all indices whose profit has been computed, and clearly $\Lambda \subseteq J$:
the richer approximation $\MISCsurr_{J}$ is thus actually returned in practical implementations instead of $\MISCsurr_{\Lambda}$. Finally,
note that many stopping criteria can be considered (and possibly used simultaneously), which typically check that computational work,
error contributions or profit estimator are below a desired threshold.

\begin{algorithm}
	\Fn{\AlCapSty{Multi-Index Stochastic Collocation}}{
          $\Lambda = \{ [\bm{1},\bm{1}] \}$\;
          $ J = \{[\bm{1},\bm{1}] \}$\tcp*[r]{set of multi-indices whose profit has been computed, $\Lambda \subseteq J$}
          $ L = \emptyset$\tcp*[r]{set of candidate multi-indices to be added to $\Lambda$}
		Compute MISC estimate $\MISCsurr_{J}$ as in Eq.~\eqref{eq:misc}\;
		\While{ stopping criteria are not met }{
                  \For(\tcp*[f]{for short, $\jj = [\aalpha , \bbeta]$}){ 
                    $\jj \in \text{Mar}(\Lambda)$ {\bf and} $\Lambda \cup \{\jj\}$ downward closed }{
                         \If(){$\jj \not \in J$}{
                         Compute MISC estimate $\MISCsurr_{ J \cup \{\jj\} }$ as in \eqref{eq:misc}\;
                         Compute error contribution $\mathcal{E}^{\MISCquad}_{\jj}$ (or $\mathcal{E}^{\MISCsurr}_{\jj}$) as in Eq.~\eqref{eq:error_contr_quad} (or Eq.~\eqref{eq:error_contr_Linf})\;
                         Compute work contribution $\mathcal{W}_{\jj}$ as in Eq.~\eqref{eq:work_contr}\;
                         Compute profit $P^{\MISCsurr}_{\jj} = \mathcal{E}^{\MISCsurr}_{\jj} / \mathcal{W}_{\jj}$ (or $P^{\MISCquad}_{\jj} = \mathcal{E}^{\MISCquad}_{\jj} / \mathcal{W}_{\jj}$)\; 
                         $ J = J\cup \{\jj\}$\;
                         $ L = L \cup \{\jj\}$\;
                         }
		}
	Choose $\ii \in L$ with the highest profit\;
	$\Lambda = \Lambda \cup \{\ii\}$\;
        $L = L \setminus \{\ii\}$\;
	} 
      }
      \Return{ $\Lambda, J, \MISCsurr_{J}$}
      \caption{MISC implementation}
      \label{algo:misc_implementation}
\end{algorithm}

\subsection{Adaptive multi-fidelity Stochastic Radial Basis Functions (SRBF)} \label{sect:SRBF}
In this section each of the components of the adaptive multi-fidelity SRBF surrogate model are discussed.
In particular, we emphasize that here the word ``stochastic'' denotes not only the fact that we are sampling parameters that are {affected by uncertainty}, but also to the fact that
the SRBF method treats one of its hyper-parameters as a random variable, as will be clear later on.

\subsubsection{SRBF surrogate model}
  
\begin{figure}[!t]
\centering
\includegraphics[width=0.495\textwidth]{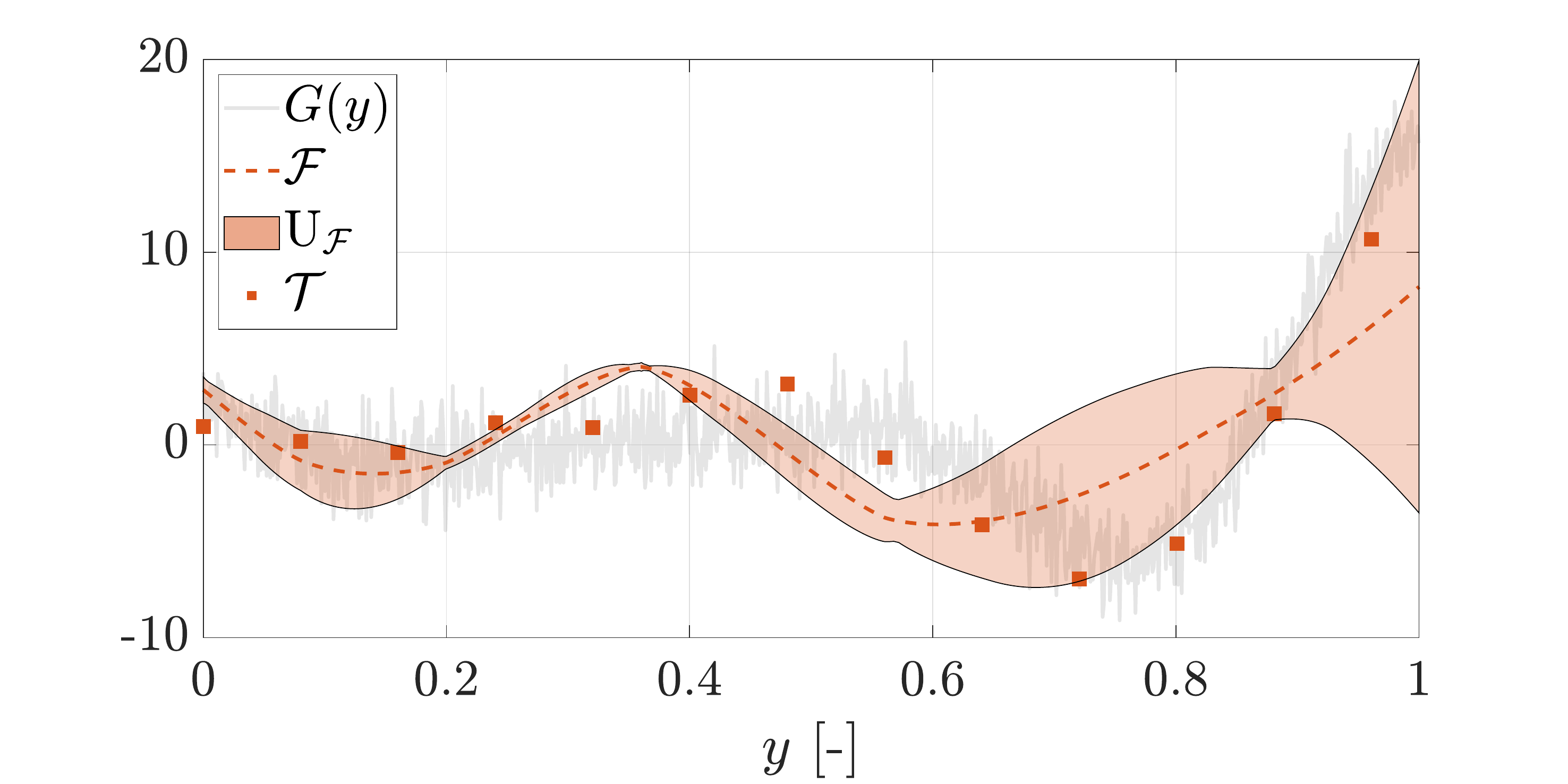}
\caption{SRBF example with least-squares regression.}\label{fig:lsrbf}
\end{figure}

Given a training set {$\mathcal{T}=\{ \left({\bf y}_i,G({\bf y}_i)\right) \}_{i=1}^{\mathcal{J}}$} and normalizing the uncertain parameters domain into a unit hypercube, the RBF prediction is here based on a power function kernel and reads
\begin{equation}\label{eq:rbf}
f\left({\bf y},\tau\right)=\sum^{\mathcal{K}}_{j=1} w_j ||{\bf y}-{\bf c}_j||^{\tau},
\end{equation}
where $w_j$ are unknown coefficients, ${\bf c}_j$ are {$\mathcal{K}$ points in $\Gamma$ called} RBF centers,
and $\tau\sim \textrm{unif}[\tau_{\min},\tau_{\max}]$ is a stochastic tuning parameter that follows a uniform distribution. {The range of $\tau$ is defined within $\tau_{\min}=1$ and $\tau_{\max}=3$, where $\tau=1$ provides a polyharmonic spline of first order (linear kernel) \cite{gutmann2001} and $\tau=3$ provides a polyharmonic spline of third order (cubic kernel) \cite{forrester2009}.} Note that the choice of the distribution for $\tau$ is arbitrary and, from a Bayesian viewpoint, this represents the degree of belief in the definition of the tuning parameter.   
The SRBF {surrogate model $\mathcal{F}\left({\bf y}\right)$} is computed as the expected value (approximated by {Monte Carlo})
of $f$ over $\tau$ \citep{volpi2015-SMO}:
\begin{eqnarray}\label{eq:SRBF}
G({\bf y}) \approx {\mathcal{F}\left({\bf y}\right)}={\mathbb{E}_{\tau}\left[f\left({\bf y},\tau\right)\right]}\approx\frac{1}{\Theta}\sum\limits_{i=1}^{\Theta}f\left({\bf y},\tau_i\right),
\end{eqnarray}
where $\Theta$ is the number of samples for $\tau$, here {set} equal to $1000$.
{To give more flexibility to the method,}
the coordinates of the RBF centers ${\bf c}_j$ are not {a-priori set to be} coincident with the training points, but {rather chosen by a} $k$-means clustering {algorithm applied to the training points, see} \citep{lloyd1982-IEEE}.
{Several values of the number of centers $\mathcal{K} \leq \mathcal{J}$ are tested and their} optimal number
{ $\mathcal{K}^*$ is chosen} by minimizing a leave-one-out cross-validation (LOOCV) metric, see \citep{li2017-SMO}.
{In details, letting $g_{i,\mathcal{K}}({\bf y})$, $i=1,\ldots,\mathcal{J}$ be the surrogate models with $\mathcal{K}$ centers trained on the whole training set $\mathcal{T}$ but the $i$-th point, $\mathcal{K}^*$ is defined as:}
\begin{equation}\label{eq:loocv}
{\mathcal{K}^*} = {\underset{{\mathcal{K}\in { C}}}{\rm argmin}} \,\mathrm{RMSE}({\mathcal{K}}), 
\end{equation}
{where $\mathcal{K} \leq \mathcal{J}, \mathcal{K} \in C \subset \mathbb{N}$} and RMSE($\mathcal{K}$) is the root mean square error of the $\mathcal{J}$ leave-one-out models $g_{1,\mathcal{K}},\ldots,g_{\mathcal{J},\mathcal{K}}$ at the point that is being left out for each $g_{i,\mathcal{K}}$:
\begin{equation}\label{eq:rmse}
  \mathrm{RMSE}({\mathcal{K}}) =
  \sqrt{\dfrac{1}{\mathcal{J}}\sum_{i=1}^{\mathcal{J}}
    \left(G(\mathbf{y}_i) - g_{{i,\mathcal{K}}}(\mathbf{y}_i)\right)^2 },
  \quad {{\bf y}_i \in \mathcal{T}.}
\end{equation}
{Clearly, once the optimal number of centers $\mathcal{K}^*$ is chosen, the whole set of points is used for the construction of the final surrogate model.}
Whenever the number of RBF centers is lower than the training set size ($\mathcal{K} < {\mathcal{J}}$),
the coefficients $w_j$ in Eq.~\eqref{eq:rbf} are determined through a least-squares regression by solving 
\begin{equation}\label{eq:lsrbf}
{\bf w}=\left( {\bf A}^{\mathsf{T}} {\bf A} \right)^{-1} {\bf A}^{\mathsf{T}} {\bf b},
\end{equation}
with ${\bf w}=[w_1, \ldots,w_{\mathcal{K}}]^{\mathsf{T}}$, ${\bf A}_{ij}=||{\bf y}_i-{\bf c}_j||^{\tau}$, {$1\leq i \leq \mathcal{J}$, $1\leq j \leq \mathcal{K}$} and
${\bf b}=[G({\bf y}_1),\ldots,G({\bf y}_{\mathcal{J}})]^{\mathsf{T}}$;
otherwise when the optimal number of RBF centers equals the training set size, exact interpolation at the training points
{($f({\bf y}_i,\tau)=G({\bf y}_i)$)}
is imposed and Eq.~\eqref{eq:lsrbf} reduces to
\begin{equation}\label{eq:interp_rbf}
{\bf w}={\bf A^{-1}b},
\end{equation}   
with ${\bf c}_j={\bf y}_j$.
{Having less RBF centers than training points and employing the least-squares approximation in Eq.~\eqref{eq:lsrbf}
  to determine the coefficients $w_j$ is particularly helpful when the training data are affected by noise.}
An example of least-squares regression is shown in Fig.~\ref{fig:lsrbf}.

The uncertainty {$U_{\mathcal{F}}\left({\bf y}\right)$}
associated with the SRBF {surrogate model} prediction is
quantified by the 95\%-confidence band of the cumulative density function (CDF) of $f({\bf y}, \tau)$ with respect to $\tau$ for fixed $\yy$ as follows
\begin{equation}\label{eq:USRBF}
{U_{\mathcal{F}}}({\bf y})={\rm CDF}^{-1}(0.975;{\bf y})-{\rm CDF}^{-1}(0.025;{\bf y}),
\end{equation}
with
\begin{equation*}
{\rm CDF}(\lambda;{\bf y}) \approx \frac{1}{\Theta}\sum\limits_{i=1}^{\Theta} H[\lambda-f({\bf y},\tau_i)],
\end{equation*}
where $H(\cdot)$ is the Heaviside step function.

\subsubsection{Multi-fidelity approach}

In this section we {restrict to the case of the CFD grid generation being} controlled by a scalar value $\alpha$, i.e., the QoI computed with the $\alpha$-th grid is denoted by $G_\alpha$, $\alpha = 1,\ldots, M$. 
{The multi-fidelity approximation {of $G$} is adaptively built following the approach introduced in \cite{serani2019-MARINE} and extended to noisy data in \cite{wackers2020-SNH}.}
Extending the definition of the training set to an arbitrary number $M$ of fidelity levels as {$\{\mathcal{T}_{\alpha}\}_{\alpha=1}^M$, with each $\mathcal{T}_{\alpha} = \{\left({\bf y}_j,G_{\alpha}({\bf y}_j)\right)\}_{j=1}^{\mathcal{J}_{\alpha}}$, the multi-fidelity approximation ${\SRBFsurr}_{\alpha}(\mathbf{y})$ of {$G(\mathbf{y})$} reads 
\begin{equation}
\label{eq:MLGeneral}
{\SRBFsurr}_{\alpha}(\mathbf{y}) :=  \mathcal{F}_1(\mathbf{y})+\sum_{i=1}^{{\alpha}-1}\epsilon_i(\mathbf{y}), 
\end{equation}
where $\mathcal{F}_1$ is the single-fidelity surrogate model associated to the lowest-fidelity training set (constructed as in Eq.~\eqref{eq:SRBF}), {and $\epsilon_i(\mathbf{z})$ is the inter-level error surrogate with associated training set {$\mathcal{E}_i=\{({\bf z},{\phi}-\mathcal{S}_{i}({\bf z}))\,|\,({\bf z}, {\phi})  \in \mathcal{T}_{i+1} \}$}}. An example of the multi-fidelity approximation with two fidelities is shown in Fig.~\ref{fig:mfm}.
\begin{figure}[!t]
\centering
\includegraphics[width=0.495\textwidth]{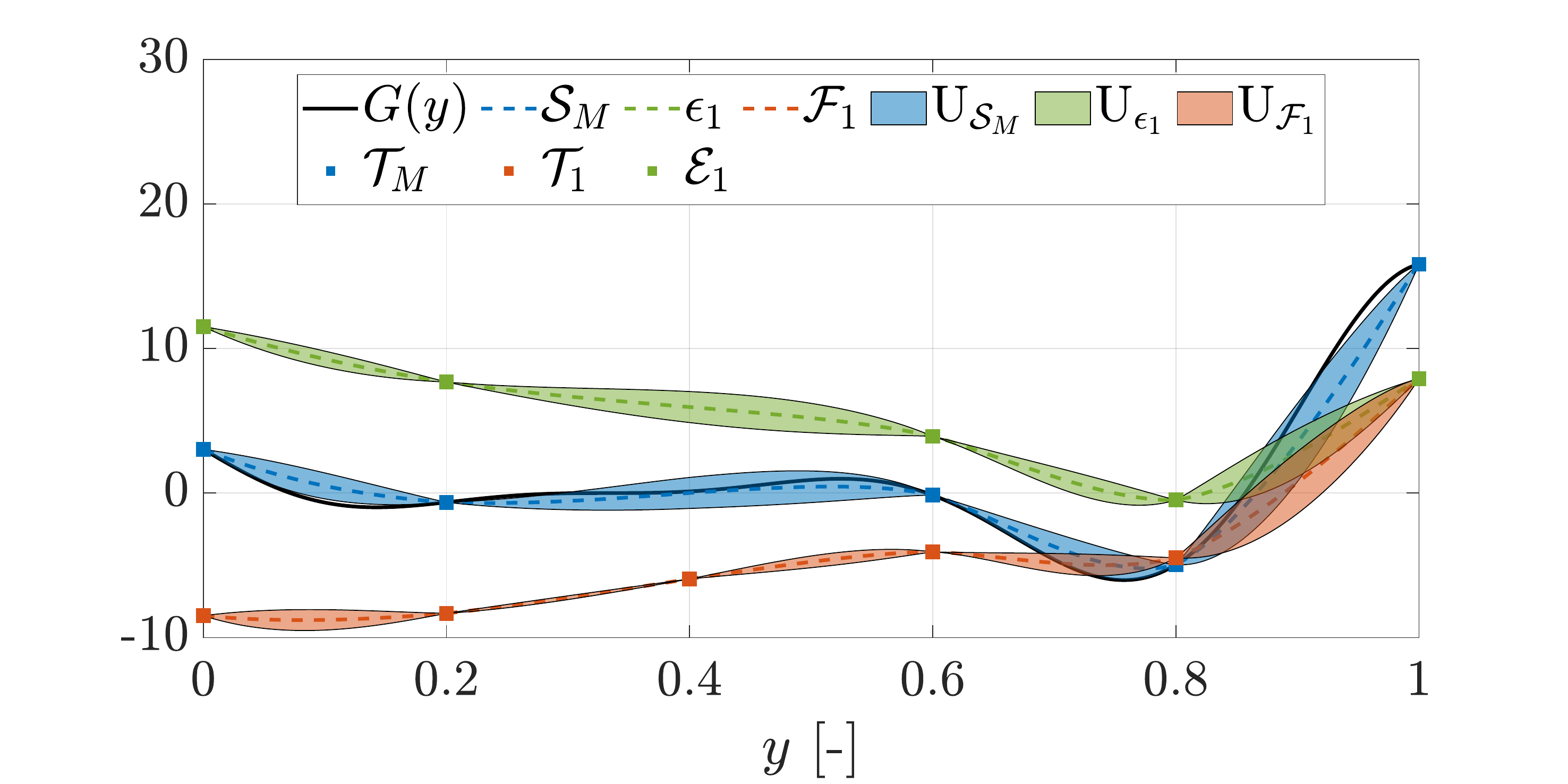}
\caption{Example of multi-fidelity surrogate with $M=2$ and exact interpolation at the training points.}\label{fig:mfm}
\end{figure}
\begin{figure}[!b]
\centering
\subfigure[]{\includegraphics[width=0.495\textwidth]{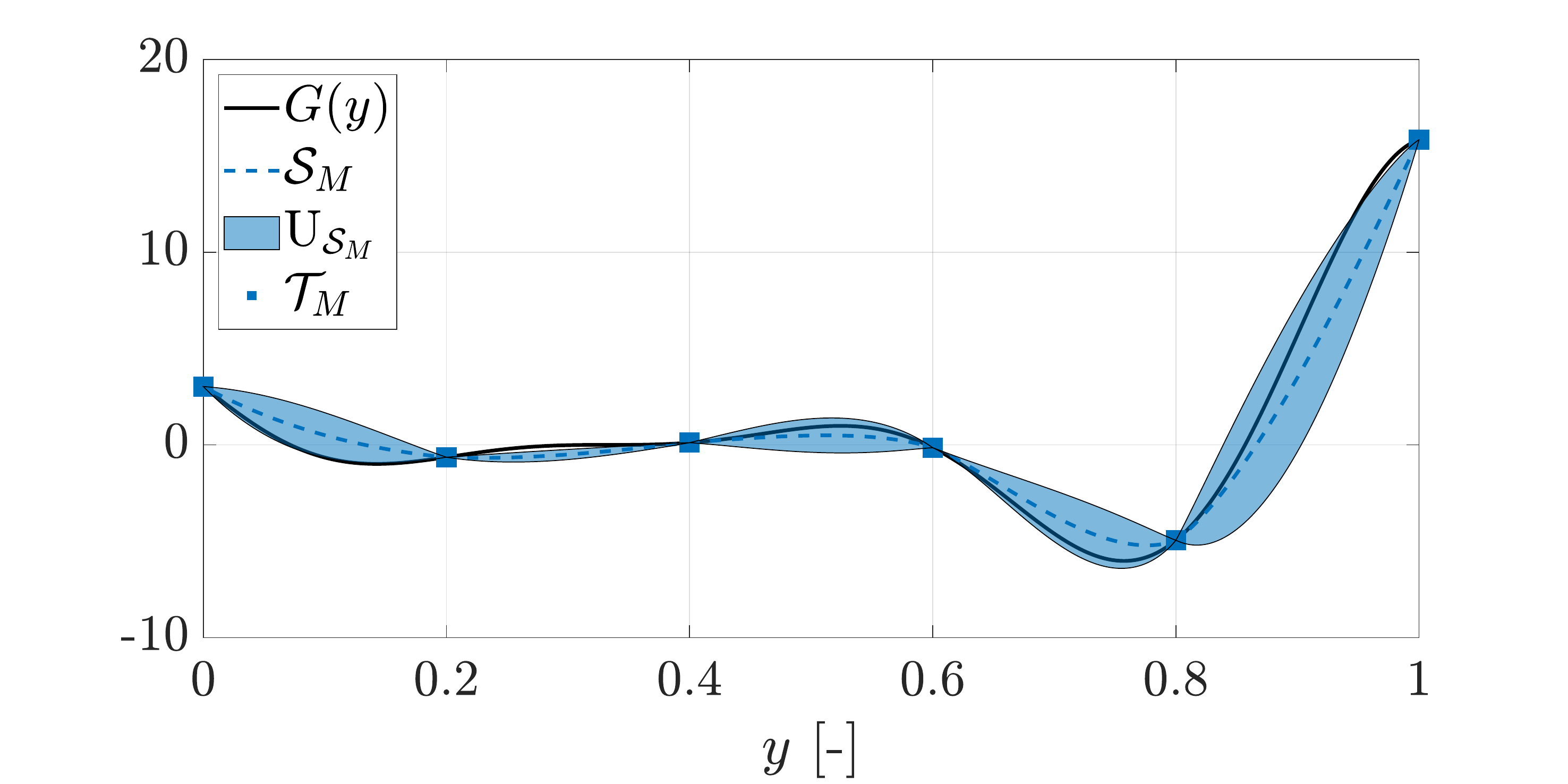}}
\subfigure[]{\includegraphics[width=0.495\textwidth]{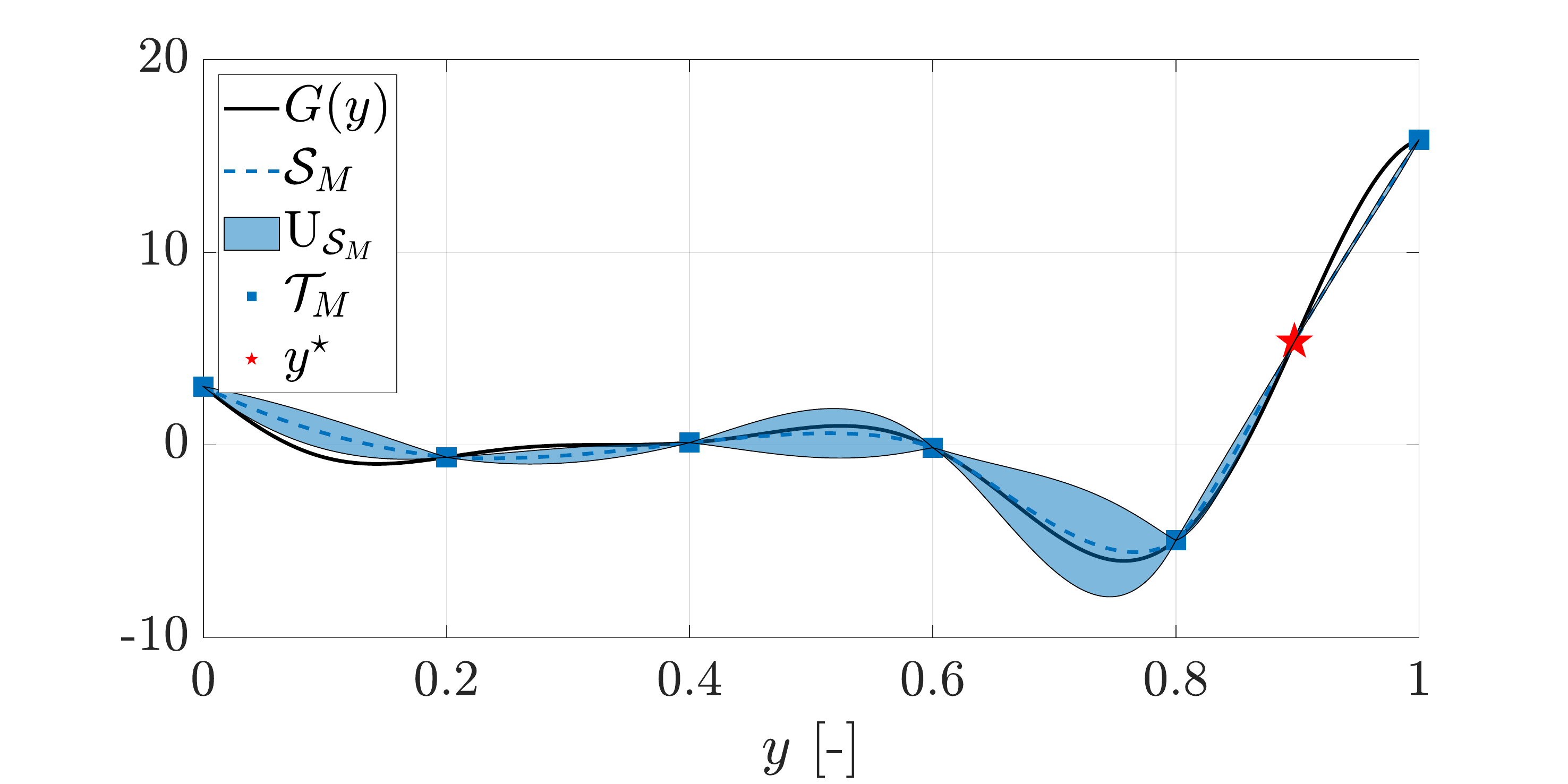}}
\caption{Example of the adaptive sampling method using {a single-fidelity training set consisting of noiseless evaluations}: (a) shows the initial SRBF with the associated prediction uncertainty and training set; (b) shows the position of the new training point, the new SRBF prediction, and its uncertainty.}\label{fig:muas}
\end{figure}

\begin{algorithm}[!t]
    \Fn{\AlCapSty{Multi-Fidelity SRBF for numerical quadrature}}{
        Define the initial training sets $\{\mathcal{T}_{\alpha}\}_{\alpha=1}^M$, {with} $\mathcal{T}_{\alpha}=\{\left({\bf y}_j,G_{\alpha}({\bf y}_j)\right)\}_{j=1}^{\mathcal{J}_{\alpha}^0}$ \; 
        $t = 0$ \;

        \While(\tcp*[f]{{$t$-th iteration of the adaptive sampling}}){stopping criteria are not met}{
            
            Perform LOOCV and find ${\mathcal{K}^{*,t}_{1}}$ as in Eq.~\eqref{eq:loocv} with: \\
                \eIf{$t=0$}{
                           ${C}=[1,\mathcal{J}^0_{1}]$  \;
                        }
                {
                          ${C}=[{\mathcal{K}_1^{*,t-1}},        {\mathcal{K}_1^{*,t-1}}+1]$ \;
                        }

            Construct the SRBF surrogate $\mathcal{F}_1(\mathbf{y})$ as in Eq.~\eqref{eq:SRBF} \tcp*{low-fidelity approximation}
            Compute the prediction uncertainty $U_{\mathcal{F}_1}(\mathbf{y})$ as in Eq.~\eqref{eq:USRBF} \;
            
            \For(\tcp*[f]{evaluate surrogates of the inter-level errors}){$\alpha=2,\dots,M$}{
                Compute the inter-level errors $\mathcal{E}_{\alpha}$ \;
                    Perform LOOCV and find {${\mathcal{K}^{*,t}_{\alpha}}$} as in Eq.~\eqref{eq:loocv} with: \\
                    \eIf{$t=0$}{
                        ${C}=[1,\mathcal{J}^0_{\alpha}]$  \;
                    }
                    {
                        ${C}=[{\mathcal{K}_{\alpha}^{*,t-1}},        {\mathcal{K}_{\alpha}^{*,t-1}}+1]$ \;
                    }
                    Construct the SRBF surrogates of the inter-level errors $\epsilon_{\alpha}(\mathbf{y})$ as in Eq.~\eqref{eq:SRBF} \;
                    Compute the prediction uncertainty $U_{\epsilon_{\alpha}}(\mathbf{y})$ as in Eq.~\eqref{eq:USRBF} \; 
                    Construct the multi-fidelity SRBF surrogate ${\SRBFsurr}_{\alpha}(\mathbf{y})$ as in Eq.~\eqref{eq:MLGeneral} \; }
                        
            Construct the multi-fidelity approximation ${\SRBFsurr}_M(\mathbf{y})$ as in Eq.~\eqref{eq:ML-MF} \tcp*{MF approximation}
            Compute the multi-fidelity prediction uncertainty $U_{{\SRBFsurr}_M}(\mathbf{y})$ as in Eq.~\eqref{eq:ML-MF} \;
            Set $\mathcal{T}_{\alpha0}=\mathcal{T}_{\alpha}$ with $ \alpha=1, \dots, M$ \tcp*{store current training sets}
                Set $\mathcal{J}^t_{\alpha0}=\mathcal{J}^t_{\alpha}$ with $ \alpha=1, \dots, M$ \tcp*{store current training sets size}
            
            \For(\tcp*[f]{perform parallel infill}){$j=1,\dots,p$}{                          
                Find ${\bf y}_j^{\star}={\underset{{\bf y} \in \Gamma}{\rm argmax}}[U_{{\SRBFsurr}_M}({\bf y})]$ \;
                Find $k_j=\mathrm{maxloc}\left[\mathbf{U}(\mathbf{y}_j^\star)\right]$  as in Eq.~\eqref{eq:choice} \;
                Update the training sets $\{\mathcal{T}_{\alpha}\}_{\alpha=1}^{k_j }\cup \{\left({\bf y}_j^{\star},{\SRBFsurr}_{\alpha}({\bf y}_j^{\star})\right)\}_{\alpha=1}^{k_j}$ {\tcp*{consider exact prediction}}                                           
                Update the training sets size  $\{\mathcal{J}_{\alpha}^{t+j}\}_{\alpha=1}^{k_j}=\{\mathcal{J}_{\alpha}^t\}_{\alpha=1}^{k_j}+1$ \;
            } 
            
            Set $\mathcal{T}_{\alpha}=\mathcal{T}_{\alpha0}$ with $\alpha=1, \dots, M$ \tcp*{restore previuos training sets}
            Set $\mathcal{J}^t_{\alpha}=\mathcal{J}^t_{\alpha0}$ with $ \alpha=1, \dots, M$ \tcp*{restore previous training sets size}
            
            \For(\tcp*[f]{perform new simulations}){$j=1,\dots,p$}{                       
                Evaluate $\{G_{\alpha}({\bf y}_j^{\star})\}_{\alpha=1}^{k_j}$ \; 
                Update the training sets $\{\mathcal{T}_{\alpha}\}_{\alpha=1}^{k_j} \cup \{\left({\bf y}_j^{\star},G_{\alpha}({\bf y}_j^{\star})\right)\}_{\alpha=1}^{k_j}$\;                                      
                Update the training sets size $\{\mathcal{J}_{\alpha}^{t+j}\}_{\alpha=1}^{k_j}=\{\mathcal{J}_{\alpha}^t\}_{\alpha=1}^{k_j}+1$\;
            }
            $t=t+1$ \tcp*{move to the next adaptive sampling iteration}                                              
        } 
    } 
    \caption{Adaptive multi-fidelity SRBF implementation}
    \label{algo:mfm_implementation}
\end{algorithm}

Assuming that the uncertainty associated to the prediction of the lowest-fidelity $U_{\mathcal{F}_1}$ and inter-level errors $U_{\epsilon_i}$ as uncorrelated, the multi-fidelity approximation ${\SRBFsurr}_M(\mathbf{y})$ of {$G(\mathbf{y})$} and its uncertainty $U_{{\SRBFsurr}_M}$ read 
\begin{equation}
\label{eq:ML-MF}
G(\mathbf{y})\approx{\SRBFsurr}_M(\mathbf{y})= \mathcal{F}_1(\mathbf{y})+\sum_{i=1}^{M-1}\epsilon_i(\mathbf{y})
~~~~~
\mathrm{and}
~~~~~
U_{{\SRBFsurr}_M}(\mathbf{y})=
\sqrt{U^2_{\mathcal{F}_1}(\mathbf{y})+\sum_{i=1}^{M-1}U^2_{{\epsilon}_i}(\mathbf{y})}.
\end{equation}
%

\subsubsection{Adaptive sampling approach}\label{sect:adaptive_SRBF}

Upon having evaluated $U_{{\SRBFsurr}_M}$, the multi-fidelity surrogate is then updated adding a new training point following a two-steps procedure: firstly, the coordinates of the new training point ${\bf y}^\star$ are identified based on the SRBF maximum uncertainty, see \cite{serani2019-IJCFD}, solving the maximization problem: 
\begin{eqnarray}\label{eq:MUAS}
{\bf y}^{\star}={\underset{{\bf y \in \Gamma}}{\rm argmax}}[U_{{\SRBFsurr}_M}({\bf y})].
\end{eqnarray}
An example (with one fidelity only) is shown in Fig.~\ref{fig:muas}. Secondly, once ${\bf y}^\star$ is identified, the training set/sets {to be} updated with {the} new training point {$  \left(\mathbf{y}^\star, G_{\alpha}(\mathbf{y}^\star)\right)$} {are $\mathcal{T}_{\alpha}$} with $\alpha=1,\dots,k$, where $k$ is defined as 
\begin{equation}\label{eq:choice}
 k=\mathrm{maxloc}\left[\mathbf{U}(\mathbf{y}^\star)\right]
 ~~~~~
\mathrm{and}
~~~~~
\mathbf{U}(\mathbf{y}^\star)\equiv\{ U_{\mathcal{F}_1}(\mathbf{y}^\star)/\gamma_1, U_{{{\epsilon}}_1}(\mathbf{y}^\star)/\gamma_2, ...,U_{{{\epsilon}}_{M-1}}(\mathbf{y}^\star)/\gamma_M \},
\end{equation}
with $\gamma_{\alpha}$ being the computational cost associated to the $\alpha$-th level.

{In the present work, the adaptive sampling procedure starts with five training points (for each fidelity level) located at the domain center and at the centers of each boundary of $\Gamma$. 
Furthermore, to avoid abrupt changes in the SRBF prediction from one iteration to the next one, the search for the optimal {number of centers for the $\alpha$-th fidelity} {$\mathcal{K}^*_\alpha$} can be constrained. Herein, 
{at every adaptive sampling iteration, the problem in Eq.~\eqref{eq:loocv} is solved assuming $\mathcal{K}$ to be either equal to the number of centers at the previous iteration or incremented by 1, i.e.\ $\mathcal{C}=[\mathcal{K}^{*,t-1}_\alpha,\mathcal{K}^{*,t-1}_\alpha+1]$,} 
except for the first iteration where no constraint is imposed.

A deterministic version of the particle swarm optimization algorithm \citep{serani2016-ASC} is used for the solution of the optimization problem in Eq.~\eqref{eq:MUAS}.

The adaptive sampling is therefore inherently sequential (the uncertainty changes every time a new point is added), but this is sub-optimal whenever the numerical simulations can be performed with an hardware capable of running $p$ simulations simultaneously. In this case, it would be ideal to identify $p$ training points where the models $G_{\alpha}$ can be run in parallel, instead of running them one after the other. To this end, we follow a parallel-infill procedure, i.e.\ we perform $p$ ``guessing steps'': the adaptive sampling procedure is repeated $p$ times replacing the evaluations of the actual model $G_{\alpha}$ with the evaluations of the multi-fidelity models $\mathcal{S}_{\alpha}$. This replacement significantly speeds up the $p$ steps, since the true models $G_{\alpha}$ are not evaluated  at this stage. Upon doing these $p$ guessing steps, the actual  $G_{\alpha}$ are evaluated all at once (i.e. in parallel) at the $p$ training points just obtained and these evaluations replace the corresponding multi-fidelity evaluations in the training set.

Finally, numerical quadrature is used on the multi-fidelity SRBF surrogate model to estimate the statistical moments of the QoI. 

\section{Numerical tests}\label{sect:numerical_tests}
In this section two numerical tests are considered. First, the performances of the MISC and SRBF method
are compared on an analytical test, and then the main problem of this work, i.e., the RoPax ferry problem mentioned in the introduction, is discussed.
In the analytic example, Taylor expansions of increasing order are considered as different fidelities
to be employed, while in the RoPax problem the fidelities are obtained by using different grid refinements, as will be clearer later.
Both problems consider uniformly distributed {uncertain parameters}.
Before entering the detailed discussion of the two tests, an overview of the {error} metrics used to carry out the comparison is given in the following.

\subsection{{Error} metrics}\label{sect:metrics}

The performance of MISC and SRBF are assessed by comparing the convergence of both methods to a reference solution according to several error metrics. 
The specific choice of the reference solution (denoted below by $G_{\text{ref}}$) for each test will be detailed in the corresponding sections. 
In the following list, we use the symbol $\mathcal{S}$ for both the MISC and SRBF surrogate models for sake of compactness,
i.e. $\MISCsurr = \MISCsurr_{\Lambda}$ for MISC (cf. Eq.~\eqref{eq:misc}) and $\mathcal{S} = {\SRBFsurr}_M$ for SRBF (cf. Eq.~\eqref{eq:ML-MF}).
  \begin{enumerate}
  \item {Relative error} between the first four centered moments (mean, variance, skewness, kurtosis) of the MISC/SRBF approximations and those of the reference solution:
  \begin{equation}\label{eq:err_moments}
     	err_i = \frac{\lvert \text{Mom}_i[\mathcal{S}] 	- \text{Mom}_i[G_{\text{ref}}] \rvert }{ \lvert \text{Mom}_i[G_{\text{ref}}] \rvert}, \quad i=1,\ldots,4,
  \end{equation}
    where $\text{Mom}_i[\mathcal{S}]$ denotes the MISC/SRBF approximation of the $i$-th centered moment (computed by the quadrature rule associated to MISC/SRBF) and $\text{Mom}_i[G_{\text{ref}}]$ the approximation of the $i$-th centered moment of the reference solution computed by a suitable quadrature rule (more details will be given later).
  \item {Relative error in discrete $L_2$ and $L_\infty$ norm} between the MISC/SRBF surrogates and the reference solution, i.e., sample mean square error and largest point-wise prediction error, respectively; the differences are evaluated at a set of $n=10000$ random points $\yy_i \in \Gamma$. In formulas, 
     \begin{equation} \label{eq:err_norm}
      err_{L_2}=\frac{\sqrt{\frac{1}{n} \sum_{i=1}^{n} \left[\mathcal{S}(\yy_i) - {G}_{\text{ref}}(\mathbf{y}_i)\right]^2 }}{\sqrt{\frac{1}{n} \sum_{i=1}^{n} {G}_{\text{ref}}(\mathbf{y}_i)^2 }}, \,\,\,
      err_{L_{\infty}}=\frac{\max_{i=1,\ldots,n}\left(\left|\mathcal{S}(\mathbf{y}_i) - {G}_{\text{ref}}(\mathbf{y}_i)\right|\right)}{\max_{i=1,\ldots,n} {G}_{\text{ref}}(\mathbf{y}_i)}.
    \end{equation}

    \item A visual comparison of the PDFs obtained by Matlab's {\it ksdensity} function, using again as input the $10000$ points used before. 
  \item Convergence of the CDF approximated by MISC/SRBF to the CDF of the reference solution,
    as measured by the Kolmogorov--Smirnov (KS) test statistic. In details, we evaluate an approximation of the quantity
    \begin{equation}\label{eq:KS_stat}
      T = \sup_{t \in \textrm{range}[\mathcal{S},G_{\text{ref}}]} \,\lvert \text{CDF}^{\mathcal{S}}(t) - \text{CDF}^{G_{\text{ref}}}(t) \rvert,
    \end{equation}
    where $\textrm{range}[\mathcal{S},G_{\text{ref}}]$ is the largest common range of values taken by $\mathcal{S}$ and $G_{\text{ref}}$,
    $\text{CDF}^{\mathcal{S}}$ and $\text{CDF}^{G_{\text{ref}}}$ are the empirical CDFs obtained by the set used before of $10000$ random samples 
    of the MISC/SRBF surrogate models and reference model, respectively. We then check that $T$ converges to zero as the surrogate models get more and more accurate.
    The values of $T$ reported in the next sections are obtained with the Matlab's {\it kstest2} function. 
  \end{enumerate}

{We emphasize that the adaptive criteria that drive the construction of the MISC and SRBF approximations need not match the error metrics above (compare them against Eqs.\ \eqref{eq:profits}, \eqref{eq:error_contr_quad} and \eqref{eq:error_contr_Linf} for MISC, and Eq.\ \eqref{eq:MUAS} for SRBF). It is actually interesting to investigate how MISC and SRBFs converge when monitoring error norms that are not aligned with the adaptivity criteria. }

\subsection{Analytical problem}\label{sect:problem_description_analytic}
\subsubsection{Formulation}

As analytical test, a two-dimensional function is chosen. This function is designed to be representative of the RoPax problem: the input {parameters} $\yy $ are independent and have a uniform distribution, and the function is non-linear, non-polynomial and monotonic. In details, it is defined as
\begin{equation*}\label{eq:AnProb}
G({\bf y}) = G(y_1,y_2)=\sin \left( \frac{\exp(y_1+y_2)}{5} \right),
\end{equation*}
with ${\bf y} \in [0,1]^2$ .
To provide a range of fidelities $G_{\alpha}({\bf y})$  for $G({\bf y})$, Taylor expansions of order $\alpha$ {of the argument of the $\sin(\cdot)$ function, that is $\frac{\exp(y_1+y_2)}{5}$,} are performed for $\alpha=1, \ldots, 6$ in the neighborhood of ${\bf y}=(0,0)$. The sixth-order Taylor expansion $G_{6}({\bf y})$ is then considered as the highest-fidelity and the first order $G_{1}({\bf y})$ as the lowest-fidelity. 
Figure~\ref{fig:anprob} shows the true function $G({\bf y})$ and the approximations $G_{6}({\bf y})$ and $G_{1}({\bf y})$. 
{We mention in-passing that the sixth-order Taylor expansion is almost indistinguishable from the true function in the range of $y_1,y_2$ considered, whereas the low-fidelity function is significantly different and does not show any change in curvature. Note that the difference between the sixth-order and the true function is actually irrelevant for our purposes since we never consider the true function in the numerical tests: errors are indeed computed with respect the sixth-order approximation, in analogy with PDE-based problems where no exact closed-formula solution is available, and the ``ground-truth'' is usually taken as a ``refined-enough'' solution.} 

\begin{figure}[tp]
\centering
\subfigure[True function]{\includegraphics[width=0.3\textwidth]{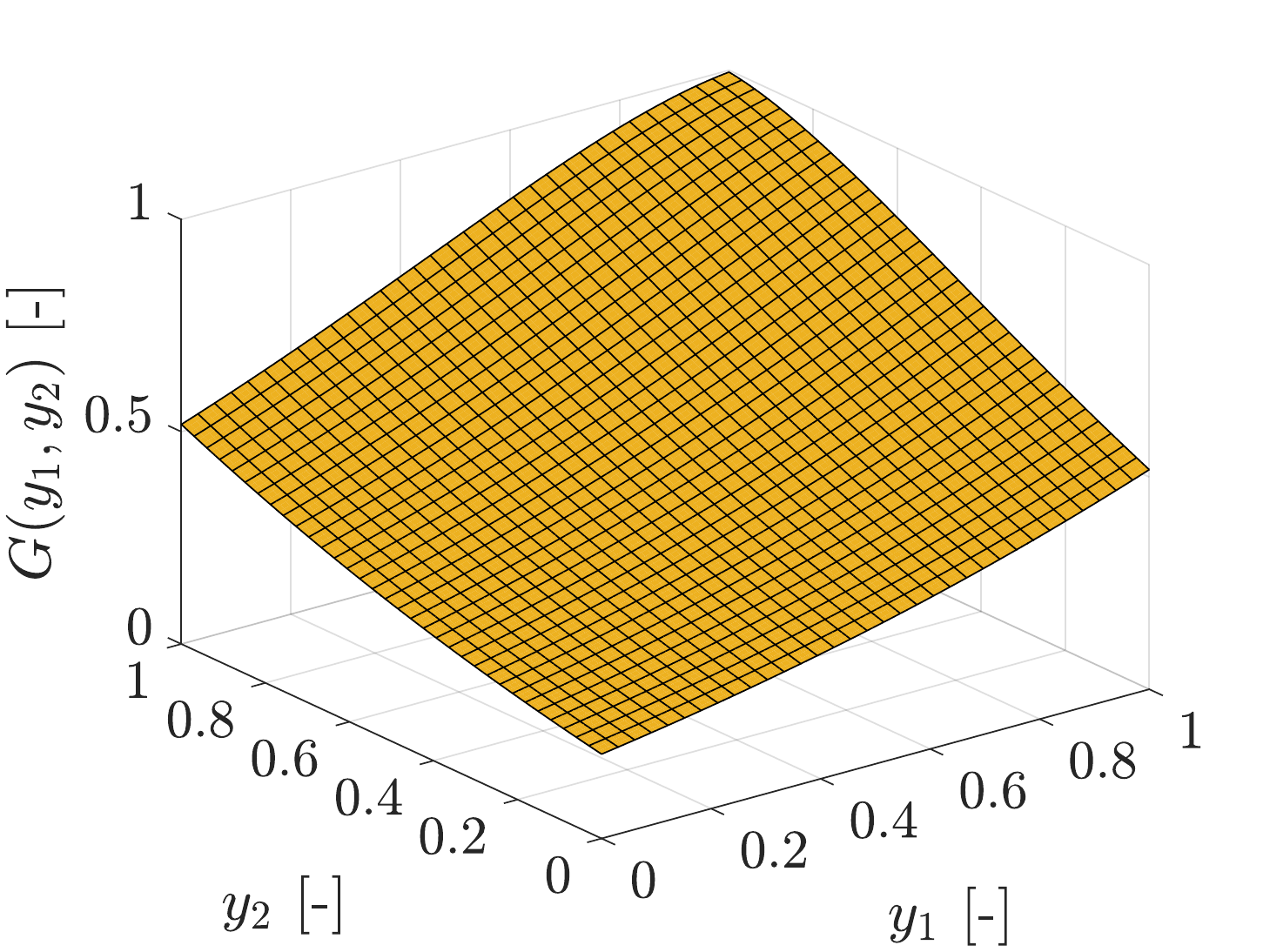}}
\subfigure[Highest-fidelity function $G_{6}({\bf y})$]{\includegraphics[width=0.3\textwidth]{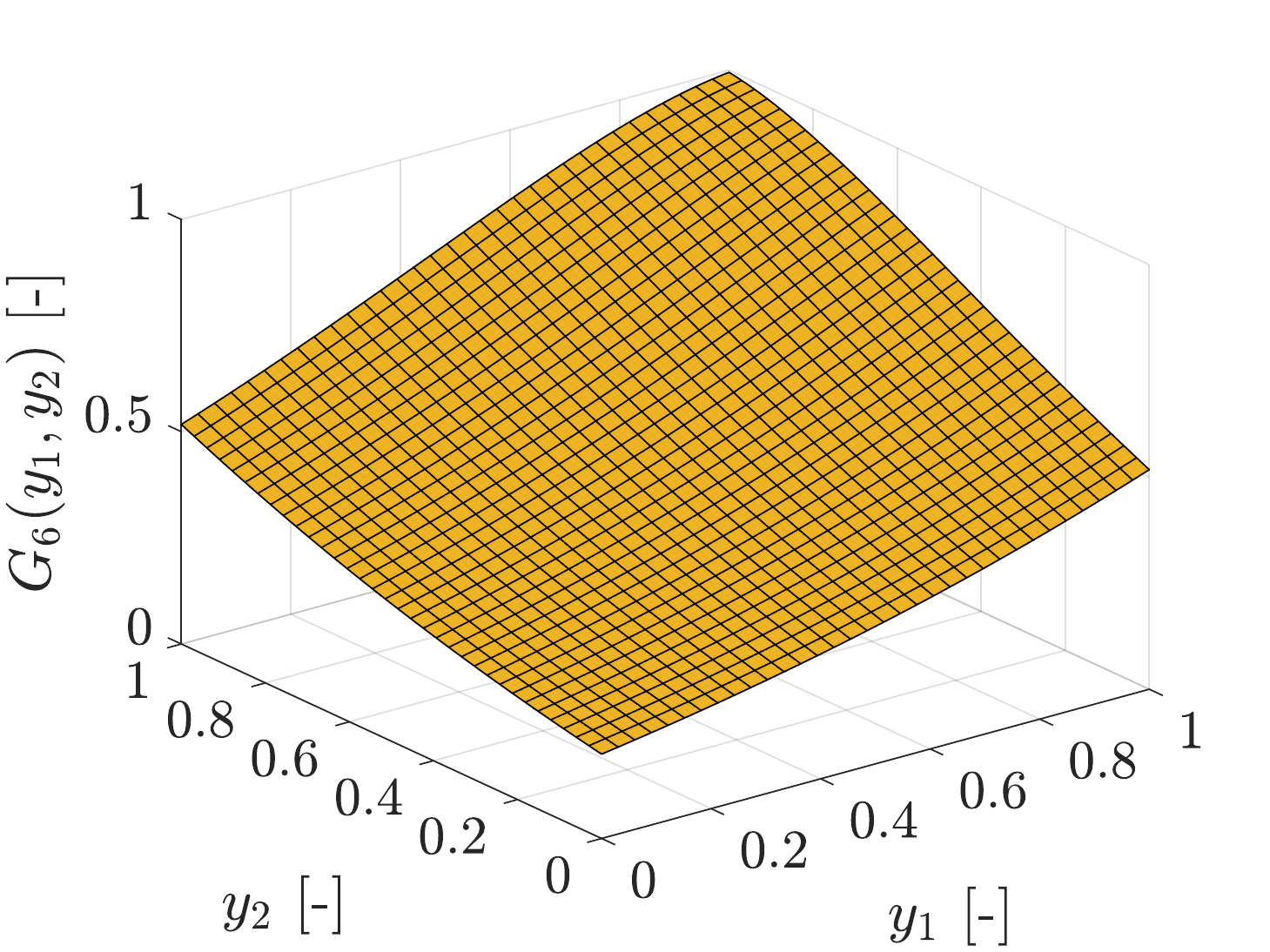}}
\subfigure[Lowest-fidelity function $G_{1}({\bf y})$]{\includegraphics[width=0.3\textwidth]{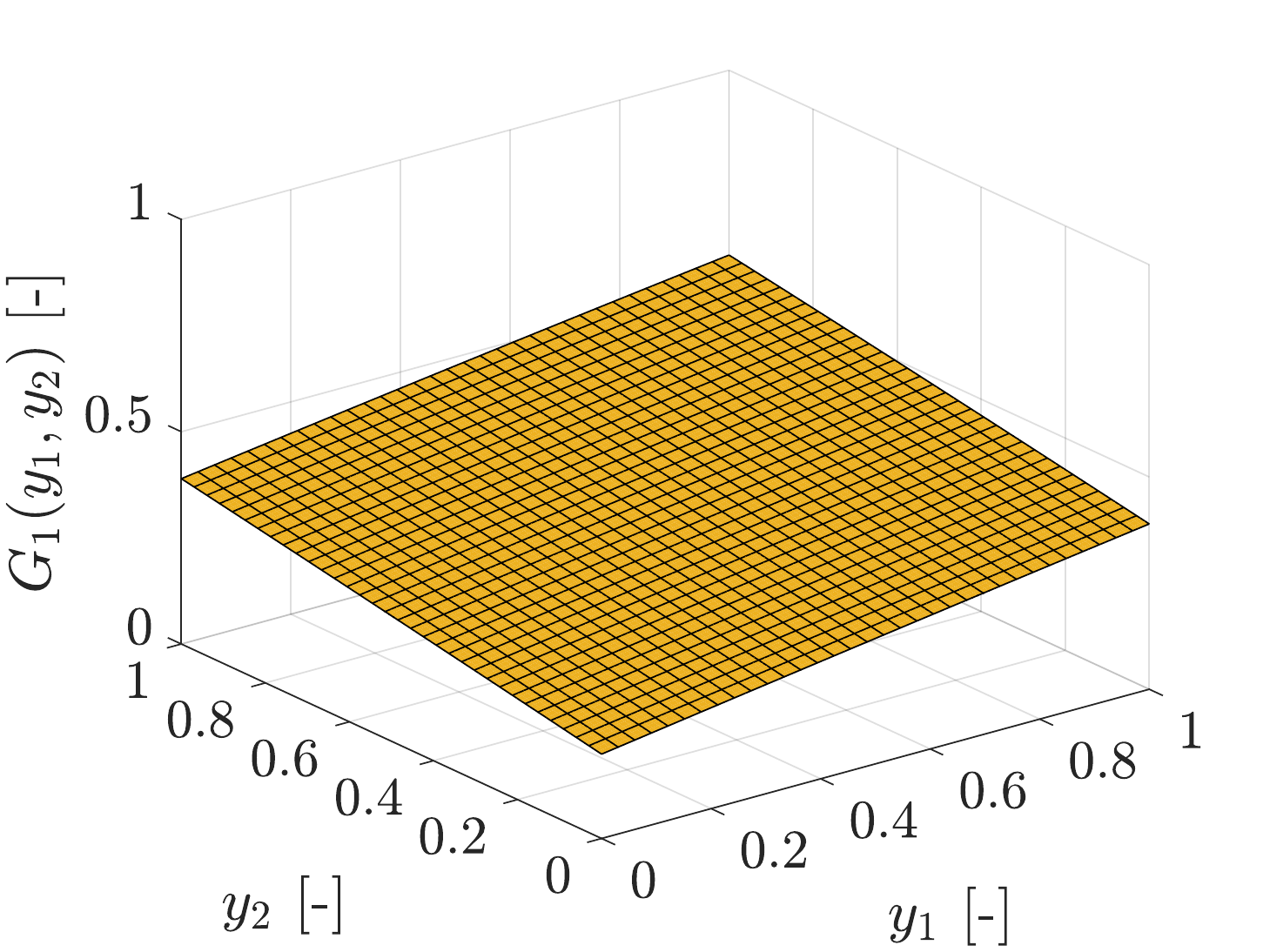}}
\caption{Analytical problem, true function and highest- and lowest-fidelity approximations.}
\label{fig:anprob}
\end{figure}

A normalized computational cost is associated to each evaluation of the $\alpha$-th Taylor expansion $G_{\alpha}({\bf y})$ as
\begin{equation*}\label{eq:AnProbCost}
\mathrm{cost}(\alpha)=8^{\alpha-1}.
\end{equation*}
This choice is {done} to keep the analogy with the RoPax problem {and will become clear} later.

\subsubsection{Numerical results}

\begin{figure}[tp] 
	\centering
	\subfigure[Results for the MISC method with profit $P^{\MISCquad}$.]{
		\includegraphics[width=0.47\textwidth]{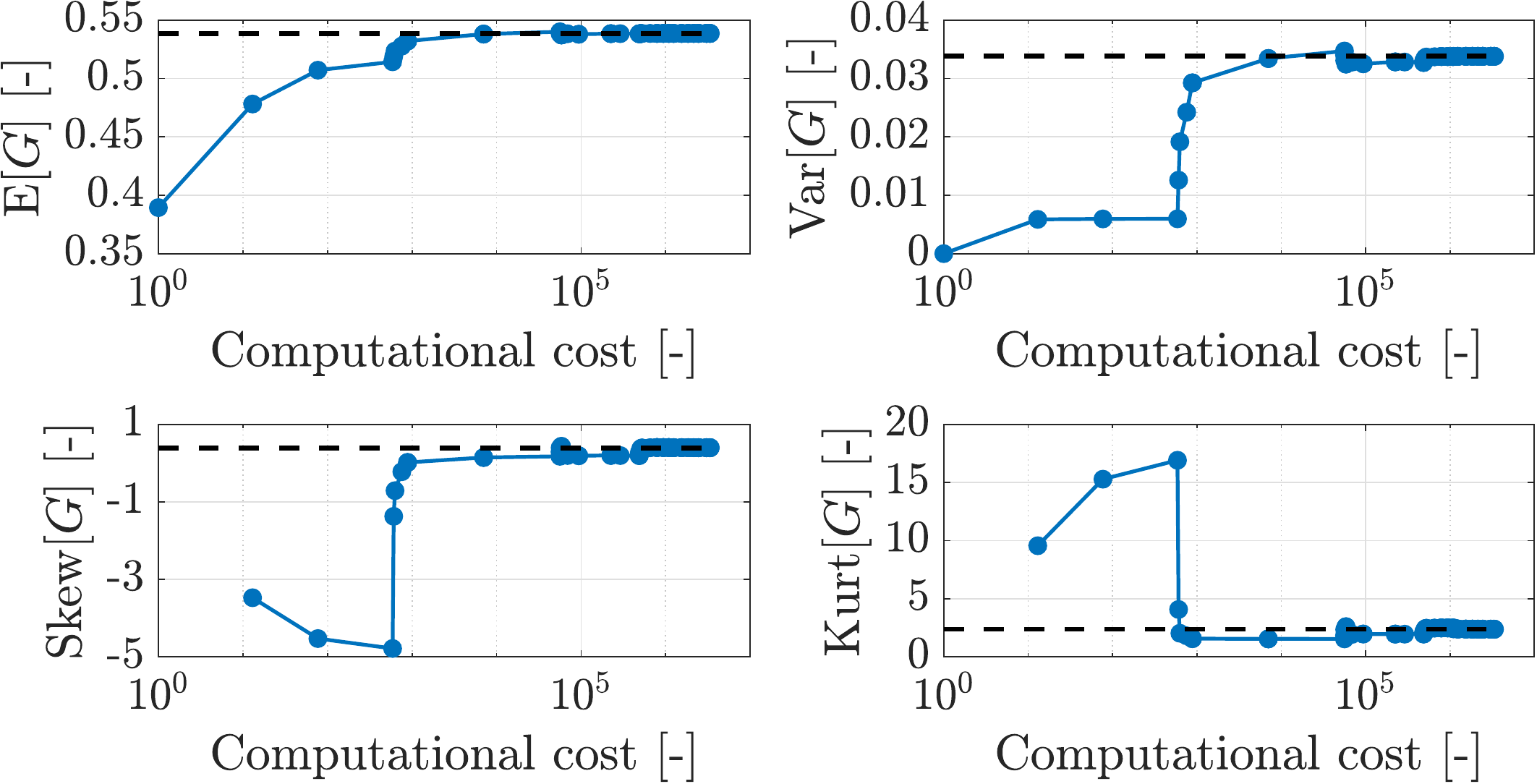}
		\includegraphics[width=0.47\textwidth]{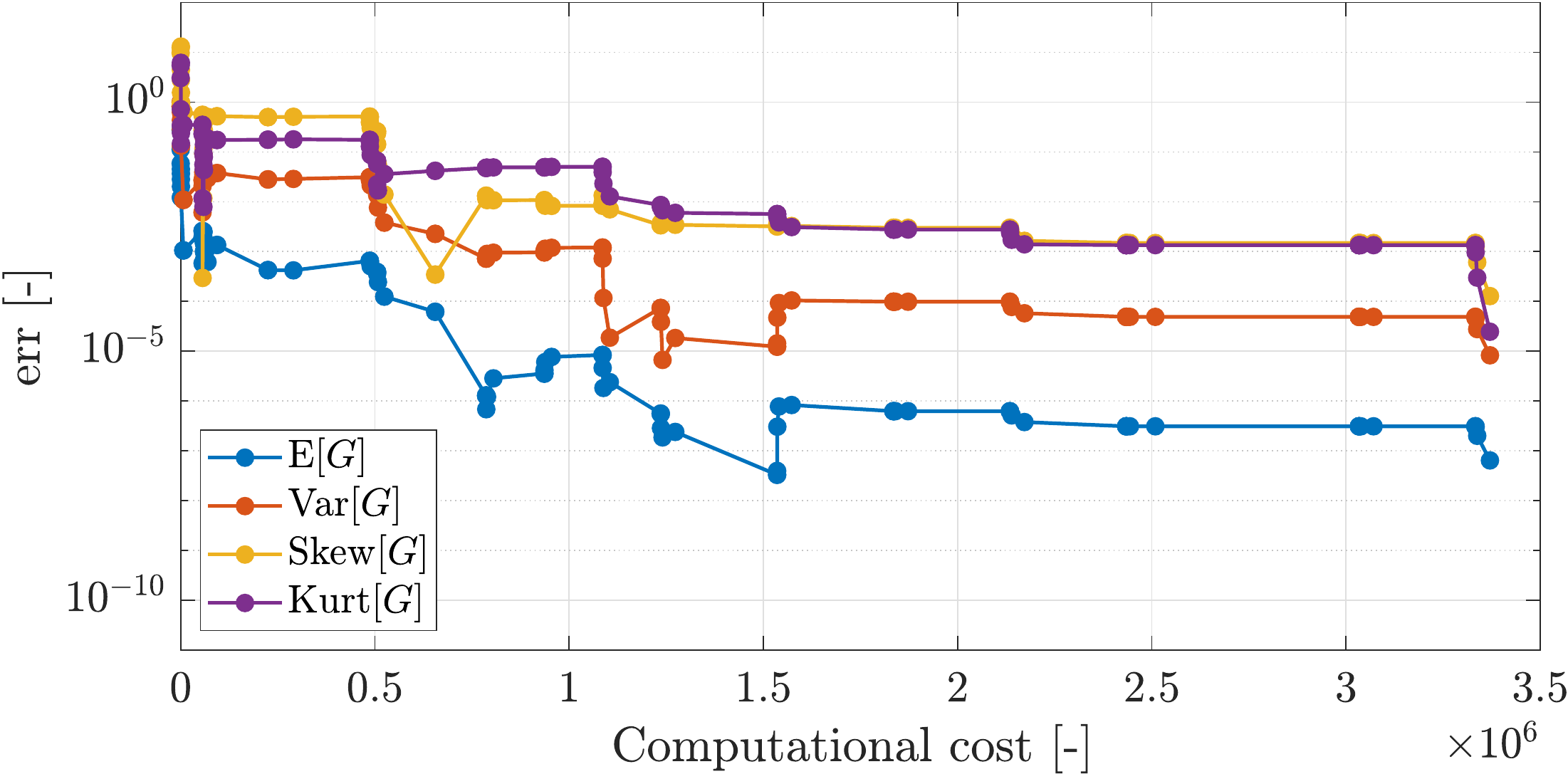}}
	\subfigure[Results for the MISC method with profit $P^{\MISCsurr}$. {For ease of comparison, we overlap the results for MISC-$P^{\MISCquad}$ (gray dashed lines) in the plot on the right}.]{
		\includegraphics[width=0.47\textwidth]{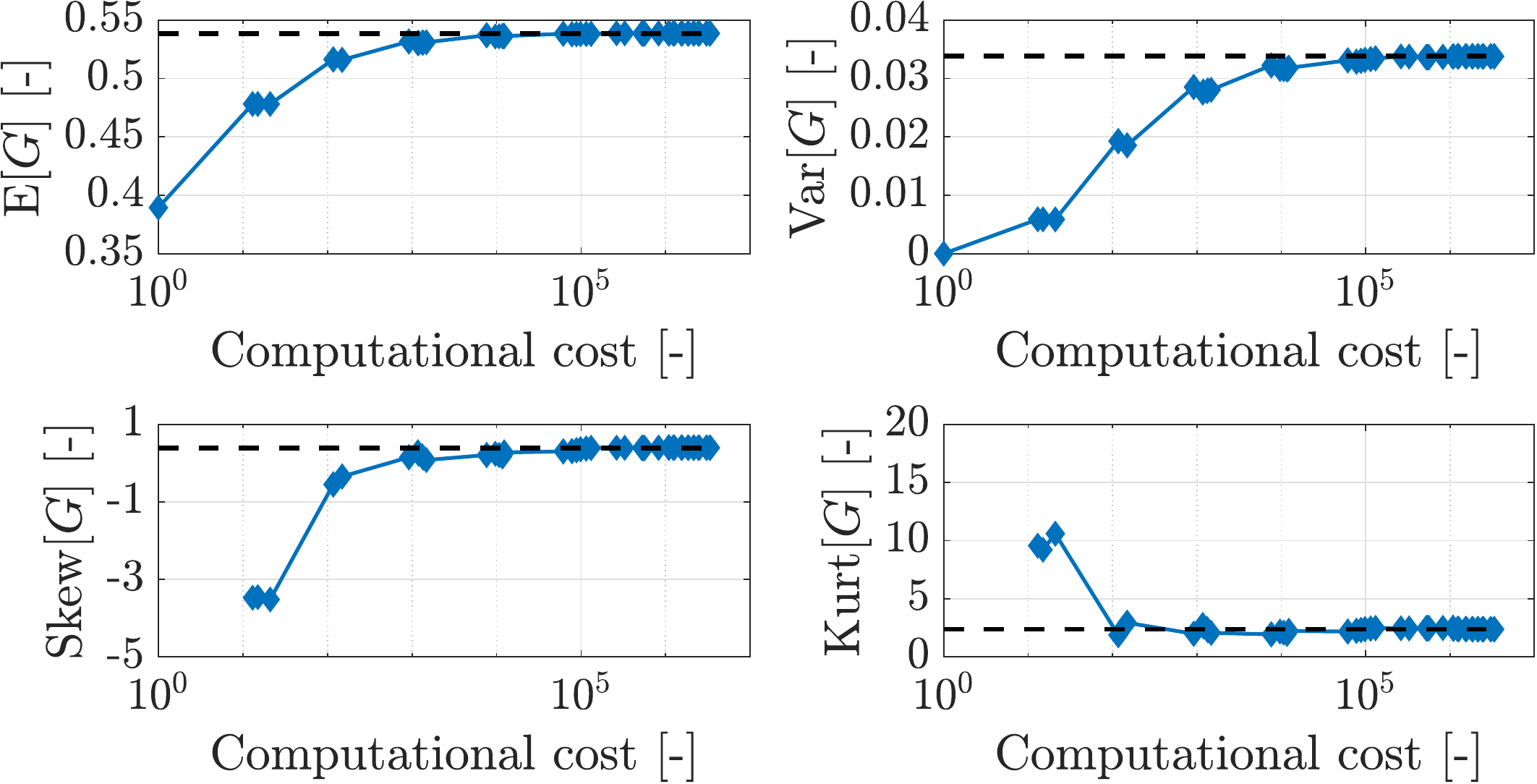}
		\includegraphics[width=0.47\textwidth]{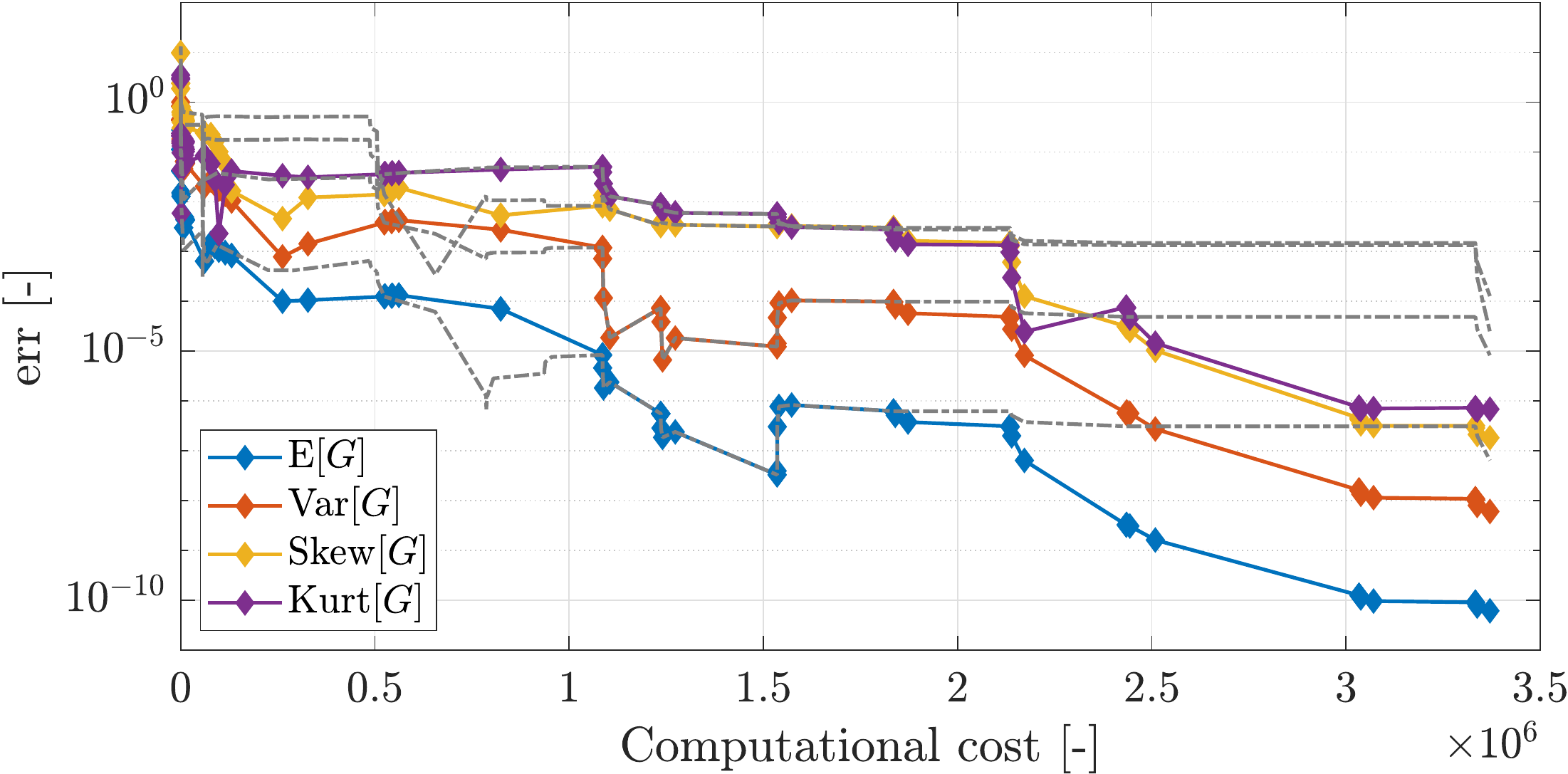}}
        \caption{Analytical problem, results for the MISC method. Left: convergence of the values of the first four centered moments.
        The black dashed line marks the reference value of the moments. 
        Right: relative error of the moments (see Eq.~\eqref{eq:err_moments}).}
	\label{fig:analytic_test_moments_MISC}
\end{figure}

\begin{figure}[tp] 
	\centering
	\includegraphics[width=0.47\textwidth]{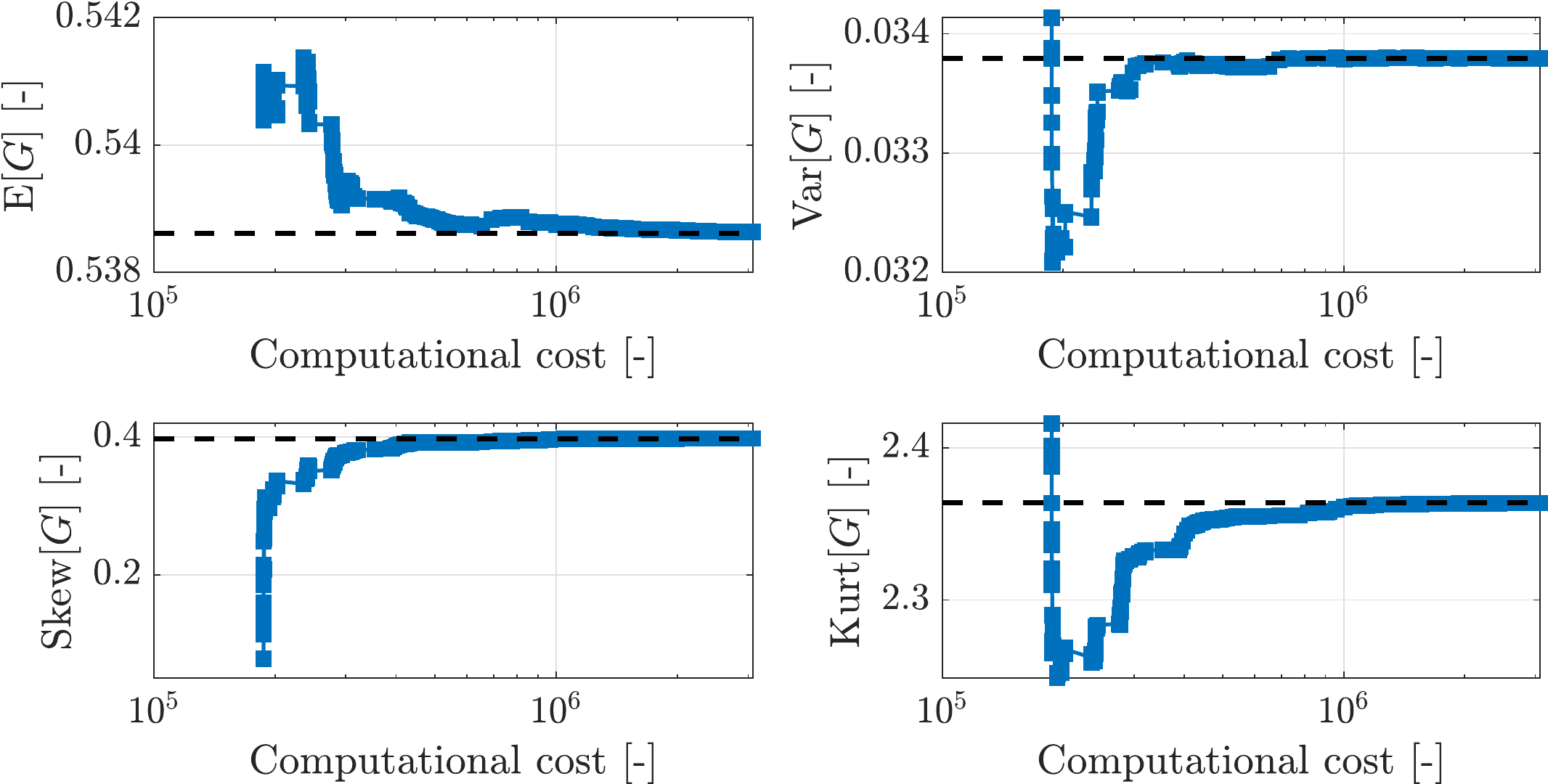}
	\includegraphics[width=0.47\textwidth]{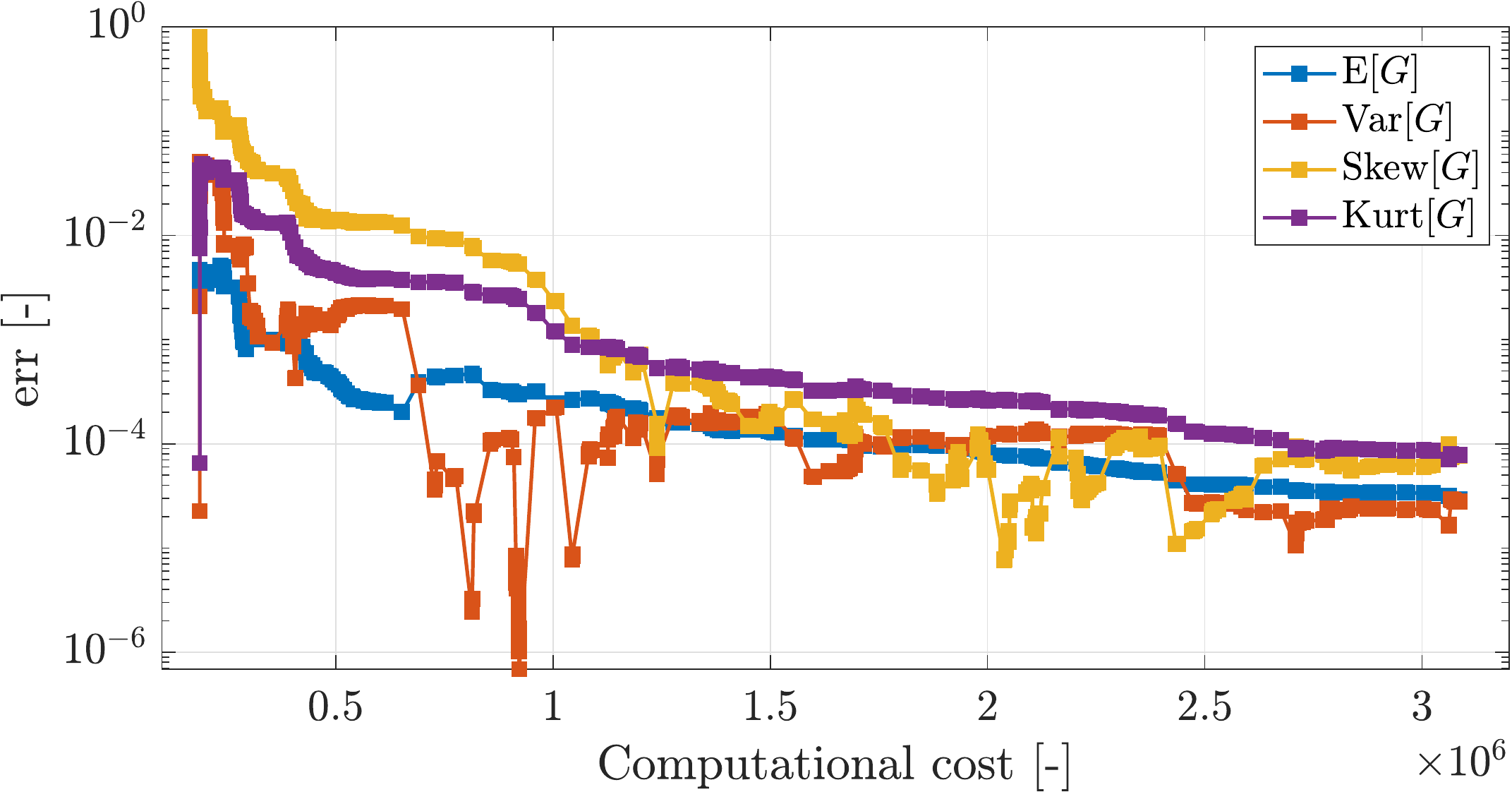}
	\caption{Analytical problem, results for the SRBF method. Left: convergence of the values of the first four centered moments. The black dashed line marks the reference value of the moments. Right: relative error of the moments (see Eq.~\eqref{eq:err_moments}).}
	\label{fig:analytic_test_moments_SRBF}
\end{figure}


We start the discussion with the comparison of the MISC/SRBF estimates of the moments with the reference values. 
The reference values are computed by an accurate sparse-grids quadrature rule with $2^{15}+1$ points where the reference {surrogate model/function} ${G}_{\text{ref}}(\mathbf{y})$ is {the highest-fidelity approximation ${G}_6(\mathbf{y})$}. {The calculations have been done using the Sparse Grids Matlab Kit\footnote{\url{https://sites.google.com/view/sparse-grids-kit}} \cite{MISCsite}}.

In Fig.~\ref{fig:analytic_test_moments_MISC} the {convergence of the} MISC estimates of the first four centered moments and their relative errors as defined in Eq.~\eqref{eq:err_moments} are reported. 
The two variants of MISC {(denoted by MISC-$P^{\MISCquad}$ and MISC-$P^{\MISCsurr}$ in the following)} introduced in Sect.~\ref{sect:misc} are tested, i.e.\ two type of profits {(see Eq.\ \eqref{eq:profits})}, $P^{\MISCquad}$ based on the quadrature error (see Eq.~\eqref{eq:error_contr_quad}) and $P^{\MISCsurr}$ based on the point-wise accuracy of the surrogate model (see Eq.~\eqref{eq:error_contr_Linf}), are considered. 
The results for the case of a quadrature-based profit $P^{\MISCquad}$ are displayed in the first set of plots (see Fig.~\ref{fig:analytic_test_moments_MISC}a left).
All the moments converge to the reference results marked with the black dashed line. 
In Fig.~\ref{fig:analytic_test_moments_MISC}a right one can observe that the error is larger the higher the order of the moment.
Remarkably, even if the adaptivity of the MISC method is driven by the improvement in the first order moment, all the moments are estimated very well. 
The second set of plots (see Fig.~\ref{fig:analytic_test_moments_MISC}b) suggests that also the version of MISC driven by the accuracy of the surrogate model $P^{\MISCsurr}$ is effective in the estimation of the moments. By comparing the two methods, one can observe that the latter one brings better results{, as the error for all the moments is always smaller}.

{
Figure~\ref{fig:analytic_test_moments_SRBF} shows the convergence of the moments for SRBF: differently from MISC, all the moments have a quite similar convergence towards their reference values, with errors all converging within the same order of magnitude. The convergence is almost monotonic for the expected value and for the kurtosis whereas some oscillations can be observed for the variance and the skewness}. Note that for this problem the SRBF method is based on interpolation at the training points (i.e.\ the weights are computed by solving Eq.~\eqref{eq:interp_rbf}), whereas in the following RoPax example the regressive approach (i.e.\ solving Eq. \eqref{eq:lsrbf}) is used, for reasons that will be clear later on. 

%

{Figure~\ref{fig:analytic_test_err_moments} shows that MISC and SRBFs achieve comparable results in terms of the relative errors of the centered moments.}
{Specifically, MISC-$P^{\MISCsurr}$ performs better in the evaluation of all the moments, whereas SRBF performs slightly better than MISC-$P^{\MISCquad}$ in the evaluation of the skewness.}
It is {also} worth noting that SRBFs start with a higher computational cost in comparison with MISC, {due to the fact that the initialization strategy requires to sample all available fidelities.}

\begin{figure} 
	\centering
	\includegraphics[width = 0.94\textwidth]{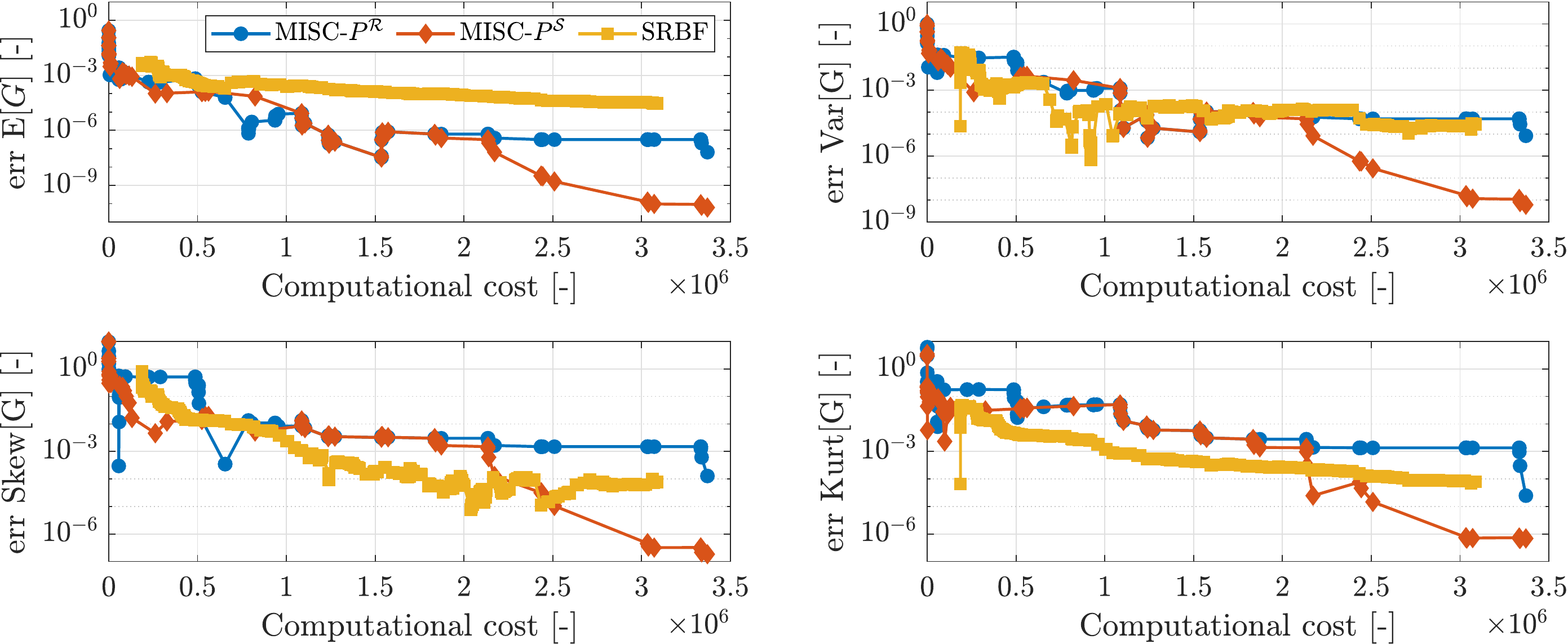}
	\caption{Analytical problem, comparison of MISC and SRBF methods: relative errors of the moments. It is a compact visualization of the results of Fig.~\ref{fig:analytic_test_moments_MISC} and \ref{fig:analytic_test_moments_SRBF}, where the results of MISC are reported only until a computational cost comparable with the one reached by the SRBF method. }
	\label{fig:analytic_test_err_moments}
\end{figure}

The results for the $L_2$ and $L_{\infty}$ norms of the MISC error (see Eq.~\eqref{eq:err_norm}) of Fig.~\ref{fig:analytic_test_norms} are in agreement with the previous findings for the estimation of the moments. 
An improvement of about two orders of magnitude is observed in favor of the surrogate-based method MISC-$P^{\MISCsurr}$ with respect to MISC-$P^{\MISCquad}$ in the final part of the convergence curve.
{When comparing these results with the convergence of the SRBF one can observe that both versions of MISC achieve better results than SRBF.}

\begin{figure}[tp] 
	\centering
		\includegraphics[width = 0.47\textwidth]{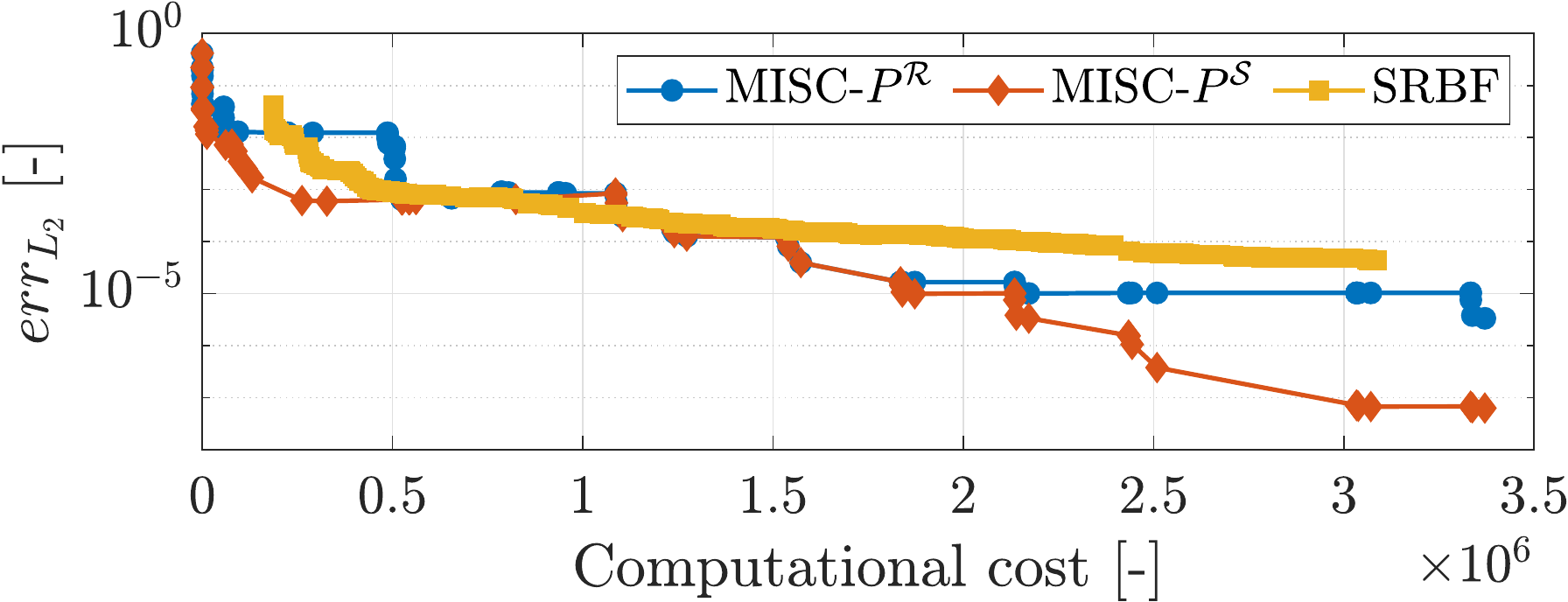}
		\includegraphics[width = 0.47\textwidth]{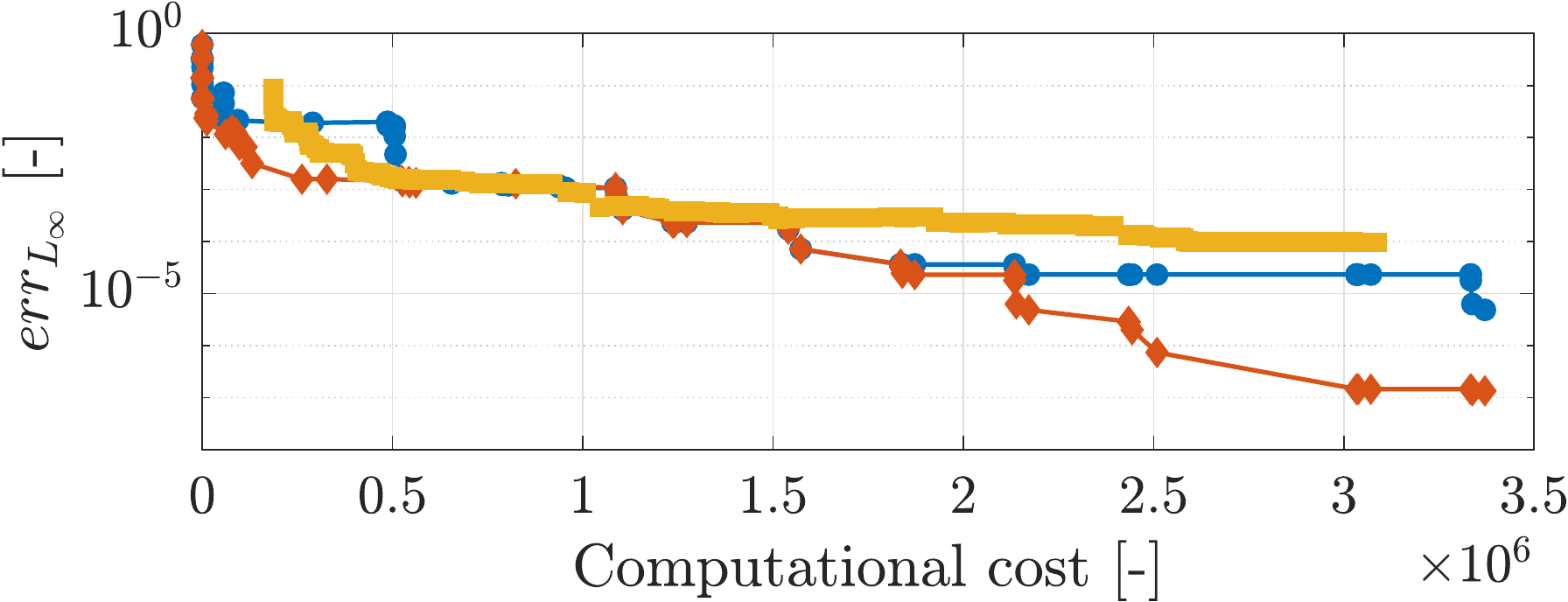}
	\caption{Analytical problem, comparison of MISC and SRBF methods: relative error of the approximation of $G$ in $L_2$ (left) and $L_{\infty}$ norm (right) (see Eq.~\eqref{eq:err_norm}). }
	\label{fig:analytic_test_norms}
\end{figure}


\begin{figure}[tp]
	\centering
	\subfigure[PDFs at the final iteration (obtained by Matlab's {\it ksdensity} enforcing positive support since $G$ takes positive values only)]{\includegraphics[width = 0.47\textwidth]{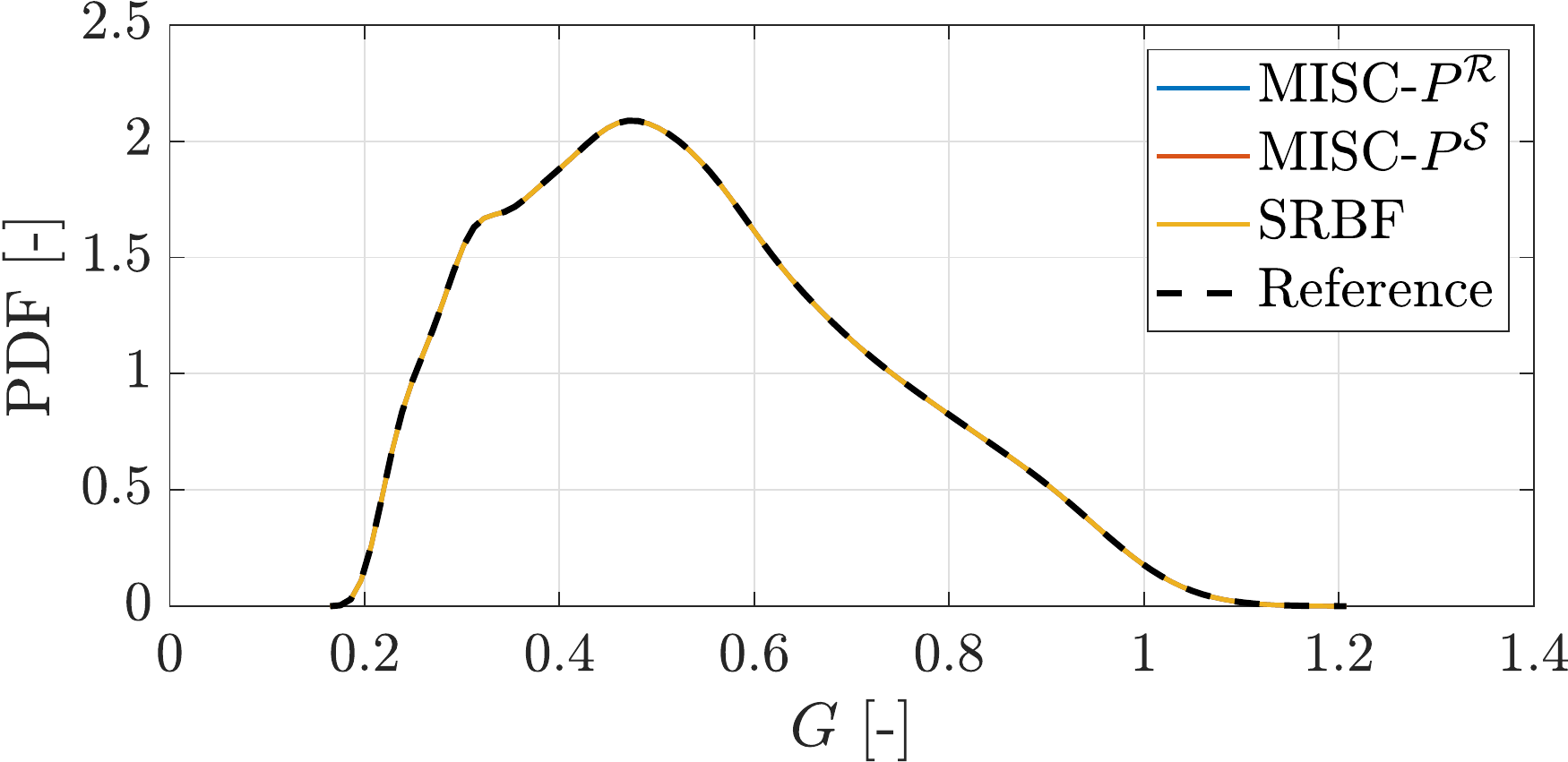}} 
	\subfigure[KS test statistic (see Eq.~\eqref{eq:KS_stat})]{\includegraphics[width = 0.47\textwidth]{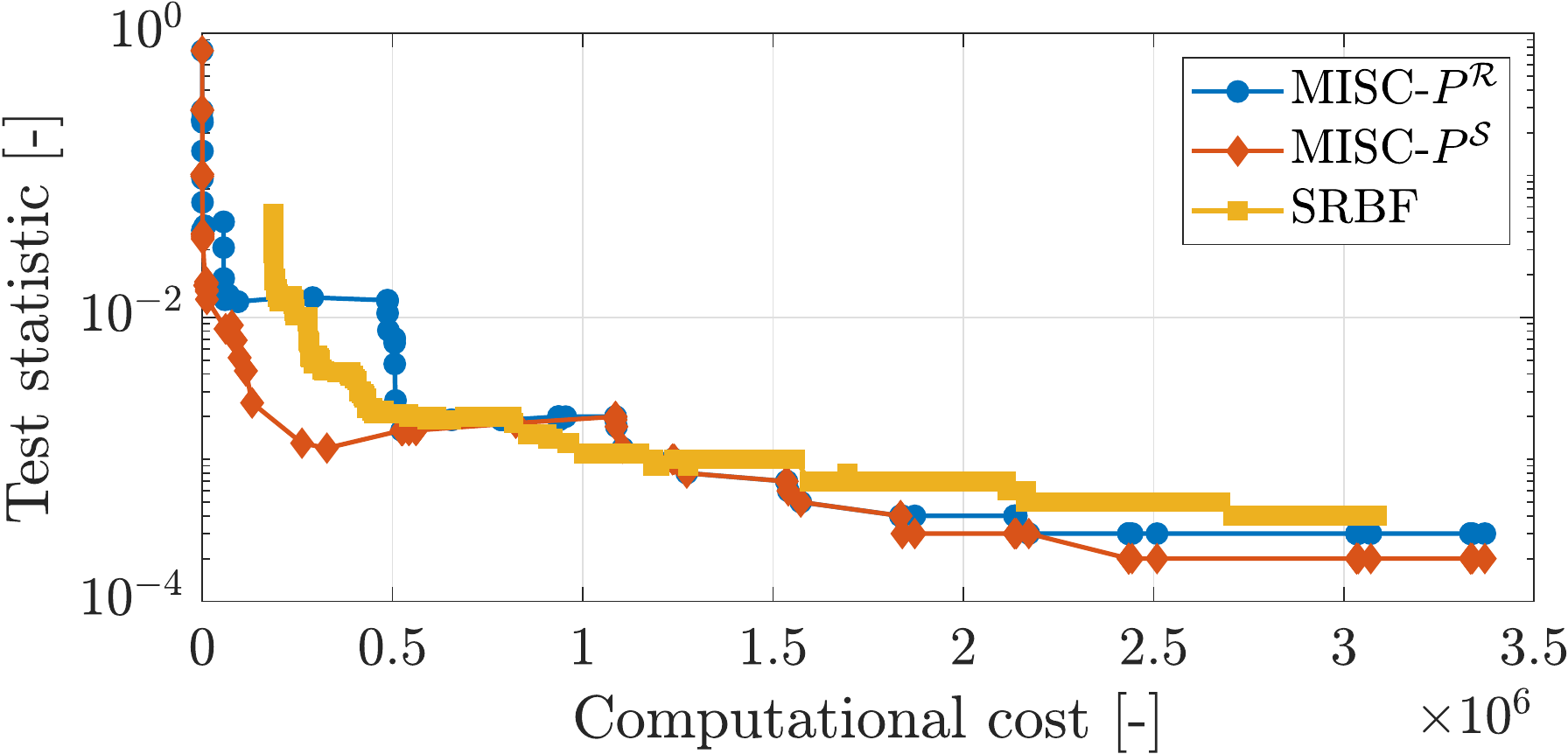}}
	\caption{Analytical problem, comparison of MISC and SRBF methods.}
	\label{fig:analytic_test_pdf}
\end{figure}


The comparison of the PDFs given in Fig.~\ref{fig:analytic_test_pdf}a shows a very good agreement of the MISC and SRBF results with the reference ones.
In Fig.~\ref{fig:analytic_test_pdf}b the results of the KS test statistics (cf. Eq.~\eqref{eq:KS_stat}) are reported: {both versions of MISC show {a slightly} better convergence of the test statistic. }


Finally, it is worth looking at the sampling performed by the two methods. In Fig.~\ref{fig:analytic_test_sampling} the samples selected by MISC-$P^{\MISCsurr}$ and SRBF are displayed. In the first case (see Fig.~\ref{fig:analytic_test_sampling}a) the samples are well distributed over the domain in a symmetric way, with no preferential directions.
The SRBF sampling is instead slightly more clustered in the regions where the high-fidelity function shows a larger curvature (cf. Fig.~\ref{fig:anprob}b and \ref{fig:analytic_test_sampling}b). 
The sampling performed by MISC-$P^{\MISCquad}$ is not shown for brevity, as it is very similar to the one for MISC-$P^{\MISCsurr}$.

\begin{figure}[tp]
	\centering
	\subfigure[MISC with profit $P^{\MISCsurr}$]{\includegraphics[scale = 0.45]{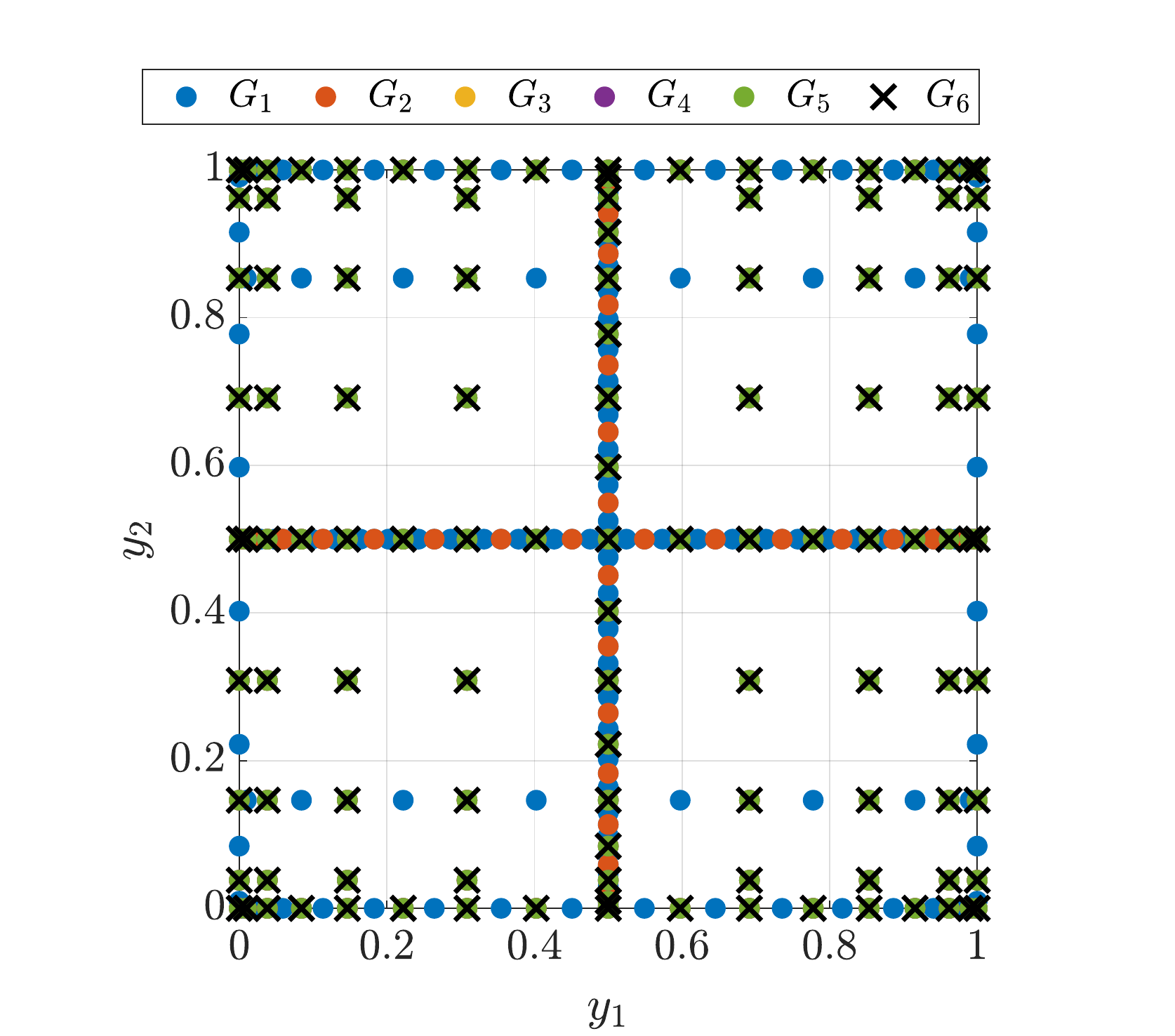}}
	\subfigure[SRBF]{\includegraphics[scale = 0.45]{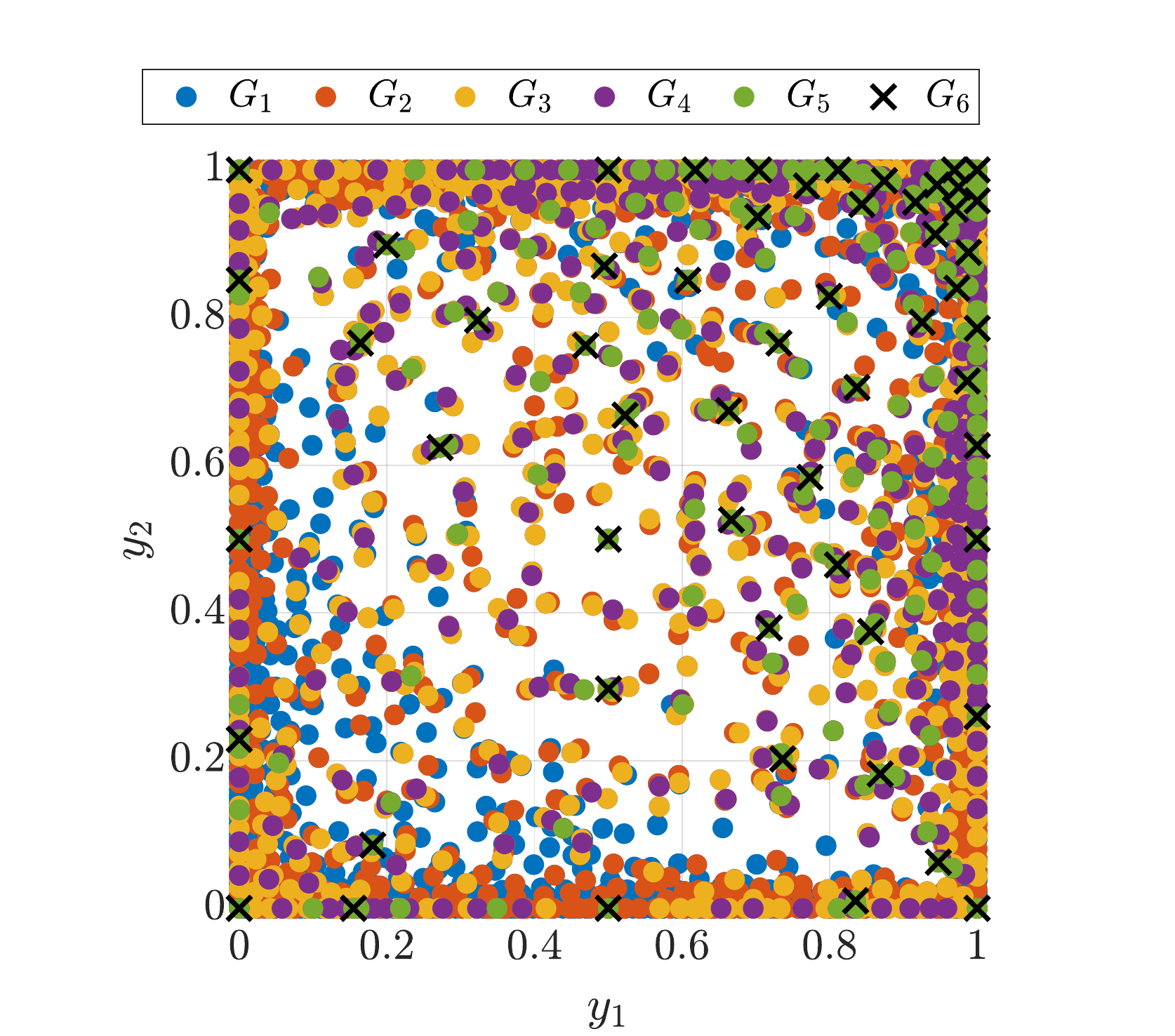}}
	\caption{Analytical problem, results for the MISC and SRBF method. Sampling points at the last iteration.}
	\label{fig:analytic_test_sampling}
\end{figure}

\subsection{RoPax resistance problem}\label{sect:problem_description_ropax}
\subsubsection{Formulation and CFD method}\label{sect:RoPax_problem}

The main problem addressed in this {work} is the forward UQ analysis of the model-scale resistance ($R_T$) of a RoPax ferry in straight ahead advancement, subject to two operational uncertainties $\mathbf{y}=[U,T]$, namely the advancement speed ($U$) and the draught ($T$), uniformly distributed within the ranges in Tab.~\ref{tab:opcon}. {The choice of using two operational parameters allows visual investigation of the results while preserving all the main difficulties that arise when solving parametric approximation problems. Furthermore, the advancement speed and draught are two operational parameters with significant practical implications: for instance Froude and Reynolds numbers vary with the advancement speed; allowed payload and block coefficient vary with the draught, etc.}
%
%

\begin{figure}[tp]
\centering
\includegraphics[width=0.80\textwidth]{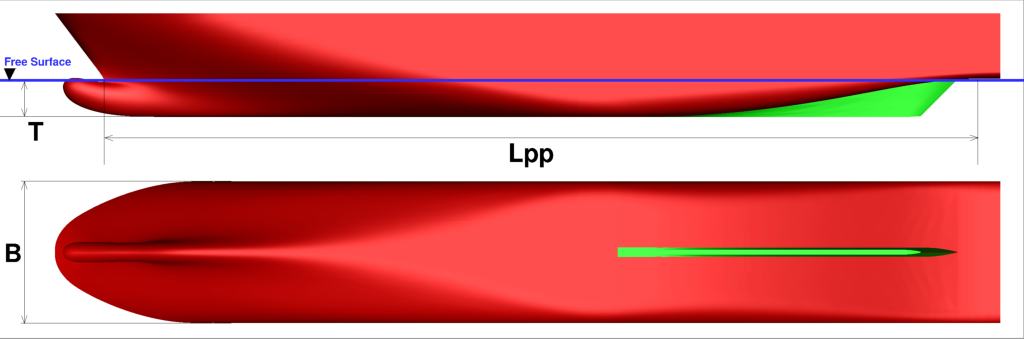}
\caption{RoPax ferry, hull form. Free surface level is reported for the nominal draught.}
\label{fig:shape}
\end{figure}

The RoPax ferry is characterized by a length between perpendicular at nominal draught ($L_{\rm PP}$) of $162.85$~m and a block coefficient $C_B=0.5677$ (see Fig.~\ref{fig:shape}). The parametric geometry of the RoPax is produced with the computer-aided design environment integrated in the CAESES software, developed by FRIENDSHIP SYSTEMS AG, and made available in the framework of the 
H2020 EU Project Holiship\footnote{\url{www.holiship.eu}}. 
The analysis is performed at model scale with a scale factor equal to $27.14$. The main dimensions and the operative conditions are summarized in Tab.~\ref{tab:opcon}. The advancement speed ranges from $12$ to $26$ knots at full scale and the draught variation is $\pm 10\%$ of the nominal draught, which corresponds to a variation of about $\pm 15\%$ of the nominal displacement. The corresponding range in Froude number ${\rm Fr}={U}/{\sqrt{g L_{\rm PP}}}$ is $[0.154,0.335]$, whereas the variation in Reynolds number (at model scale) is ${\rm Re}=\rho U L_{\rm PP}/\mu=U L_{\rm PP}/\nu \in [6.423\cdot10^6,1.392\cdot10^7]$, where $\rho=998.2$ kg/m$^3$ is the water density, $\nu=\mu/\rho=1.105\cdot10^{-6}$ m$^2$/s the kinematic viscosity and $g=9.81$ m/s$^2$ the gravitational acceleration.

\begin{table}[tp]
	\begin{center}
		\caption{Main geometrical details and operative conditions of the RoPax ferry (model scale $1:27.14$).}
		\label{tab:opcon}
		\begin{tabular}{lcccc}
			\toprule
			\textbf{Description} & \textbf{Symbol} & \textbf{Full Scale} & \textbf{Model Scale} & \textbf{Unit}   \\
			\midrule
			Length between perpendiculars & $L_{\rm PP}$ & $162.85$          & $6.0$              & m     \\
			Beam                          & $B$          & $29.84$           & $1.0993$             & m    \\
			Block coefficient             & $C_B$        & $0.5677$          & $0.5677$           & --   \\
			Nominal displacement          & $\nabla$     & $19584.04$        & $0.9996$           & m$^3$ \\
			Nominal draught               & $T_n$        & $7.10$            & $0.261660$ & m \\
			Draught range                 & $T$  & $[ 6.391, 7.812]$ & $[0.2355, 0.2878]$ & m \\
			Speed range                   & $U$  & $[6.173, 13.376]$ & $[1.185, 2.567]$ & m/s   \\
			Froude range                  & Fr & $[0.154, 0.335]$ & $[0.154, 0.335]$ & --   \\
			Reynolds range                & Re & $[9.081\cdot10^8, 1.968\cdot10^9]$ & $[6.423\cdot 10^6, 1.392\cdot 10^7]$ & --   \\
			\bottomrule
		\end{tabular}                                           
	\end{center}
\end{table}

The hydrodynamics performance of the RoPax ferry is assessed by the RANS code $\chi$navis developed at CNR-INM. The main features of the solver are summarized here; for more details the interested reader is referred to \cite{dimascio2007,dimascio2009,broglia2014,broglia2018} and references therein. $\chi$navis is based on a finite volume discretization of the RANS equations, with variables collocated at the cell centers. Turbulent stresses are related to the mean velocity gradients by the Boussinesq hypothesis; the turbulent viscosity is estimated by the Spalart-Allmaras turbulence model~\cite{spalall}. Wall-functions are not adopted, therefore the wall distance $y^+=1$ is ensured on the no-slip wall. Free-surface effects are taken into account by a reliable single-phase level-set approach.

The computational domain extends to 2$L_{\rm PP}$ in front of the hull, 3$L_{\rm PP}$ behind, and 1.5$L_{\rm PP}$ {sideway}; a depth of 2$L_{\rm PP}$ is imposed (see Fig.~\ref{fig:bc}a). On the solid walls (in red in the figure), the velocity is set equal to zero, whereas zero normal gradient is enforced on the pressure field; at the (fictitious) inflow boundary (in blue in Fig.~\ref{fig:bc}a), the velocity is set to the undisturbed flow value and the pressure is extrapolated from inside; the dynamic pressure is set to zero at the outflow (in yellow), whereas the velocity is extrapolated from inner points. On the top boundary, which remains always in the air region, fluid dynamic quantities are extrapolated from inside (in purple).
Taking advantage of the symmetry of the flow relative to the $y=0$ plane, computations are performed for half ship only, {and the} usual symmetry boundary conditions are enforced on the symmetry longitudinal plane (in green).

\begin{figure}[tp]
\centering
\subfigure[Boundary conditions]{\includegraphics[width=0.45\textwidth]{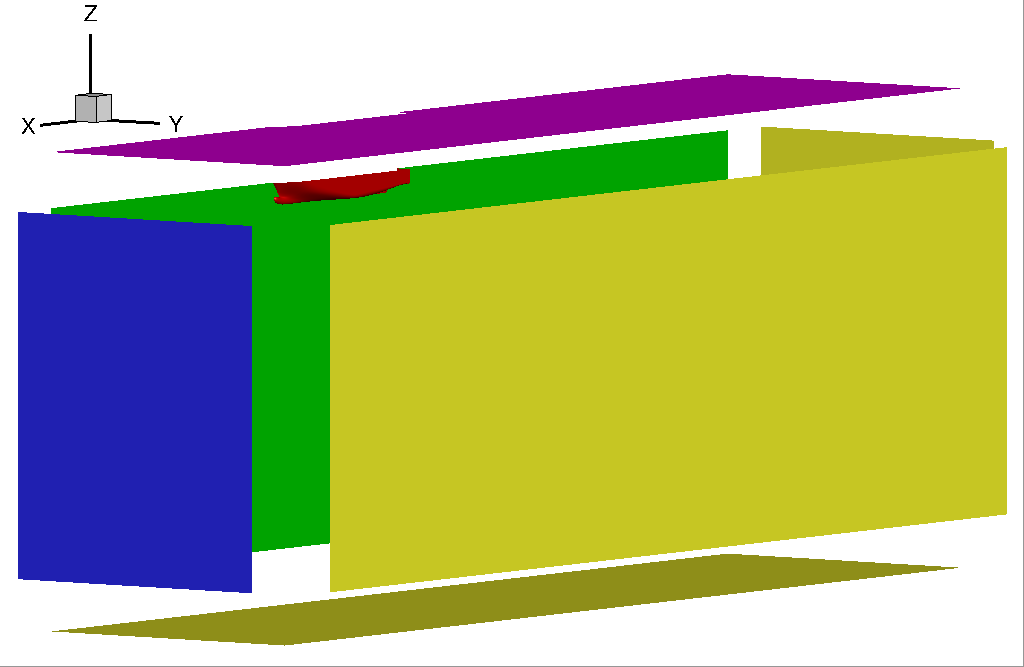}}
\subfigure[Computational grid] {\includegraphics[width=0.45\textwidth]{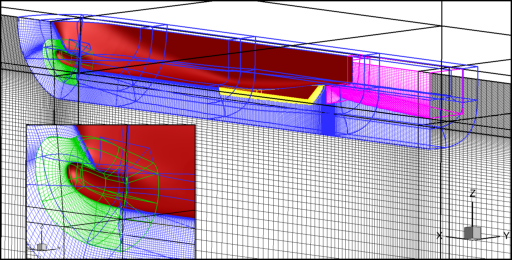}}
\caption{RoPax ferry, numerical setup. Boundary conditions and computational grid.}
\label{fig:bc}
\end{figure}
\begin{figure}[tp]
\centering
\includegraphics[width=0.8\textwidth]{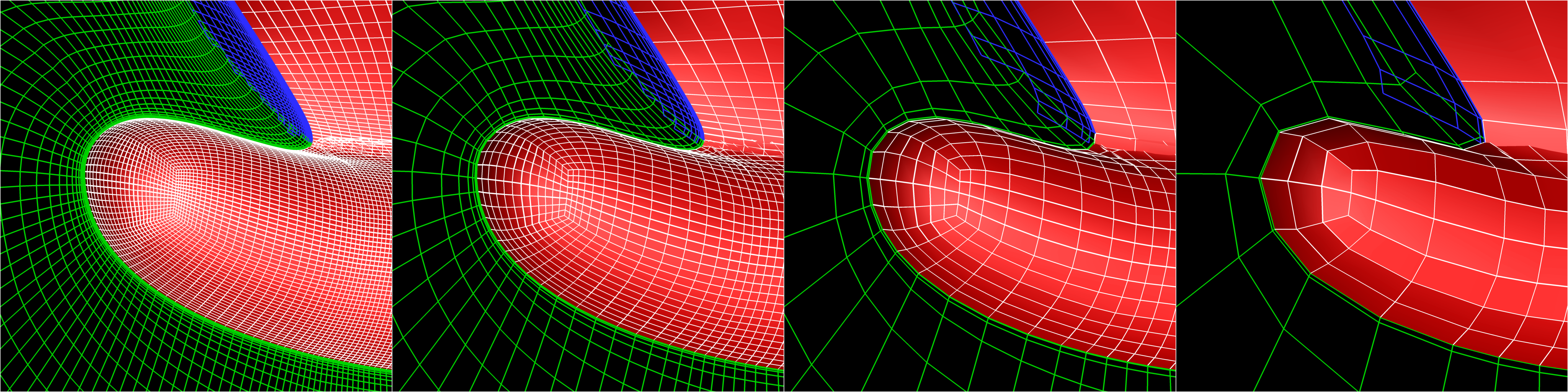}
\caption{RoPax grids, detail of the bow region, left to right: $\mathcal{M}_4$, $\mathcal{M}_3$, $\mathcal{M}_2$, $\mathcal{M}_1$.}
\label{fig:grids}
\end{figure}
\begin{figure}[tp]
\centering
\includegraphics[width=0.8\textwidth]{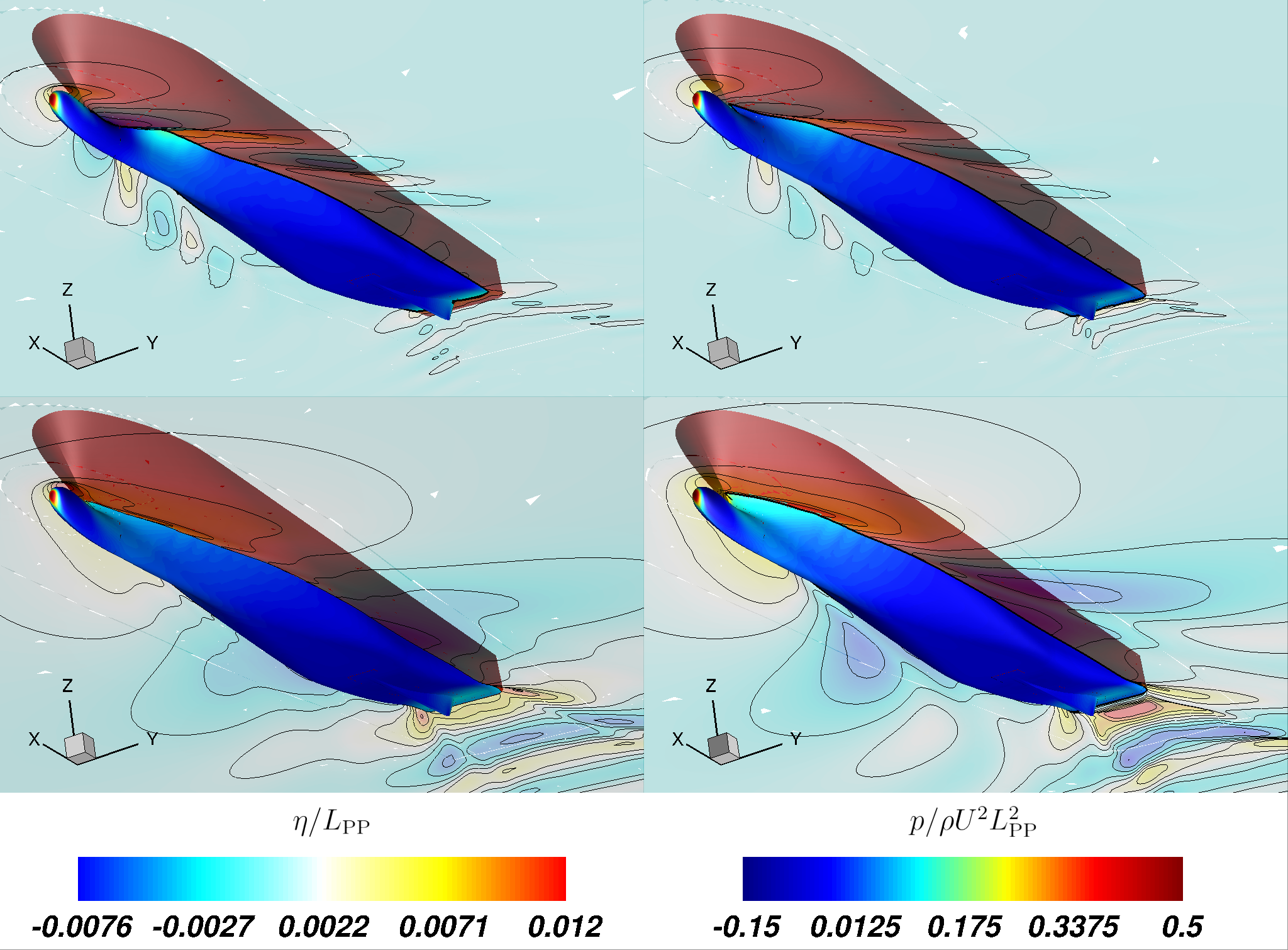}
\caption{RoPax ferry, $\mathcal{M}_4$ CFD results in terms of non-dimensional wave pattern {(left)} and surface pressure {(right)} for:
	${\rm Fr}=0.193$, $T=3.9249\cdot10^{-2}L_{\rm PP}$ and $T=4.7971\cdot10^{-2}L_{\rm PP}$, top row left and right;
	${\rm Fr}=0.335$, $T=3.9249\cdot10^{-2}L_{\rm PP}$ and $T=4.7971\cdot10^{-2}L_{\rm PP}$, bottom row left and right.}
	\label{fig:RoPax1}
\end{figure}

\begin{figure}[tp]
\centering
\includegraphics[width=0.8\textwidth]{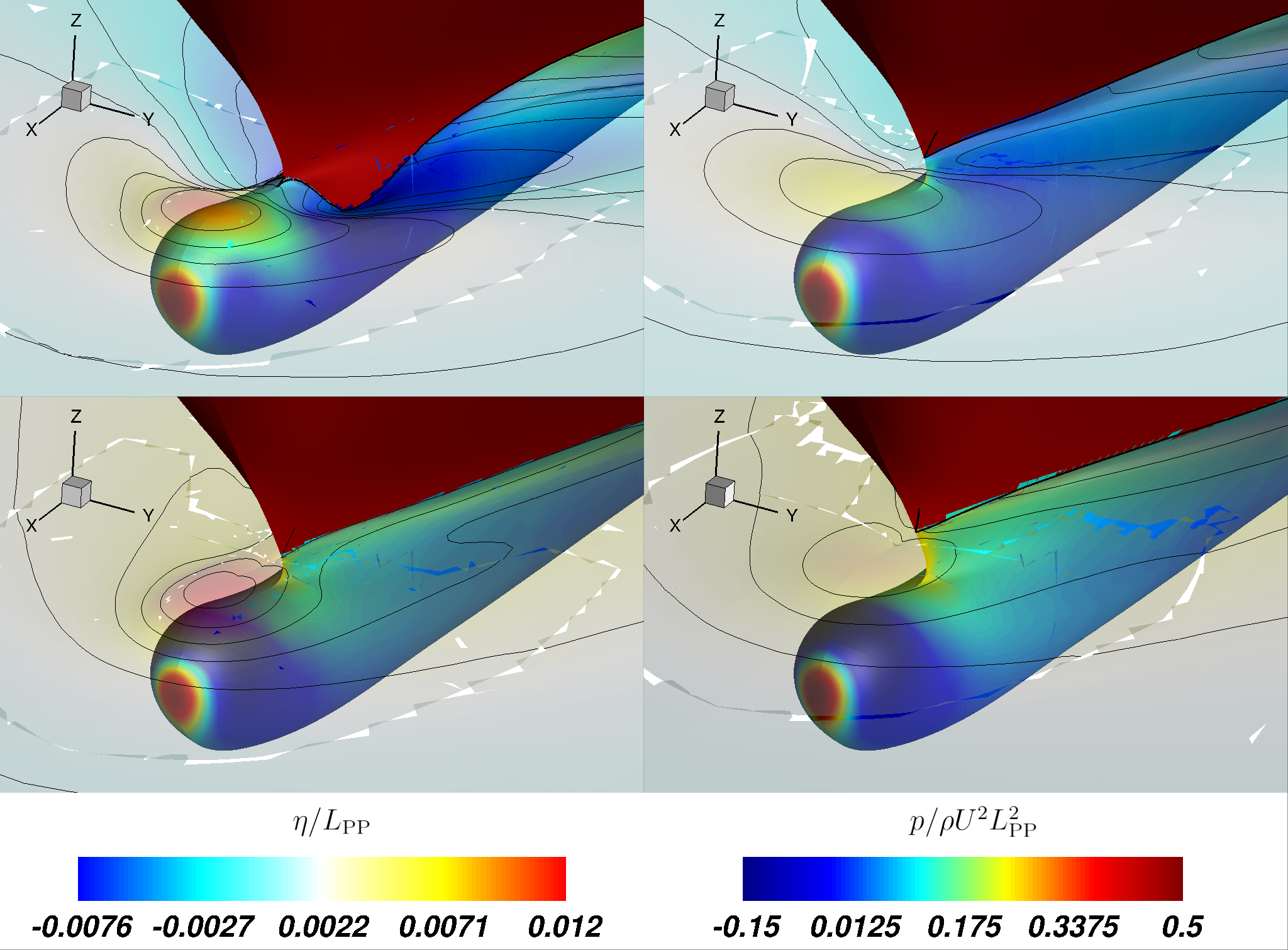}
\caption{RoPax ferry, enlarged view of the bow region as in Fig.~\ref{fig:RoPax1}.}
\label{fig:RoPax2}
\end{figure}

The computational grid is composed by 60 adjacent and partially overlapped blocks; Fig.~\ref{fig:bc}b shows a particular of the block structures in the region around the ship hull and the computational grid on the symmetry plane. Taking the advantage of a Chimera overlapping approach, the grids around the skeg and around the bow are generated separately from the grid around the hull; a background Cartesian grid is then built and the whole grid is assembled by means of an in-house overlapping grid pre-processor. The final grid counts for a total of about 5.5M control volumes for half the domain. The numerical solutions are computed by means of a full multi-grid--full approximation scheme (FMG--FAS), with four grid levels (from coarser to finer: $\mathcal{M}_1$, $\mathcal{M}_2$, $\mathcal{M}_3$, and $\mathcal{M}_4$), each obtained from the next finer grid with a coarsening ratio equal to 2, along each curvilinear direction. In the FMG--FAS approximation procedure, the solution is {first} computed on the coarsest grid level {and then} approximated on the next finer grid 
by exploiting all the coarser grid levels available with a V-Cycle. The process is repeated up to the finest grid level.
For the present UQ problem the number of grid volumes is 5.5M for $\mathcal{M}_4$, 699K for $\mathcal{M}_3$, 87K for $\mathcal{M}_2$, and 11K for $\mathcal{M}_1$.
To provide an idea about the different grid resolutions between the grid levels, Fig.~\ref{fig:grids} shows a particular of the grid
at the bow region for $\mathcal{M}_4$, $\mathcal{M}_3$, $\mathcal{M}_2$ and $\mathcal{M}_1$ grids. 

{Since the grids are obtained by a dyadic derefinement, the following normalized computational costs can be assigned to each grid:}
\begin{equation}
  \label{eq:cost}
  \text{cost}({\alpha}) = 8^{\alpha-1}
\end{equation}
with $\alpha=1,\dots,4$.
In the FMG-FAS scheme the computation on the $\alpha$-th grid level involves computations on all the coarser grids 
. However, with the estimation in Eq.~\eqref{eq:cost}, only the cost of the highest-fidelity level samples is taken into account, i.e. the computations on the coarser grids are considered negligible.

Fig.~\ref{fig:RoPax1} shows an overview of the numerical solutions obtained for different conditions in terms of wave pattern and pressure on the hull surface; wave height (as elevation with respect to the unperturbed level) and surface pressure are reported in non-dimensional values, making the height non dimensional with $L_{\rm PP}$ and the pressure with $\rho U^2 L_{\rm PP}^2$, as usual. 
A clear, and obvious, Froude number dependency is seen for the wave patterns; at the lower speed shown, the free surface is weakly perturbed ({note that} the same color range has been used for all the panels), whereas, at higher Froude, a clear Kelvin pattern is seen. Also, at higher speed, the formation of a well-defined transom wave system is observed, including the presence of the classical rooster tail.
It is also worth to observe the influence of the draught on the wave system; in particular at the lowest speed and smaller draught reported, the rear part of the bulbous is partially dry (better seen in the enlarged views presented in Fig.~\ref{fig:RoPax2}).
The region of very low pressure caused by the flow acceleration around the bow is obviously the cause.
For all cases, the high pressure in the stagnation point at the bow prevents the bow to be outside the water, as it is at the rest conditions at least for the nominal and the smaller draughts (see Fig.~\ref{fig:shape}).
At the higher speed, the larger draught condition causes a stronger rooster tail system at the stern, with higher crest and trough.

Figure~\ref{fig:RoPax_noise} shows the complete FMG--FAS cycle for the minimum and the maximum Froude numbers. The final evaluation of the $R_T$ for each grid is performed averaging the $R_T$ value among the last 100 iterations of the cycle. These are highlighted by the gray areas. 
Even if a general second order convergence has been verified (not shown here for the sake of conciseness), it is evident that, although the FMG--FAS switches to a finer grid when the solver residual are lower than the defined convergence threshold (not shown here), the value of $R_T$ is clearly not converged yet, at least for the coarsest grid level. This has been observed mostly for low Froude numbers. Therefore, the final value of $R_T$ can significantly deviate for simulations performed on the same grid but with slightly different conditions, thus producing evaluations affected by numerical noise. This will have an impact on the following UQ analysis. 

\begin{figure}[tp]
	\centering
	\subfigure[Complete FMG--FAS cycle]{\includegraphics[width=0.6\textwidth]{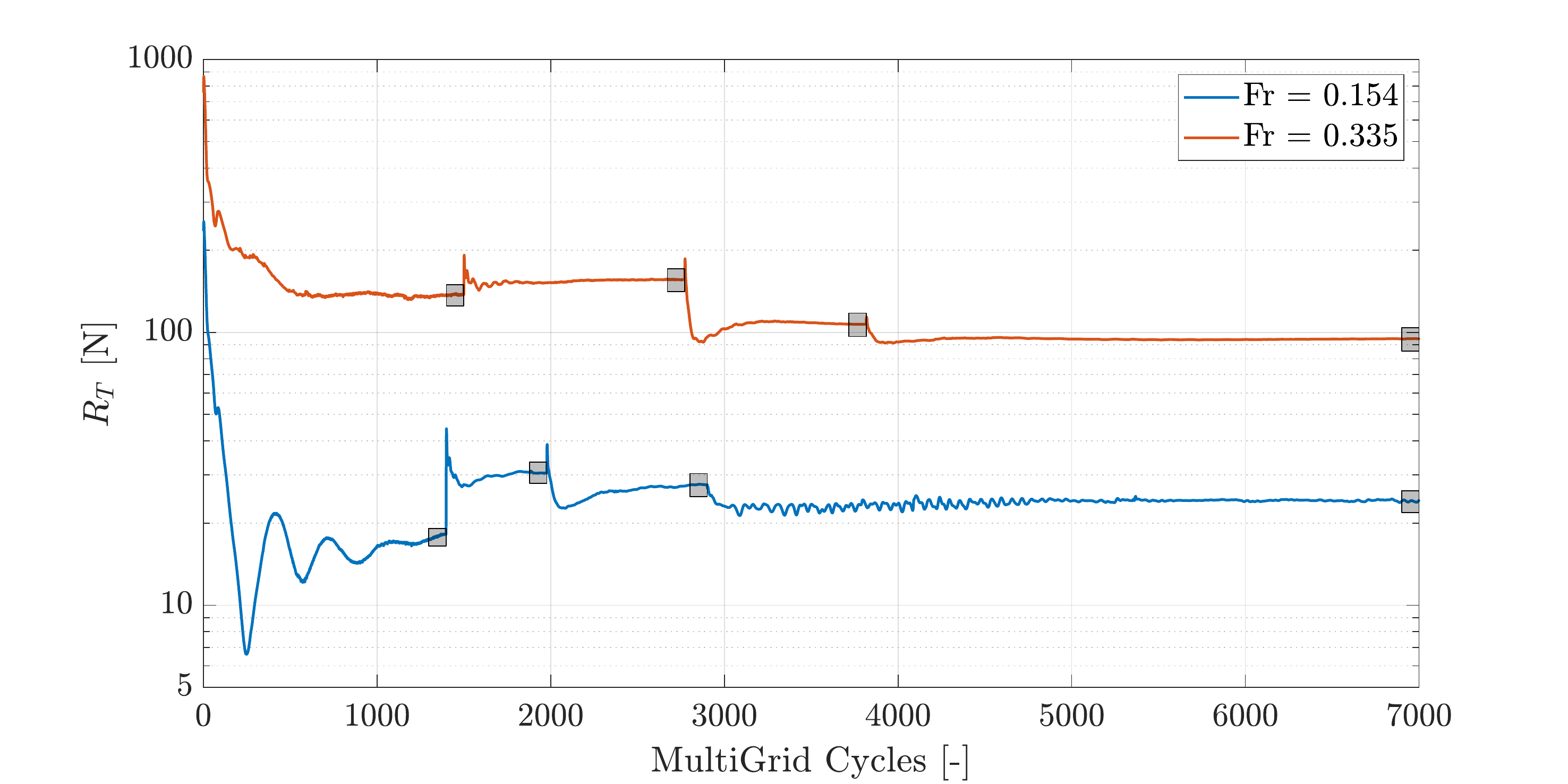}}\\
	\subfigure[Detail of the shift from $\mathcal{M}_1$ to $\mathcal{M}_2$]{\includegraphics[width=0.24\textwidth]{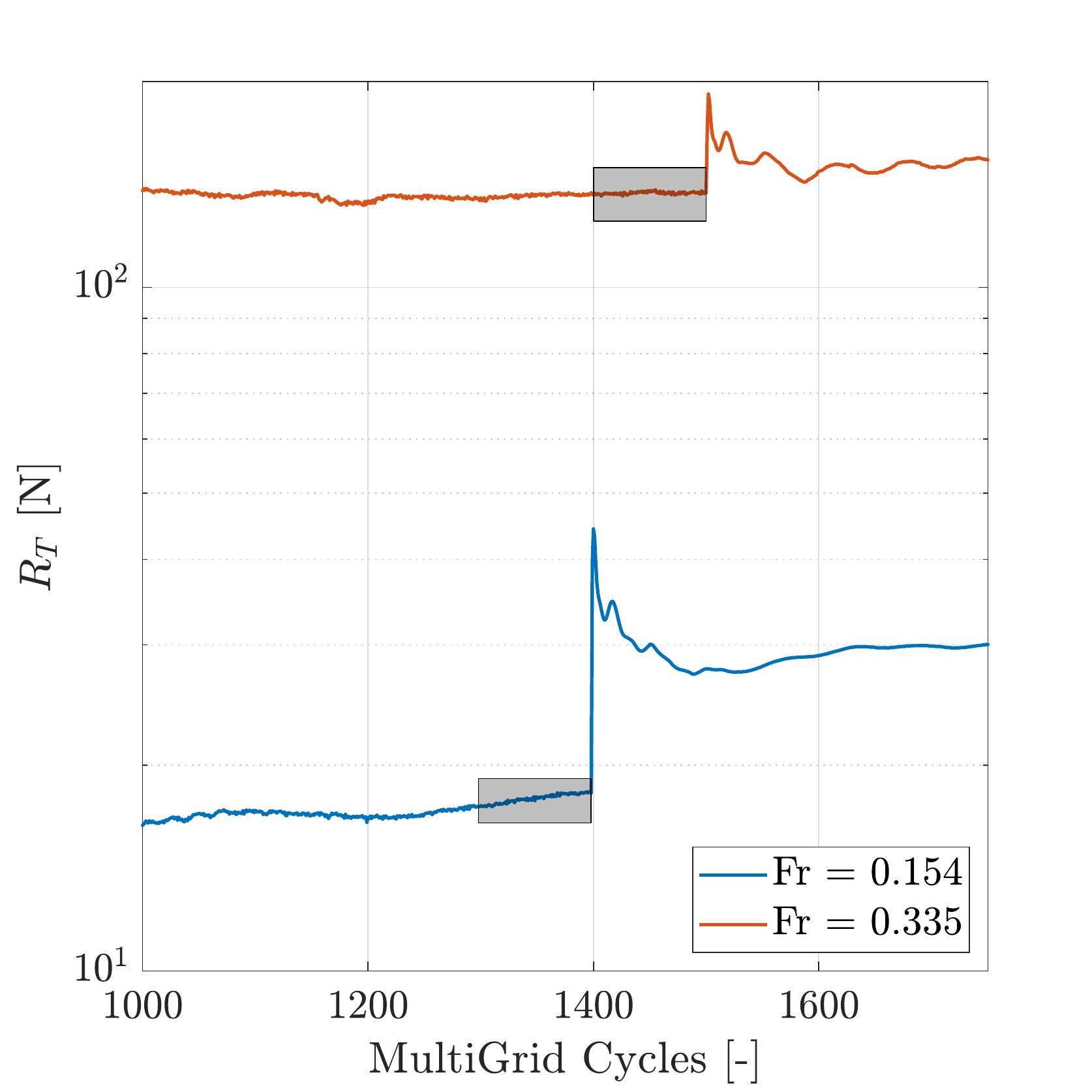}}
    \subfigure[Detail of the shift from $\mathcal{M}_2$ to $\mathcal{M}_3$]{\includegraphics[width=0.24\textwidth]{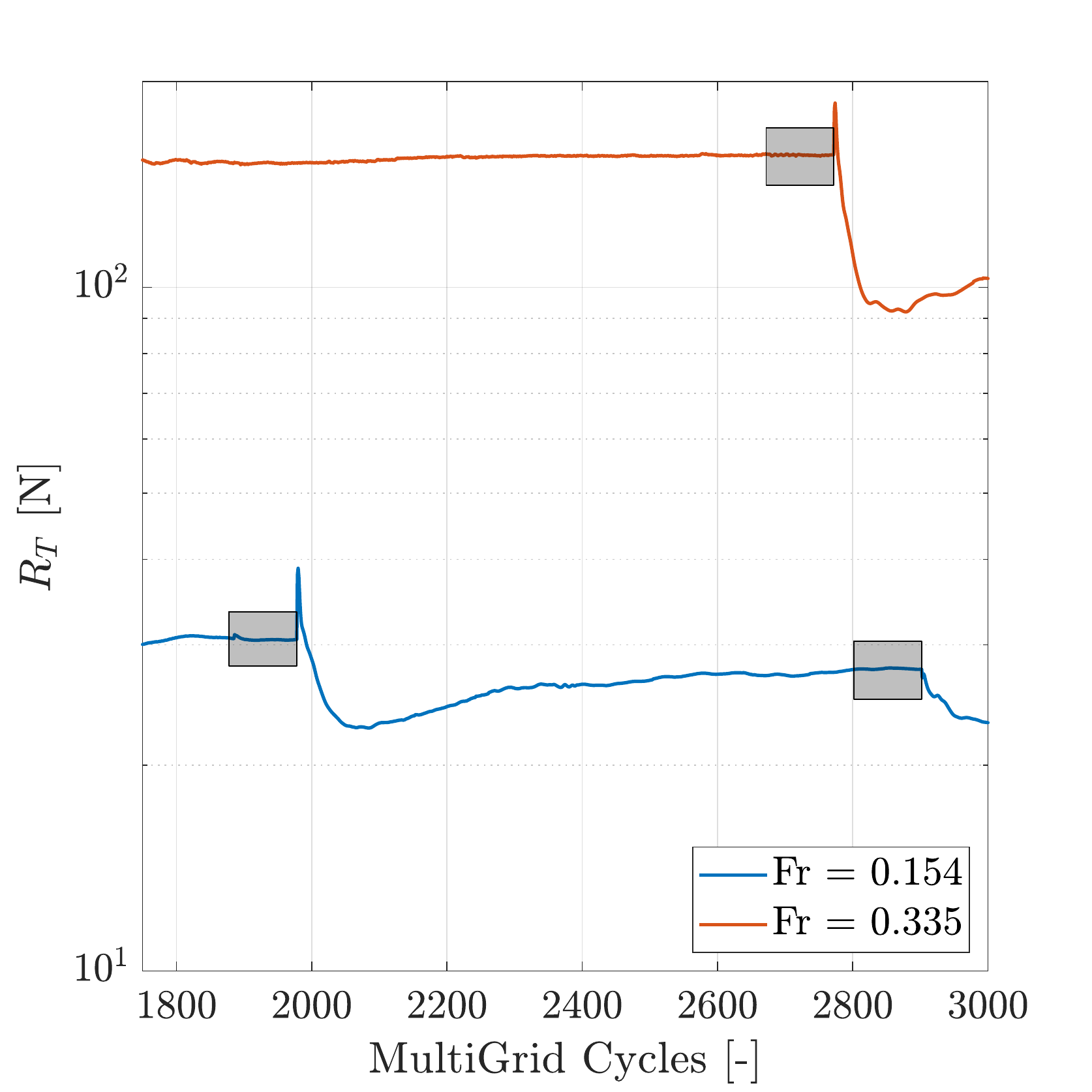}}
    \subfigure[Detail of the shift from $\mathcal{M}_3$ to $\mathcal{M}_4$]{\includegraphics[width=0.24\textwidth]{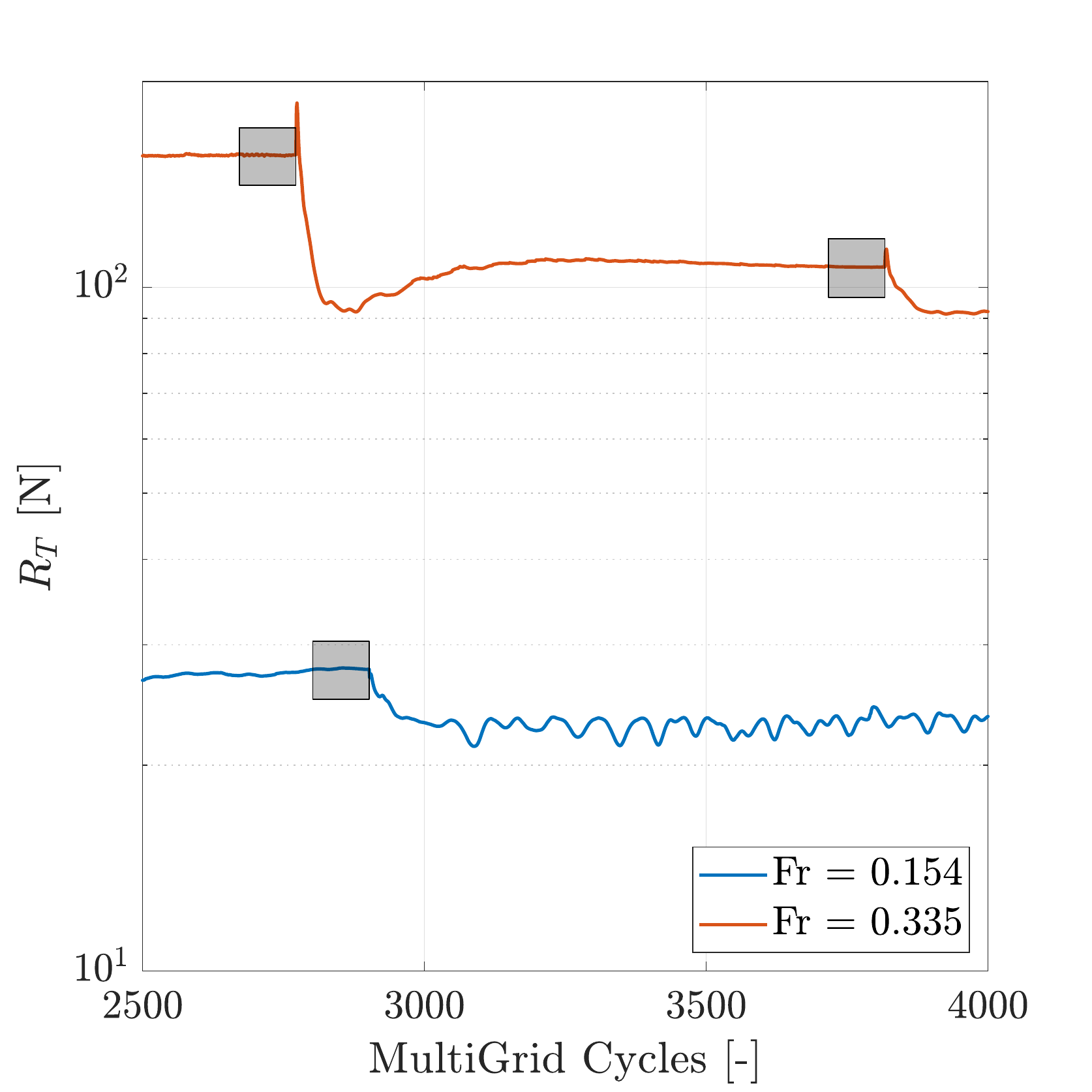}}
    \subfigure[Detail of the last iterations on $\mathcal{M}_4$]{\includegraphics[width=0.24\textwidth]{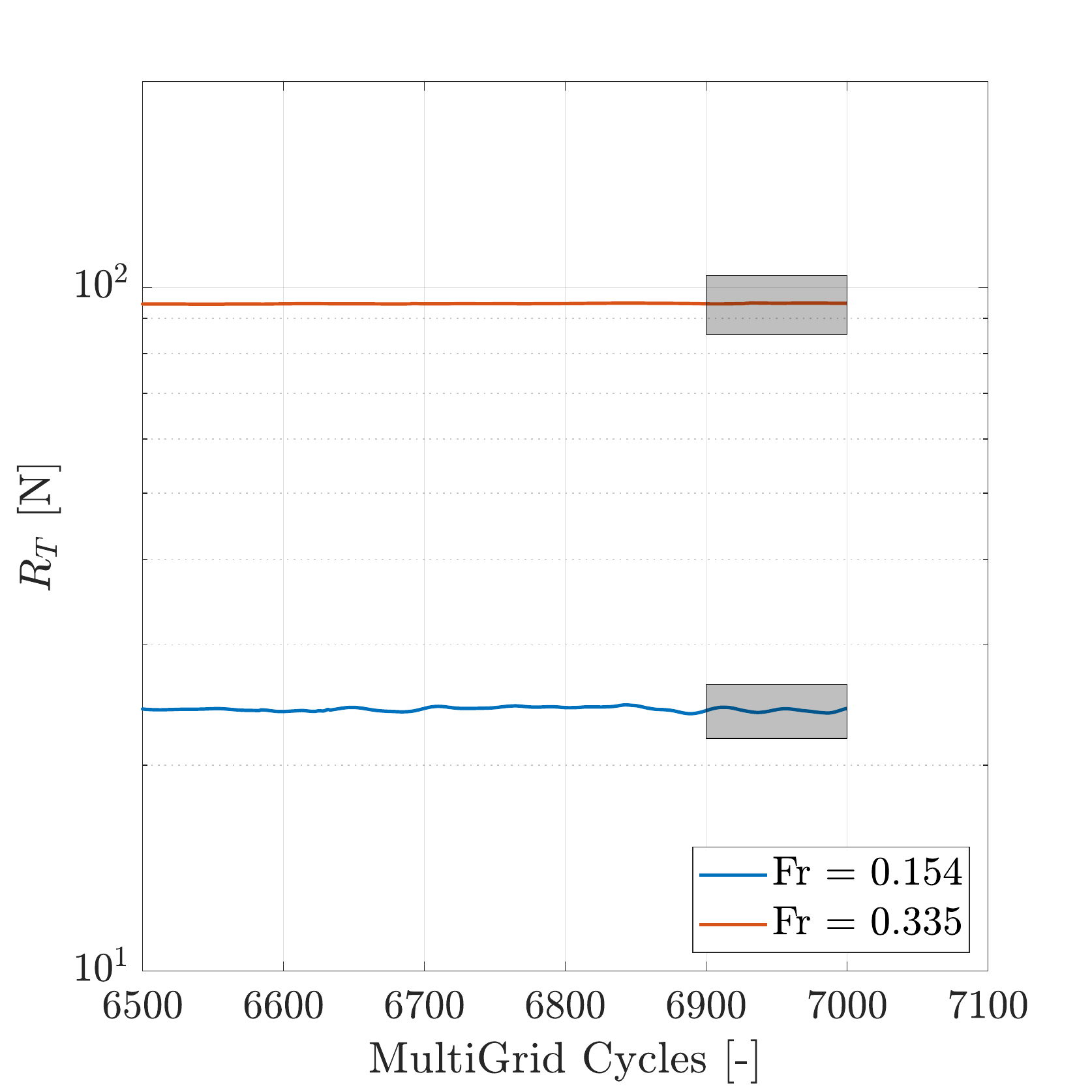}}
	\caption{RoPax problem, FMG--FAS cycle and highlight (grey areas) of the iterations considered to evaluate the $R_T$ for each grid for the minimum and maximum Froude numbers.}
	\label{fig:RoPax_noise}
\end{figure}

\begin{figure}[tp]
	\centering
	\subfigure[Low{est}-fidelity model]{\includegraphics[scale = 0.45]{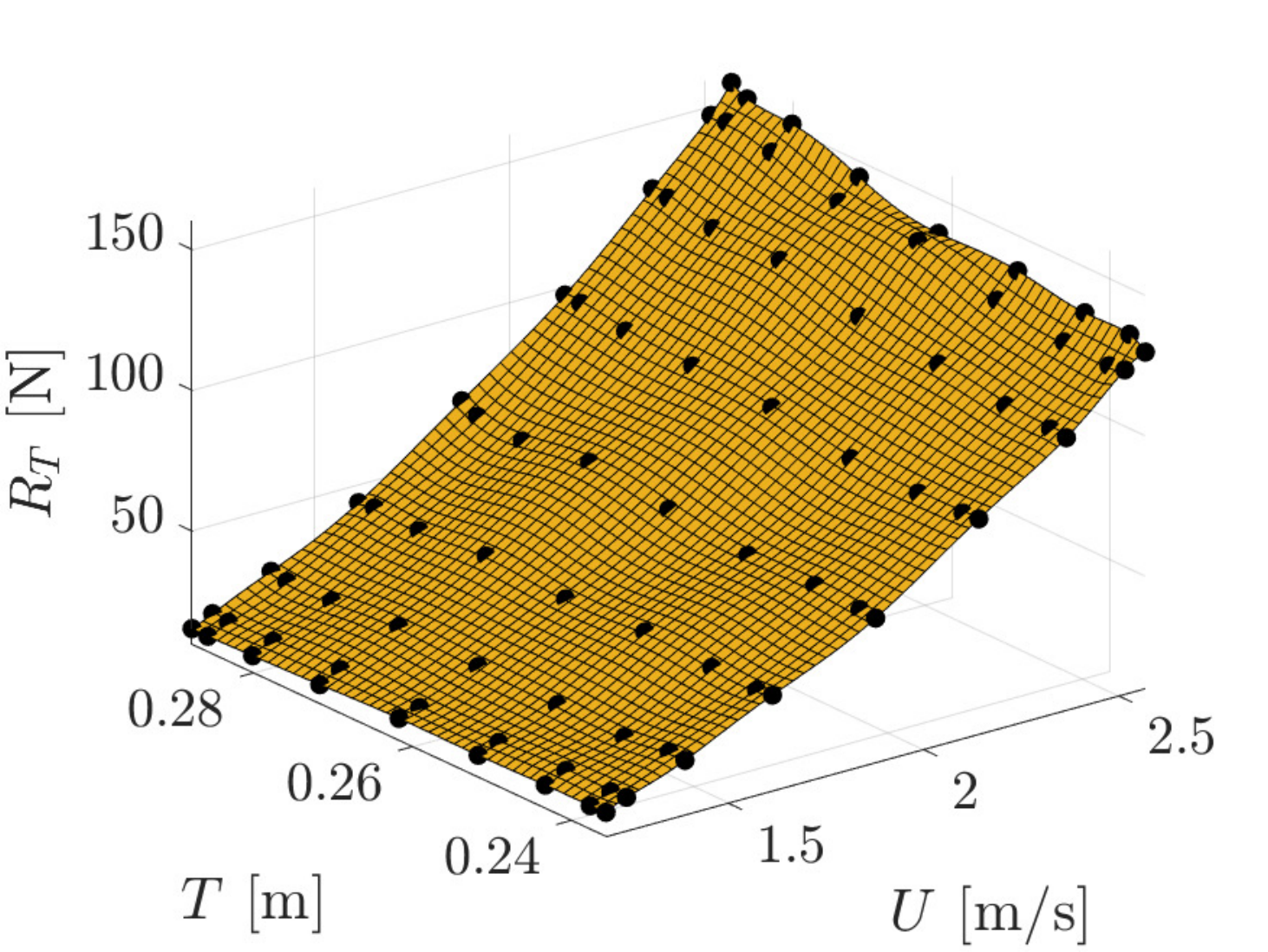}}
	\subfigure[High{est}-fidelity model]{\includegraphics[scale = 0.45]{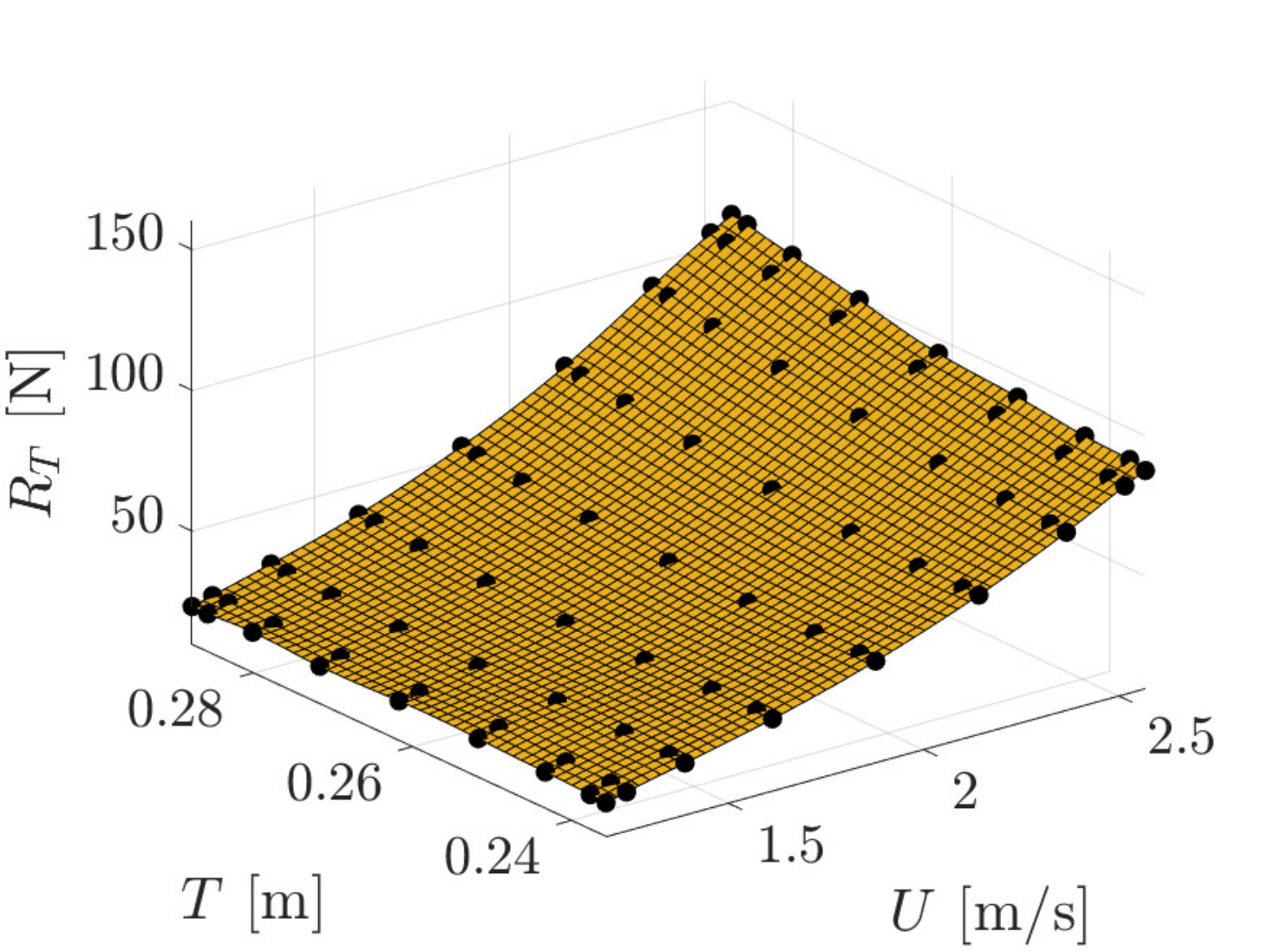}}\\
	\subfigure[MISC model]{\includegraphics[scale = 0.45]{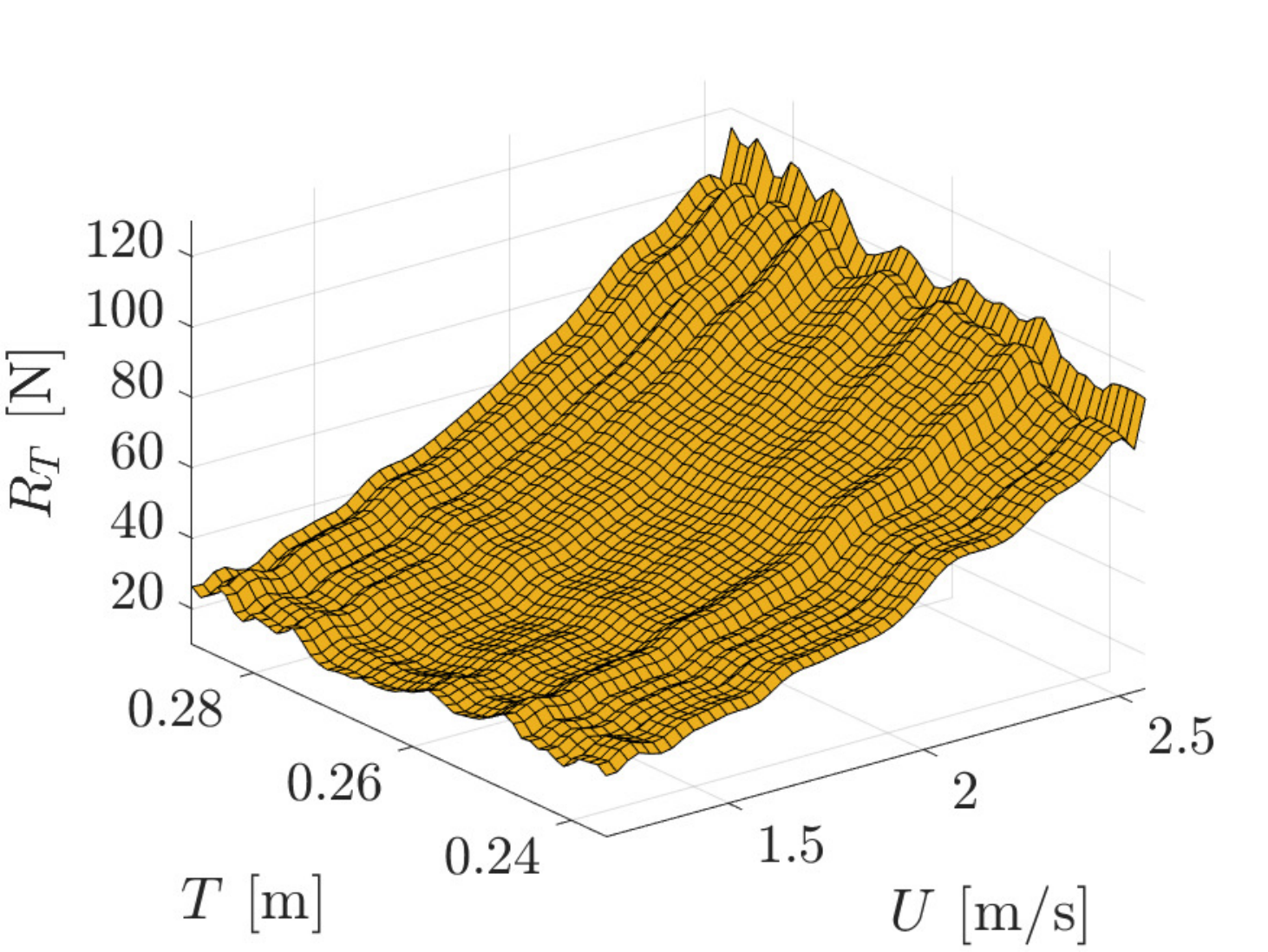}}
	\subfigure[SRBF model]{\includegraphics[scale = 0.45]{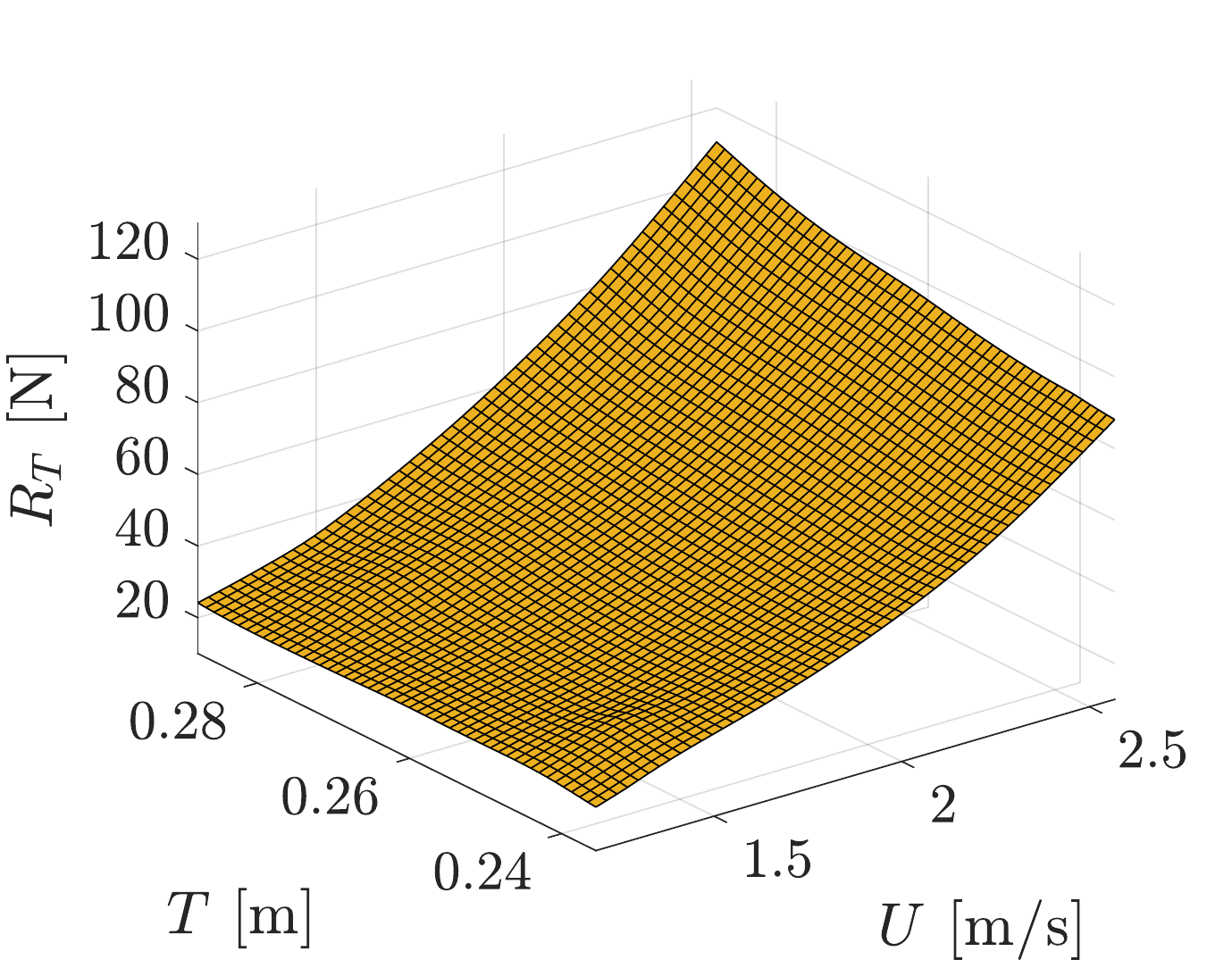}}
	\caption{RoPax problem, surrogate models. }
	\label{fig:RoPax_surrogate_models}
\end{figure}

\subsubsection{Numerical results}\label{sect:RoPax_results}

\begin{figure}[tp] 
	\centering
	\includegraphics[width=0.47\textwidth]{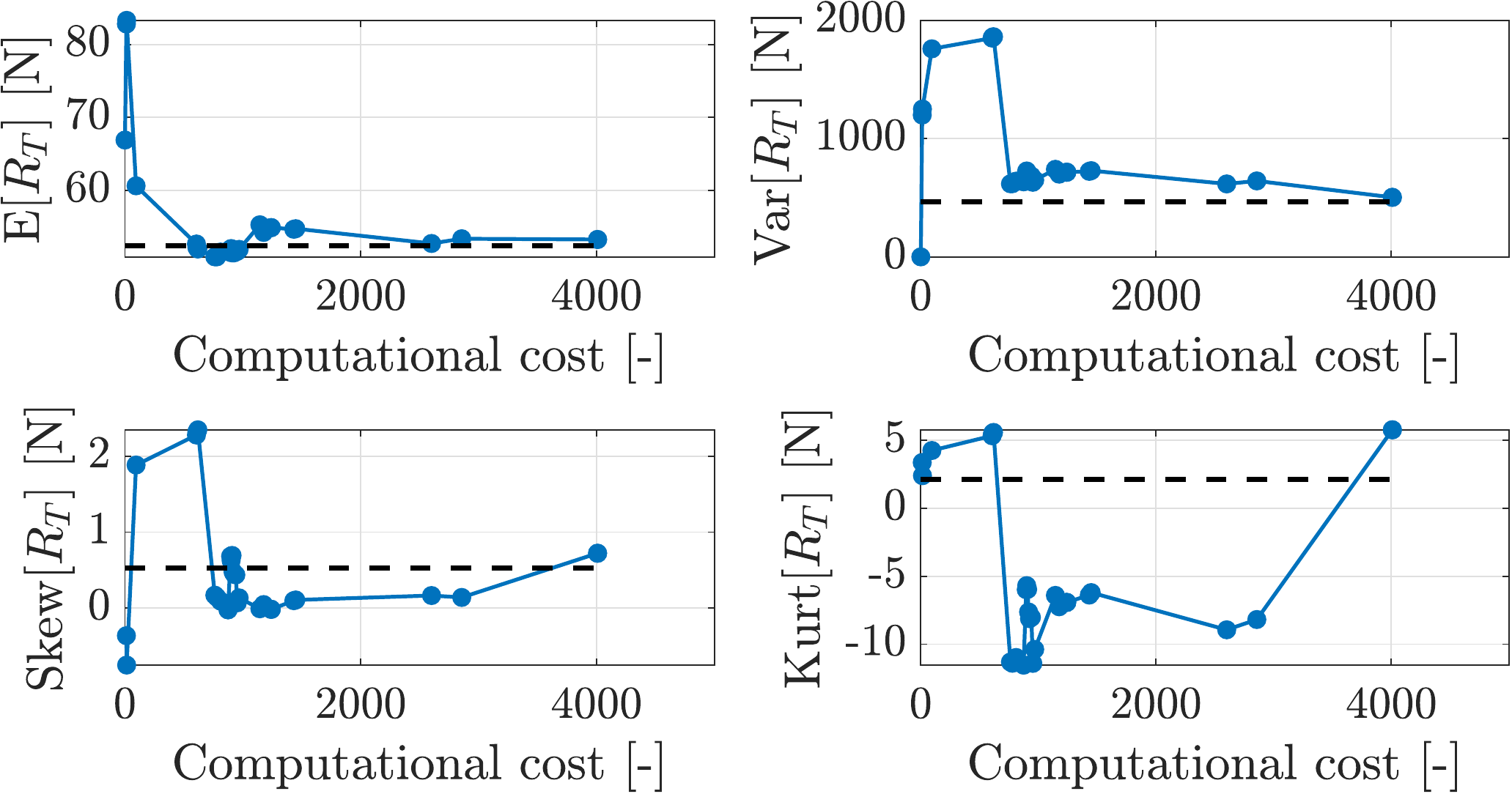}
	\includegraphics[width=0.47\textwidth]{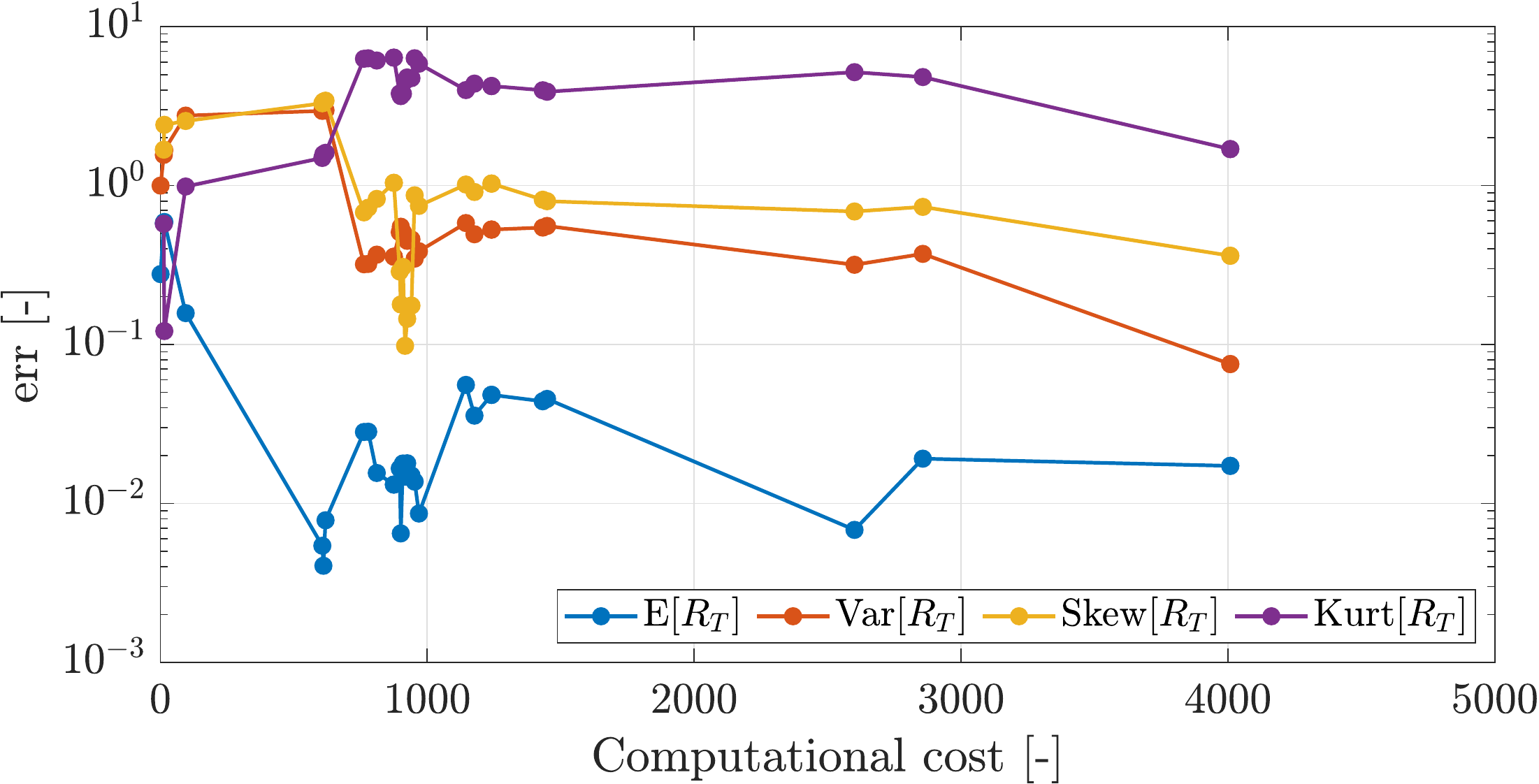}
	\caption{RoPax problem, results for the MISC method. Left: convergence of the values of the first four centered moments. The black dashed line marks the reference value of the moments. 
		Right: relative error of the moments (see Eq.~\eqref{eq:err_moments}).}
	\label{fig:RoPax_moments_MISC}
\end{figure}

Hereafter a detailed comparison of the performance of the MISC and SRBF methods is provided.
{The reference {surrogate model} ${G}_{\text{ref}}(\mathbf{y})$ is obtained {considering} highest-fidelity simulations only.} 
In details, an isotropic tensor grid consisting of {$9 \times 9$} CC points (see Eq.~\eqref{eq:CC_points}) is constructed {over $\Gamma$}
and the corresponding simulations on the grid $\mathcal{M}_4$ are performed. The resulting surrogate model is an interpolatory model,
based on global tensor Lagrange polynomials, which is shown in  Fig.~\ref{fig:RoPax_surrogate_models}b.
{Figure \ref{fig:RoPax_surrogate_models}a shows instead the surrogate obtained with simulations on $\mathcal{M}_1$,
  at the same CC points. Notice that both surrogates are affected by the noise,
  and more specifically, the noise is more evident in the lowest-fidelity surface
  which is significantly less smooth than the highest-fidelity surrogate.}
Reference values for the centered moments are then computed applying the tensor quadrature formula associated to the CC points to the highest-fidelity simulations.

First, the performance of the MISC method is discussed. Only the results
for the version of MISC with quadrature-based profits $P^{\MISCquad}$ are reported here,
{since} this approach outperforms the {version of MISC with surrogate-based profits $P^{\MISCsurr}$, for reasons related with the presence of the numerical noise that will be clarified in a moment. 
{In Fig.~\ref{fig:RoPax_moments_MISC} the values of the approximations of the first four centered moments of {$R_T$} obtained with MISC at different computational costs are displayed on the left, while the relative errors are shown on the right.
  Upon inspection of these plots, we can conclude that the quality of the estimates decreases with increasing order of the moments.
  In particular, the expected value and the variance seem to converge reasonably well (although the estimate of the expected value
  seems to hit a temporary plateau, after having obtained a good estimate at a low computational cost),
  whereas the kurtosis is strongly underestimated and its approximation results to be very poor.}

{To explain this behavior, we have a closer look at the MISC quadrature formula \eqref{eq:MISC_quad}.
  In particular, let us recall that the computation of the first four centered moments {implicitly uses}
  surrogate models for $r$th powers of the quantity of interest $R_T^r$, $r=1,\ldots,4$, (see Sect.~\ref{sect:misc}). 
  These models are displayed in Fig.~\ref{fig:RoPax_surfaces_MISC}.
  The first one, corresponding to $r=1$, is quite rough: the surface shows an oscillatory behavior and the expected monotonicity of $R_T$ with respect to $U$ and $T$ is destroyed. This is due to the already discussed presence of numerical noise, which particularly affects the low-fidelity simulations. 
Indeed, MISC intensively samples low-fidelities by construction, see Fig.~\ref{fig:RoPax_points_MISC}a,
where we report the evaluations allocated on each fidelity as the iterations proceed. 
In particular, most of the low-fidelity simulations are added from iteration 13 on (see Fig.~\ref{fig:RoPax_points_MISC}b):
this explains that the estimate of $\mathbb{E}[R_T]$ reaches reasonable values at early iterations,
i.e. when there is still a balance of low- and higher- fidelity simulations, and its convergence deteriorates later, i.e. when low-fidelity simulations are the majority. 
Given that the numerical noise introduces spurious oscillations in the surrogate model already for $r=1$,
such oscillations can only be amplified for $r>1$, as can be observed in Fig.~\ref{fig:RoPax_surfaces_MISC}b,c,d.
Hence, the computation of moments suffers more from the noise the higher the order.}

This observation then suggests {that} a way to mitigate the impact of such oscillations in the computation of 
{statistical moments is to employ a} method that does not require higher order surrogate models.
In this work, {we propose to compute such moments by taking Monte Carlo samples of the surrogate model of $R_T$,
  and approximate the moments from these values, with the usual sample formulas.}
The results reported in Fig.~\ref{fig:RoPax_MISC_MC}
have been obtained taking the average of 10 repetitions with 10000 samples
(again, note that this computation is not expensive since it only requires evaluations of the MISC surrogate model)
and are quite promising: the benefits increase for higher and higher order moments,
and in particular, the improvement in the estimate of the kurtosis is quite impressive;
the results of MISC and Monte Carlo quadrature in the case $r=1$ are {instead} substantially equivalent,
{which is to be expected since} they both work with the same surrogate model.
{This strategy thus mitigates the impact of the noise on the computation of moments.
  However, it is not entirely satisfactory, since the choice of the number of samples to be employed is non trivial:
  on the one hand, we need a sufficiently large number of samples to ensure accuracy of the estimates, on the other hand,
  taking too many samples results in resolving the spurious oscillations.
  In other words, the chosen number of samples should give the best compromise between these two aspects, and some trial-and-error
  tuning, or some accurate a-priori analysis should be carried out. Deeper studies of this matter will be subject of future work. }   

{Further, at this point it is also clear why the version of MISC with profit $P^{\MISCquad}$, i.e. based on the quadrature, gives better results than the version with profit $P^{\MISCsurr}$, i.e. based on the quality of the surrogate model. Indeed, in case of noisy simulations, {the} adaptivity criterion based on the profit $P^{\MISCsurr}$ 
results in capturing even more the spurious oscillations due to the numerical noise, since it is based on the direct point-wise evaluation of the noisy surrogate model.}
\begin{figure}[tp] 
	\centering
	\subfigure[Surrogate model for $R_T$]{\includegraphics[scale = 0.45]{Fig_test_ferry/RoPax_surface_MISC}}
	\subfigure[Surrogate model for $R_T^2$]{\includegraphics[scale = 0.45]{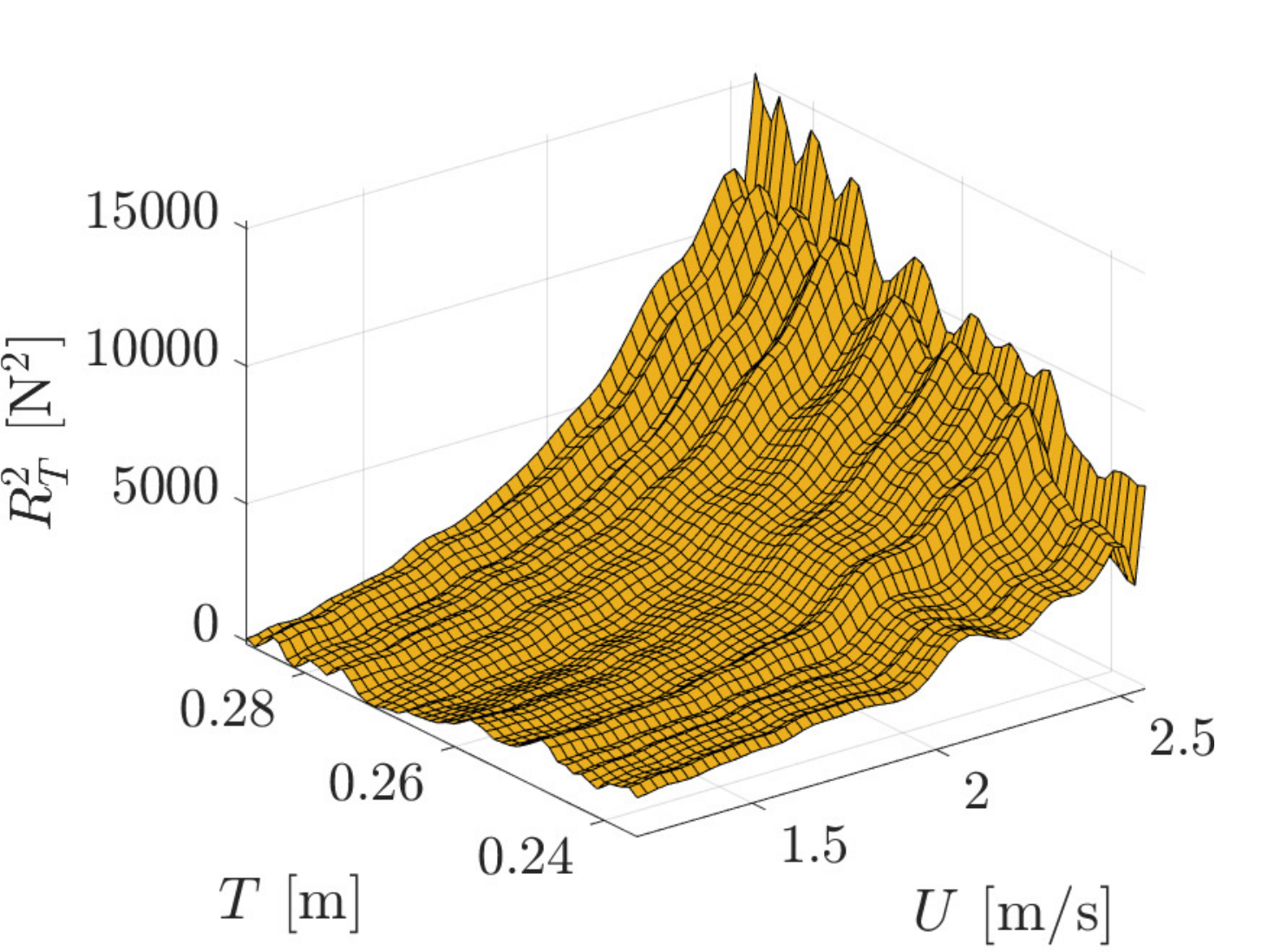}}\\
	\subfigure[Surrogate model for $R_T^3$]{\includegraphics[scale = 0.45]{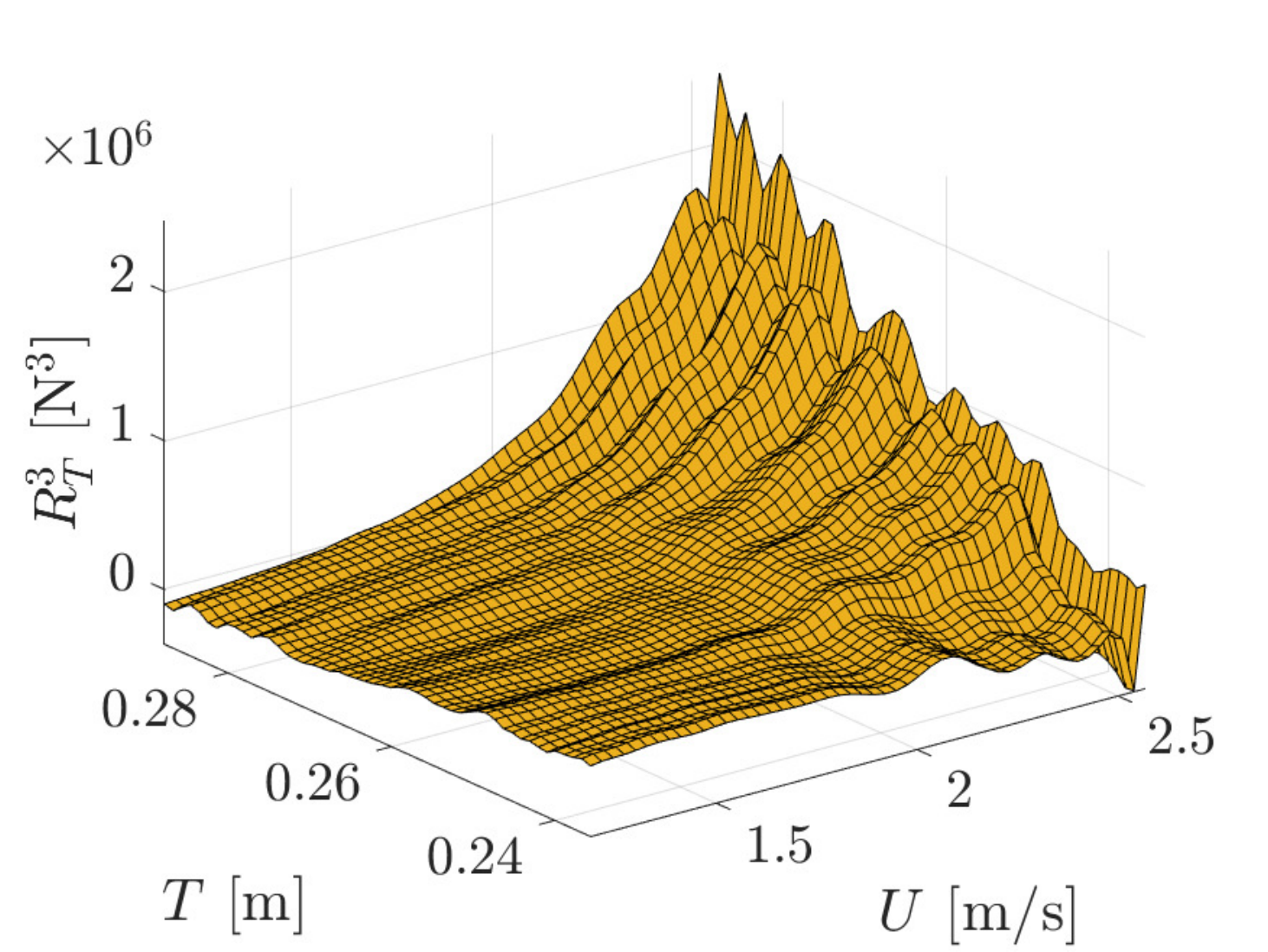}}
	\subfigure[Surrogate model for $R_T^4$]{\includegraphics[scale = 0.45]{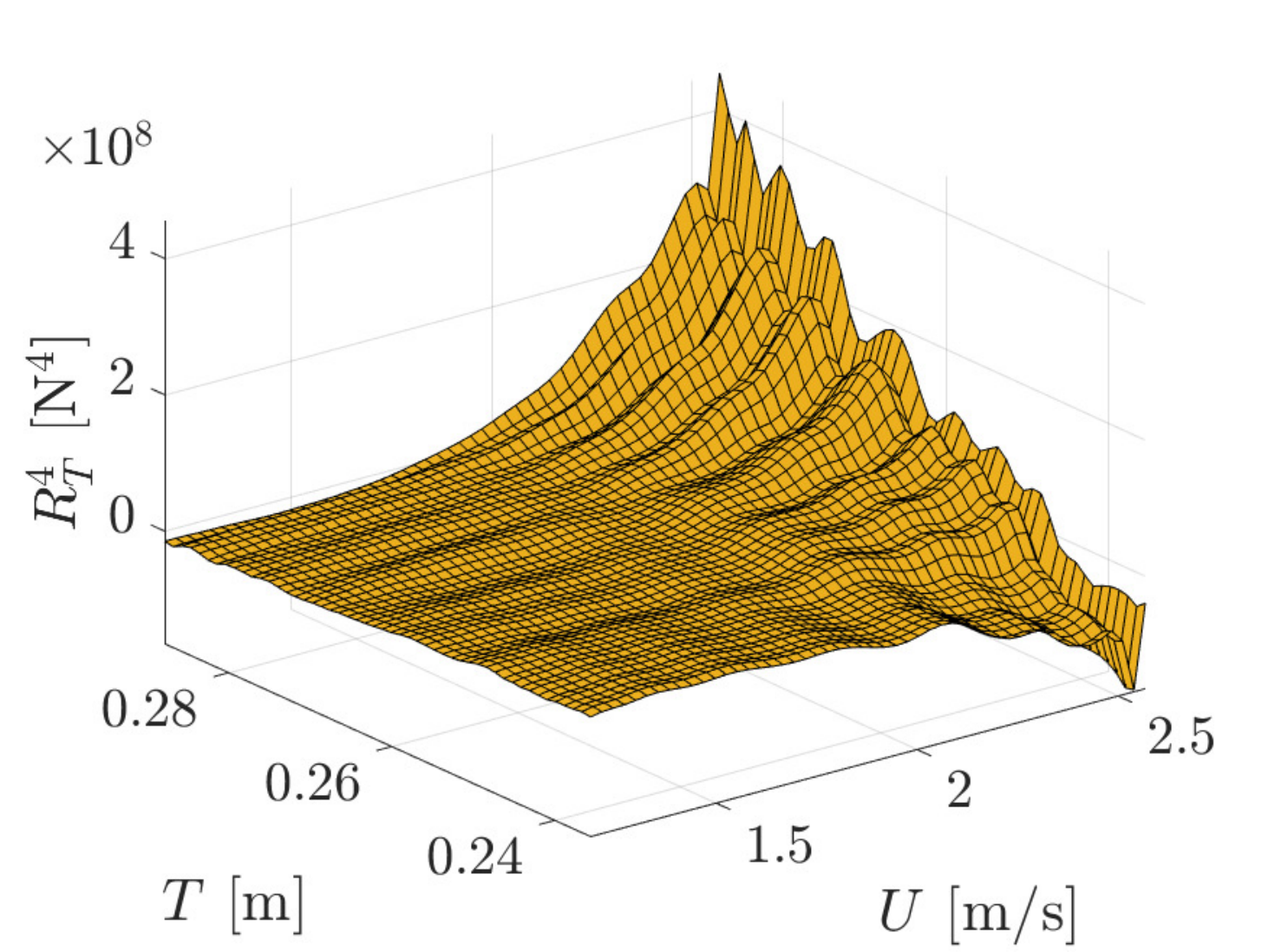}}
	\caption{RoPax problem, surrogate models obtained with the MISC method.}
	\label{fig:RoPax_surfaces_MISC}
\end{figure}

\begin{figure}[tp] 
	\centering
	\subfigure[Total number of simulations]{\includegraphics[scale = 0.38]{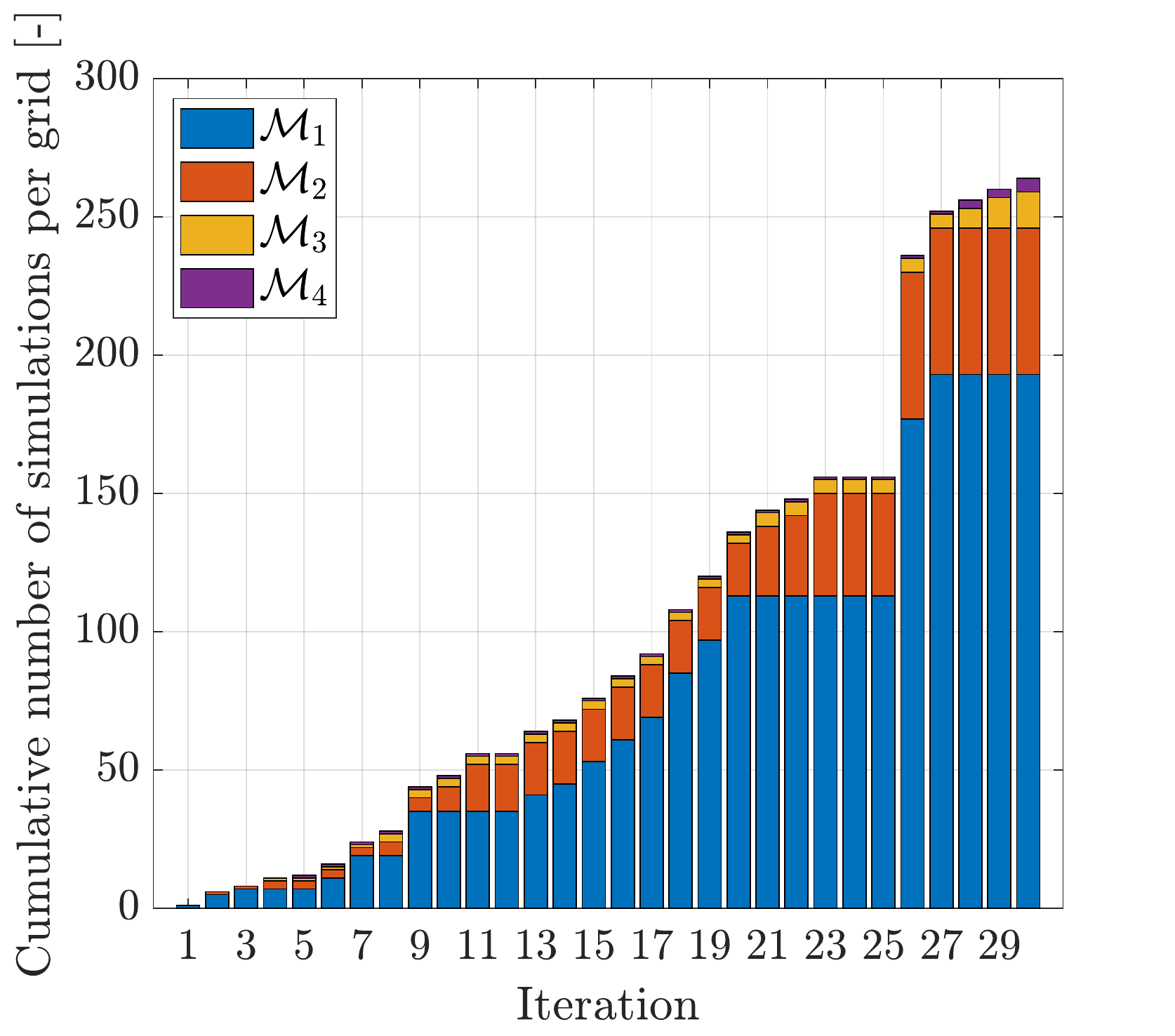}}
	\subfigure[Number of new simulations]{\includegraphics[scale = 0.38]{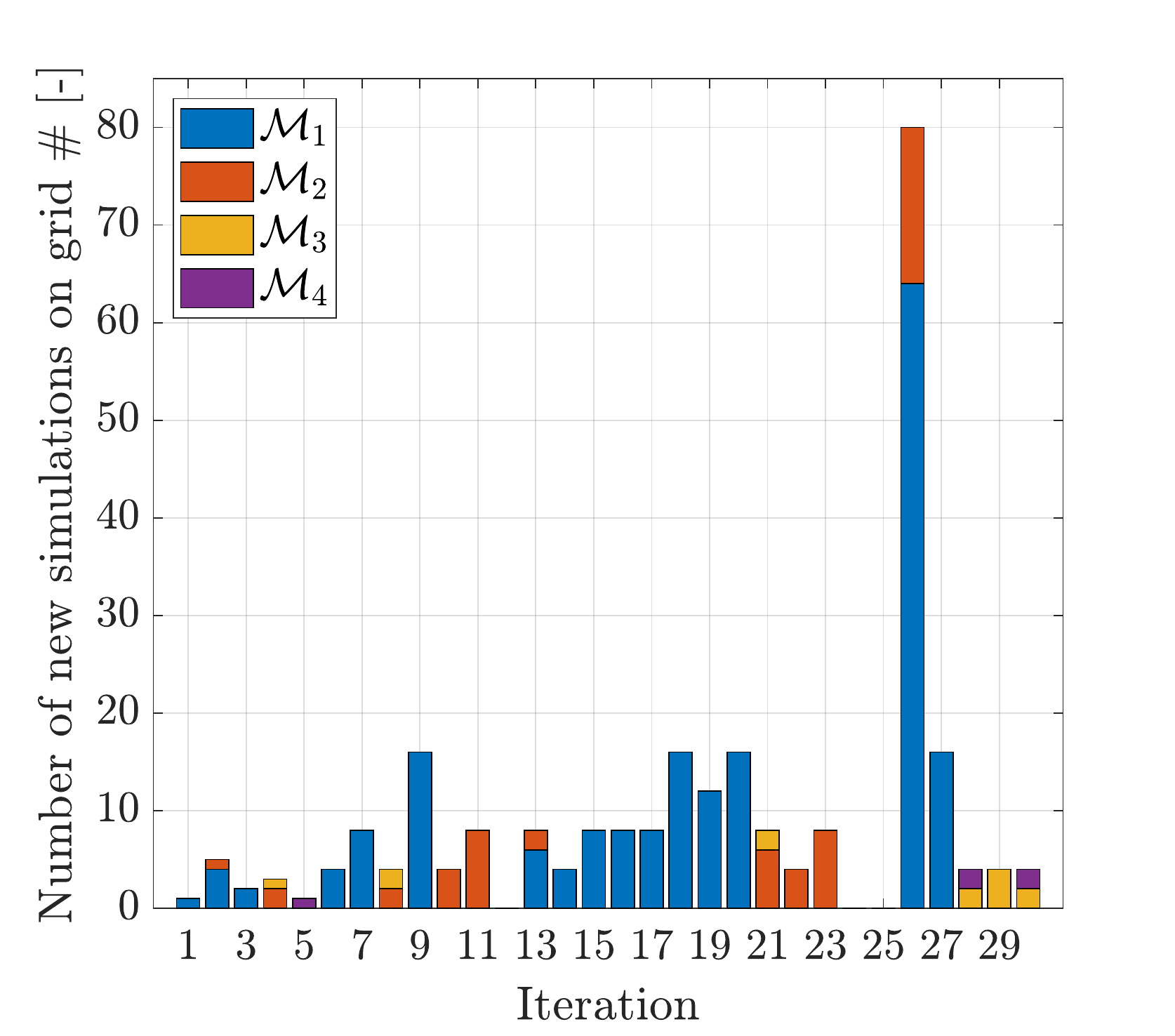}}
	\caption{RoPax problem, number of simulations per grid selected by the MISC method at each iteration.}
	\label{fig:RoPax_points_MISC}
\end{figure}

\begin{figure}[tp] 
	\centering
	\includegraphics[width = 0.94\textwidth]{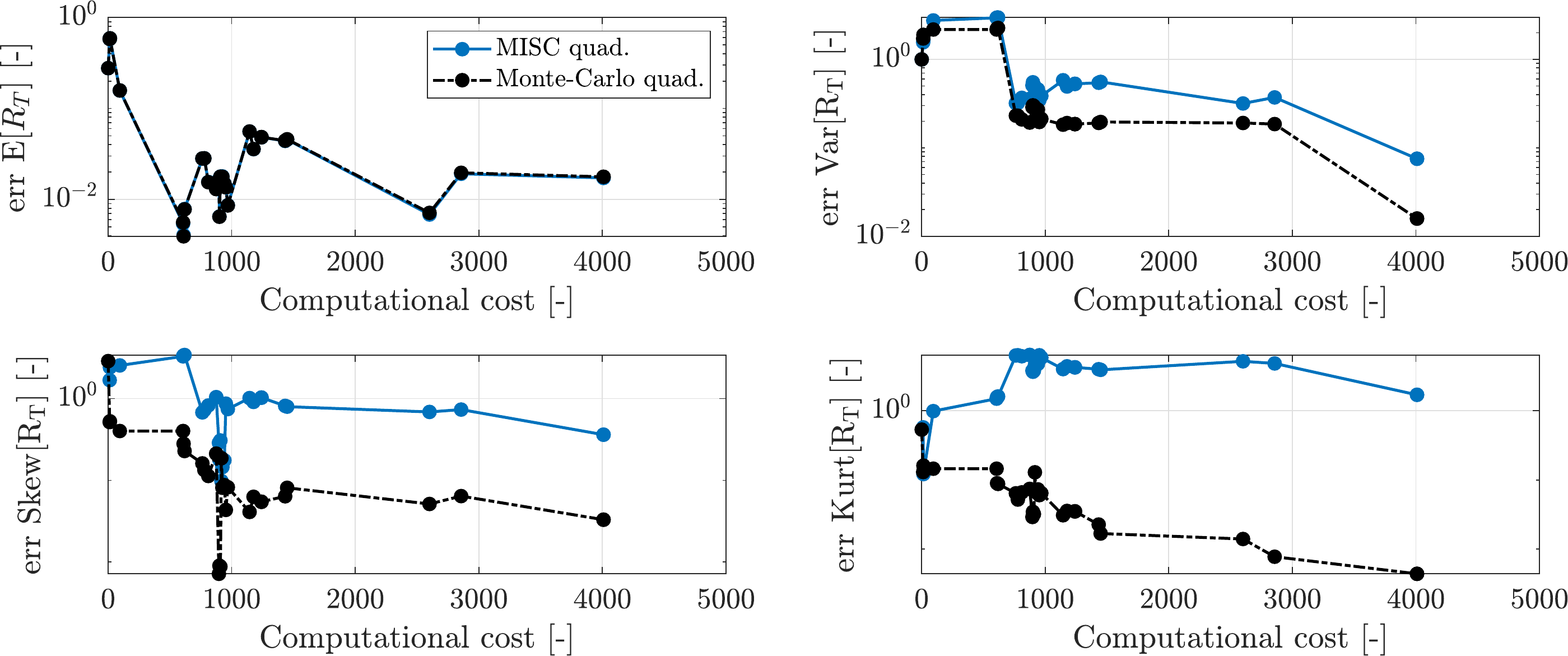}
	\caption{RoPax problem, results for the MISC method: relative error of the moments (see Eq.~\eqref{eq:err_moments}). A comparison between the results obtained with the MISC quadrature formula \eqref{eq:MISC_quad} and the Monte Carlo quadrature.}
	\label{fig:RoPax_MISC_MC}
\end{figure}

\begin{figure}[tp] 
    \centering
    \subfigure[Total number of simulations]{\includegraphics[scale = 0.38]{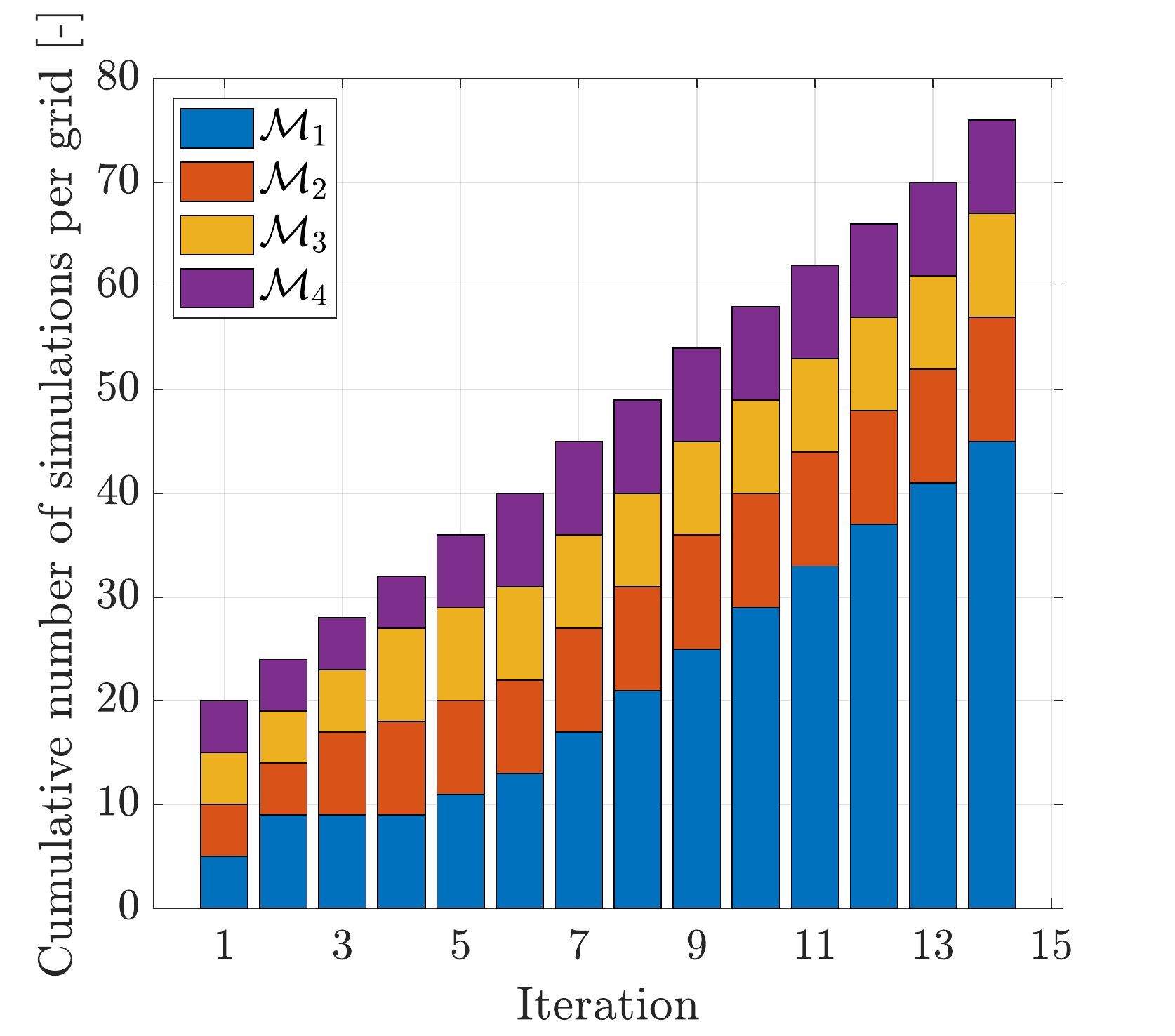}}
    \subfigure[Number of new simulations]{\includegraphics[scale = 0.38]{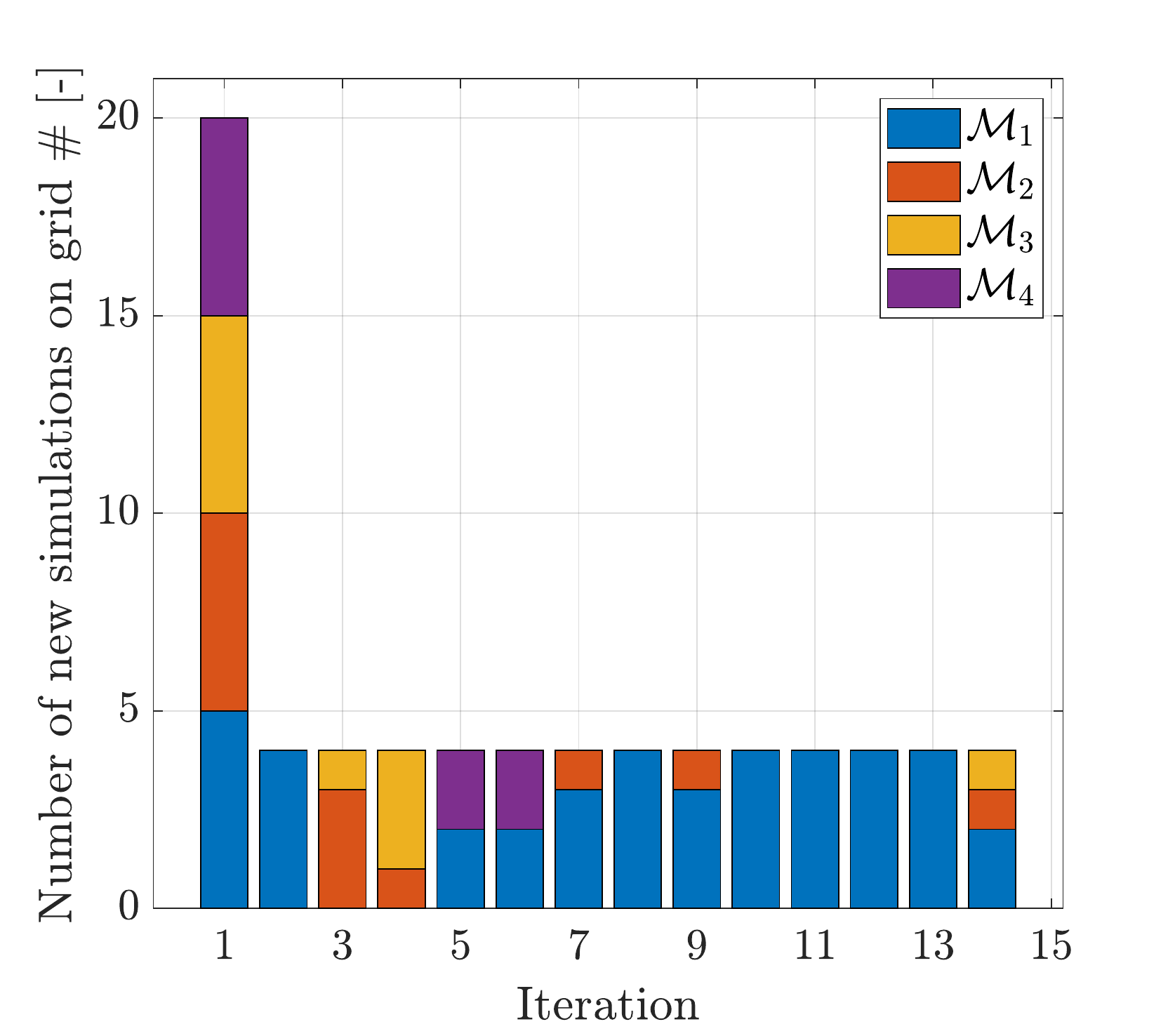}}
    \caption{RoPax problem, number of simulations per grid selected by the SRBF method at each iteration.}
    \label{fig:RoPax_points_SRBF}
\end{figure}

{Next, we move to SRBF. As already mentioned, in this application we use
  a regression approach to compute the weights of the surrogate model
  (i.e. solving Eq. \eqref{eq:lsrbf}), motivated by the fact that the evaluations of
  the CFD solver are noisy as just discussed.}
The SRBF surrogate model at the last iteration of the adaptive sampling procedure is shown in Fig.~\ref{fig:RoPax_surrogate_models}d. The surface is smoother than the one produced by MISC (cf. Fig.~\ref{fig:RoPax_surrogate_models}c), although a small bump is present in the bottom part. This figure thus shows that SRBF is in general effective in filtering-out the numerical noise in the training set. 

Figure \ref{fig:RoPax_points_SRBF} shows that SRBF spent about $50\%$ of the final computational cost at the first iteration, then requiring simulations on the finest grids only at iterations 5 and 6. In all the other iterations mainly low-fidelity simulations are performed. This sampling behavior is due to the high values of prediction uncertainty that are found in the corners of the {parametric} domain, because the topology of the initial training leads to extrapolation in those zones. Such corner regions are those with the highest estimated prediction uncertainty, and the adaptive sampling procedure requires all the fidelities before exploring other regions.

Figure \ref{fig:RoPax_moments_SRBF} shows the convergence of the first four centered moments of $R_T$ and their relative errors obtained with SRBF. The method initially converges rapidly to  the reference values but then the metrics start oscillating. {Similarly to the analytical problem, the errors of all the moments have a quite similar convergence.} This is particularly evident in the last iterations of the adaptive sampling process, see Fig.~\ref{fig:RoPax_moments_SRBF}a. 
Figure~\ref{fig:RoPax_moments_SRBF}b shows the detail of the last iterations of the adaptive sampling. The oscillatory behavior is evident and mostly associated with the {intensive sampling of the lowest fidelity happening in correspondence with computational costs between 4638 and 4655}. {In this range the expected value and the skewness oscillates more than the other moments, indicating that the surrogate model is oscillating around the training data.} 

\begin{figure}[tp] 
	\centering
	\subfigure[Convergence over the full range of available computational costs.]{
		\includegraphics[width=0.47\textwidth]{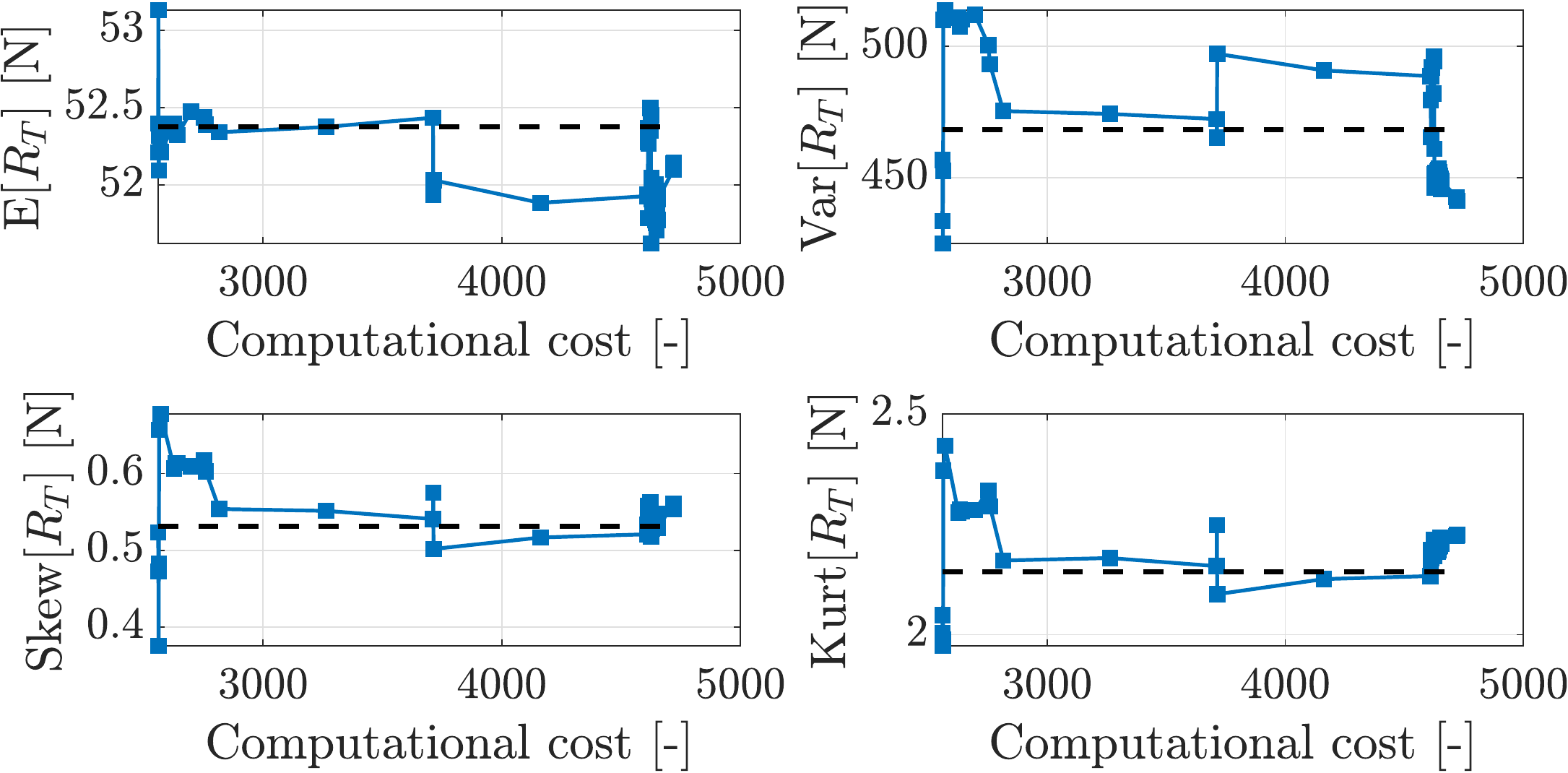}
		\includegraphics[width=0.47\textwidth]{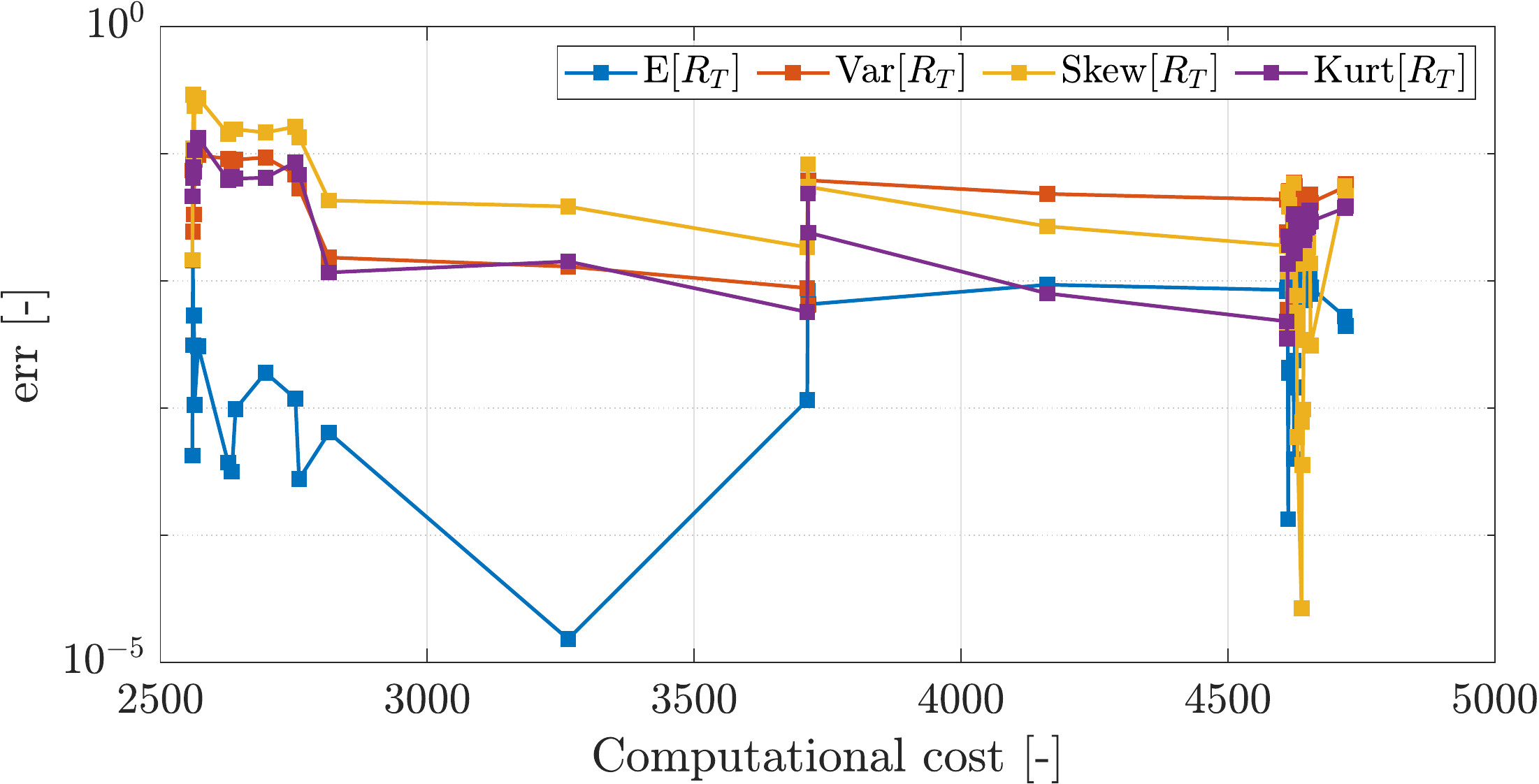}}
	\subfigure[Zoom on the final part of the convergence.]{
		\includegraphics[width=0.47\textwidth]{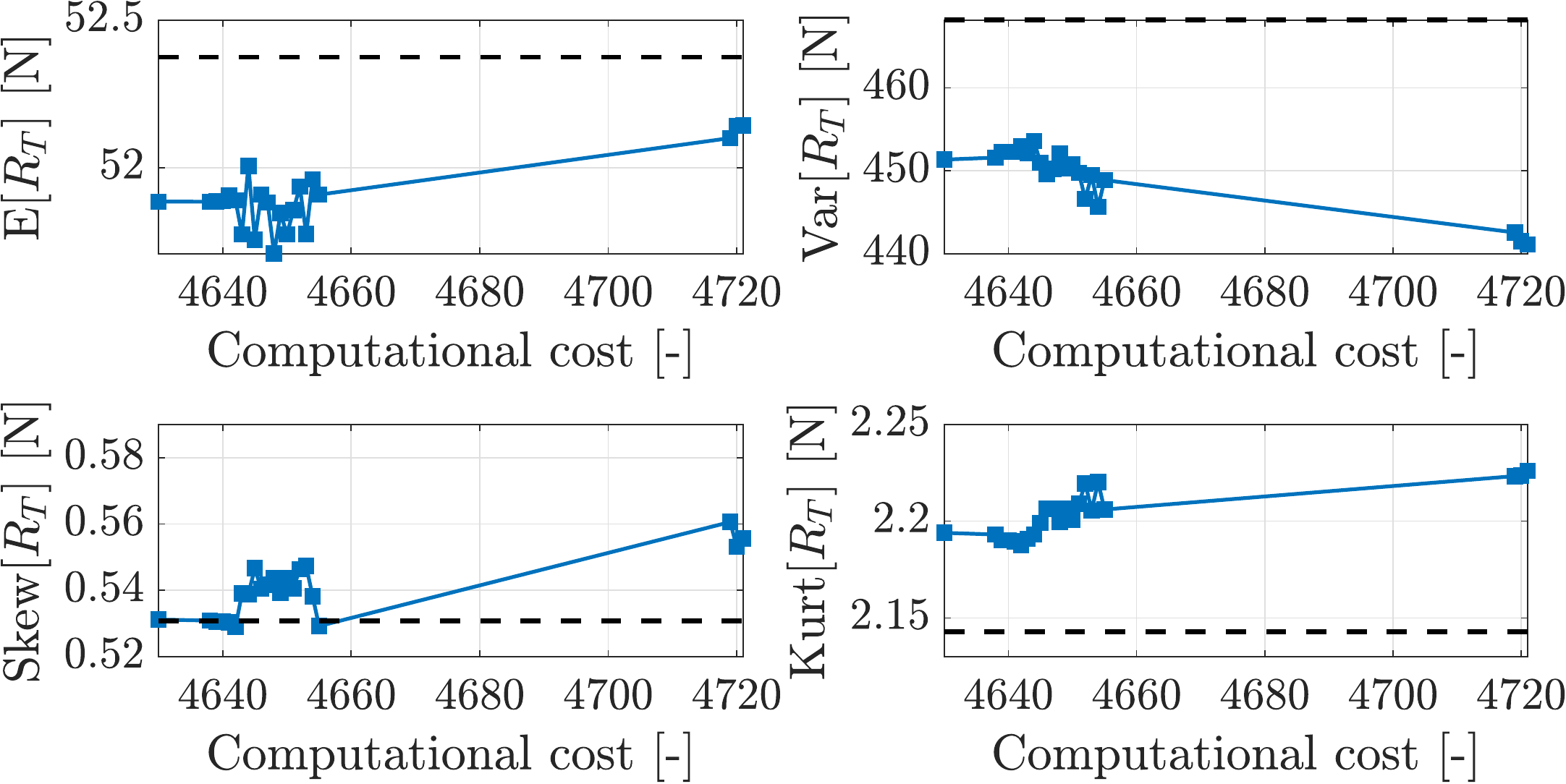}
		\includegraphics[width=0.47\textwidth]{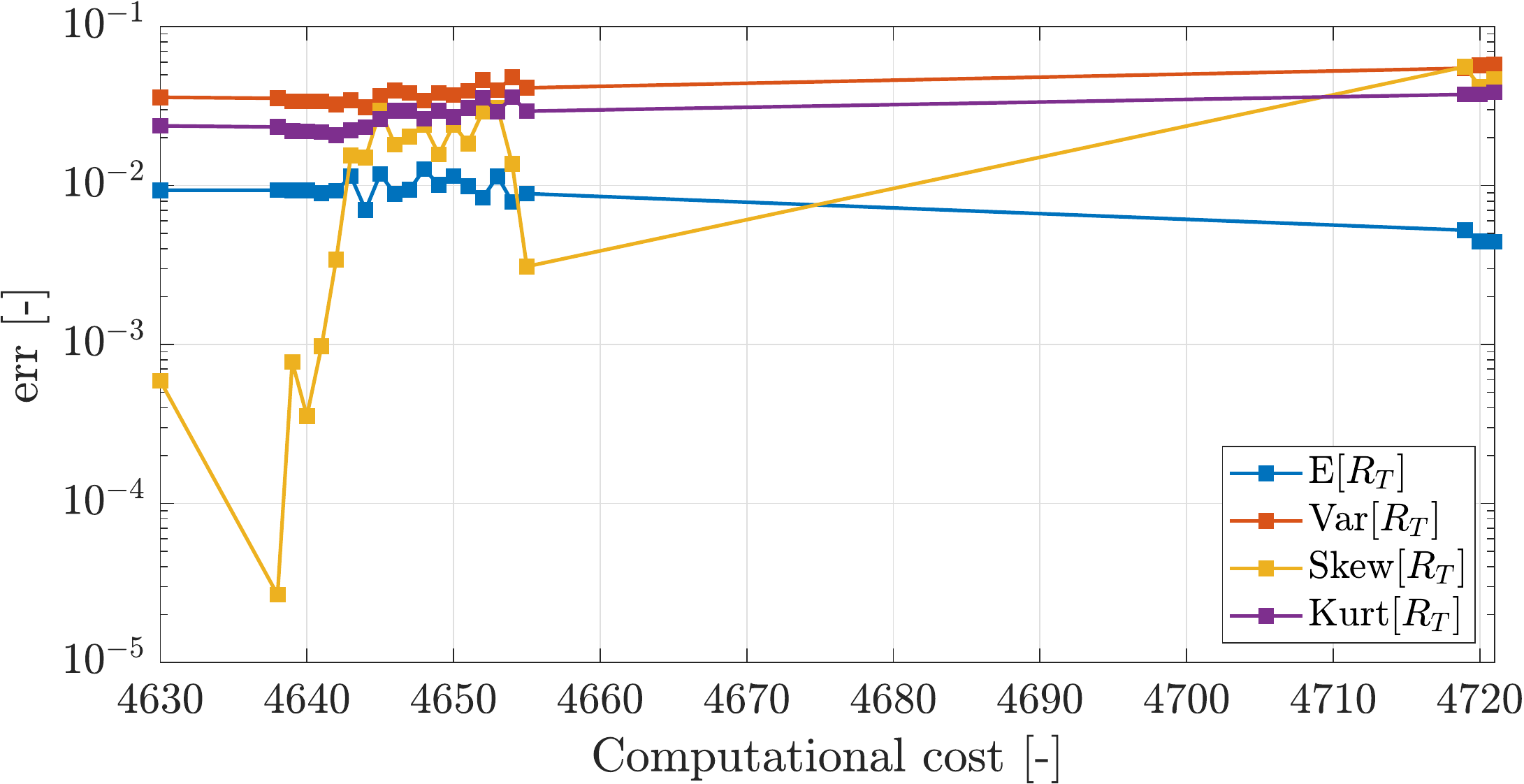}}
	\caption{RoPax problem, results for the SRBF method. Left: values of the first four centered moments. The black dashed line marks the reference value of the moments. 
		Right: relative error of the moments (see Eq.~\eqref{eq:err_moments}).}
	\label{fig:RoPax_moments_SRBF}
\end{figure}

{To conclude the discussion on the convergence of the moments,} in Fig.~\ref{fig:RoPax_err_moments} the convergence of the relative errors of the moments obtained with the two methods is compared (of course we consider MISC results when using Monte Carlo quadrature). Both MISC and SRBF achieve similar values: SRBF reach smaller errors in all moments but the kurtosis, but the convergence trend is bumpier than for MISC, and has a larger initial computational cost. 
{Especially for SRBF, the convergence of the moments is oscillating, nevertheless the oscillations fall within a range that can be considered reasonable from a practical viewpoint as it is comparable with  the numerical uncertainties and/or noise of the solver.}

\begin{figure}[tp] 
	\centering
	\includegraphics[width = 0.94\textwidth]{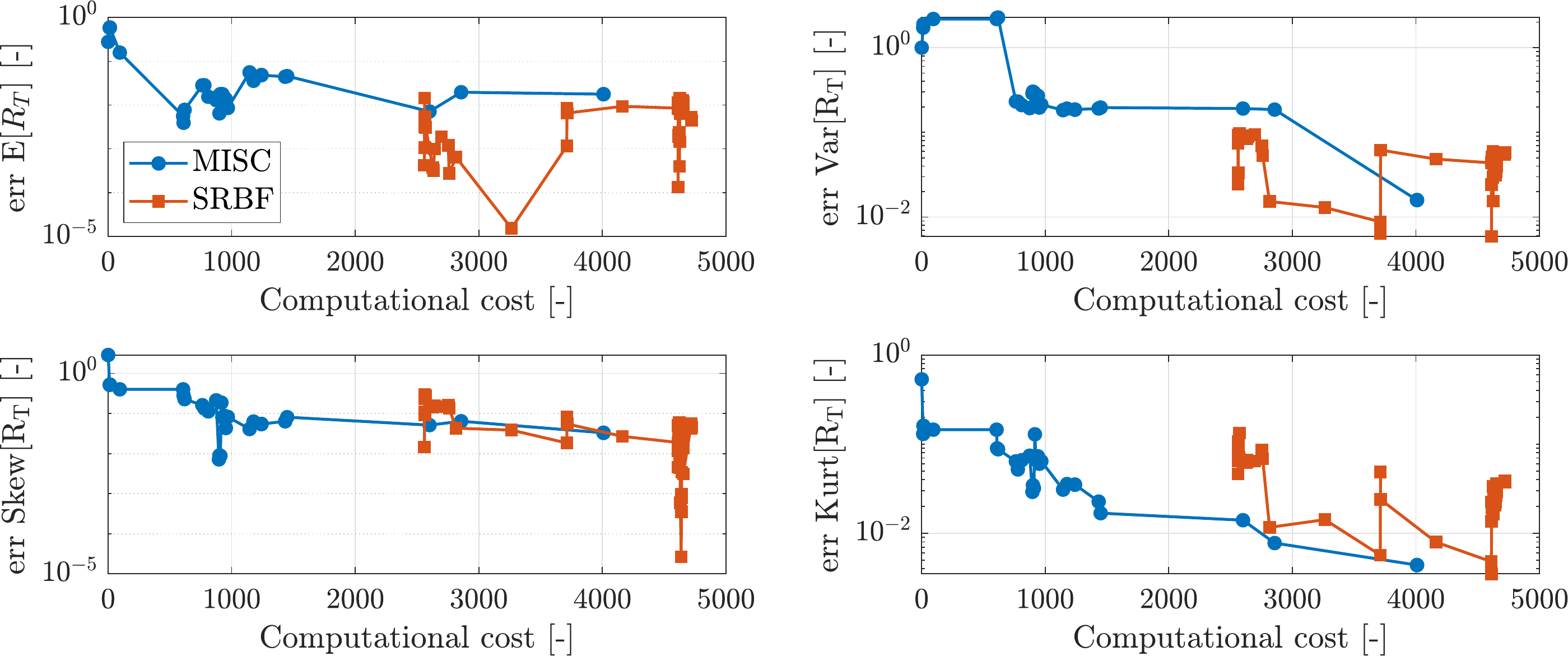}
	\caption{RoPax problem, comparison of MISC and SRBF results: relative error of the moments. It is a compact visualization of the results of Fig.~\ref{fig:RoPax_MISC_MC} and \ref{fig:RoPax_moments_SRBF}.}
	\label{fig:RoPax_err_moments}
\end{figure}

In Fig.~\ref{fig:RoPax_norms} the relative error in $L_2$ and $L_{\infty}$ norms of the estimates of the advancement resistance are plotted.
These metrics confirm that the MISC method reaches reasonable estimates with a low computational cost, whereas the SRBF method returns slightly better results but requires an higher computational cost. 
In the case of SRBF, it is worth noting that in the last iterations the relative errors of variance, skewness, and kurtosis increase whereas the $L_2$ metric decreases.
{This apparent contradiction is discussed comparing the convergence of the variance and the $L_2$ metric.
Figure~\ref{fig:RoPax_srbf_varl2}a and d show the convergence of the difference of the variances between the multi-fidelity surrogate model prediction $R_T$ and the reference value $R_T^*$, i.e. of $\Delta \text{Var}[R_T]=\text{Var}[{R}_T]-\text{Var}[R_T^*]$, and of the $L_2^2(R_T)$ metric, along with the summands of their decompositions: 
$\Delta \text{Var}[R_T]=\mathbb{E}[{R}_T^2]-\mathbb{E}[R_T^{*2}]-(\mathbb{E}[{R}_T]^2-\mathbb{E}[R_T^*]^{2}), \quad  L_2^2(R_T)=\mathbb{E}[\left({R}_T-{R}_T^*\right)^2]=\mathbb{E}[{R}_T^2]+\mathbb{E}[R_T^{*2}]-2\mathbb{E}[{R}_T R_T^*].$
In order to have $\Delta \text{Var}[R_T]$ going to zero it should happen that its two components, $\mathbb{E}[{R}_T^2]-\mathbb{E}[R_T^{*2}]$ and $\mathbb{E}[{R}_T]^2-\mathbb{E}[R_T^*]^{2}$, go to zero remaining equal in size; however, this does not happen, see Fig.~\ref{fig:RoPax_srbf_varl2}b and e. 
Conversely, the two components of the $L_2^2(R_T)$ metric, i.e. $\mathbb{E}[{R}_T^2]+\mathbb{E}[R_T^{*2}]$ and $2\mathbb{E}[{R}_T R_T^*]$, are
always almost equal in size (see Fig.\ \ref{fig:RoPax_srbf_varl2}c and f), therefore the $L_2^2(R_T)$ metric goes to zero.
Figure~\ref{fig:RoPax_srbf_varl2}e shows a zoom on the last iterations of the $\Delta\text{Var}[R_T]$ convergence: the component $\mathbb{E}[{R}_T]^2-\mathbb{E}[R_T^2]^2$ is converging to zero faster than $\mathbb{E}[{R}_T^2]-\mathbb{E}[R_T^{*2}]$, 
therefore their difference $\Delta \text{Var}[R_T]$ does not converge to zero.}

\begin{figure}[tp] 
	\centering
	\includegraphics[width = 0.47\textwidth]{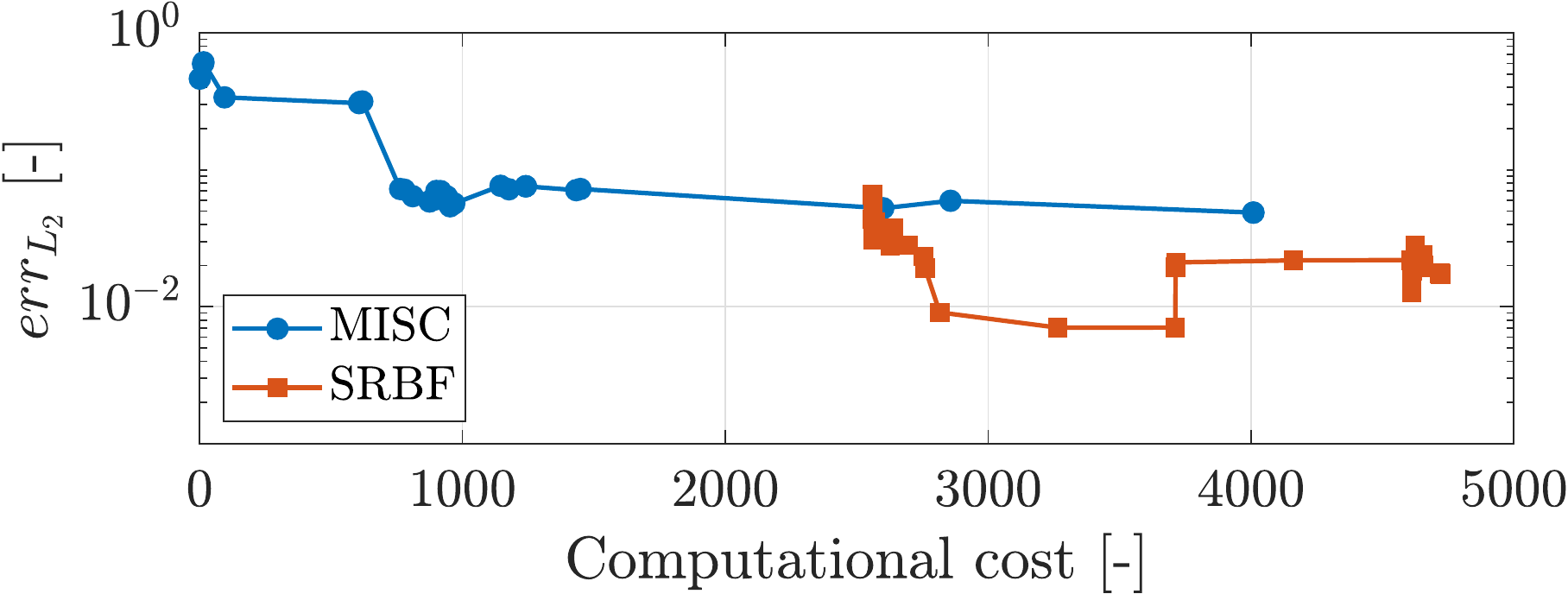}
	\includegraphics[width = 0.47\textwidth]{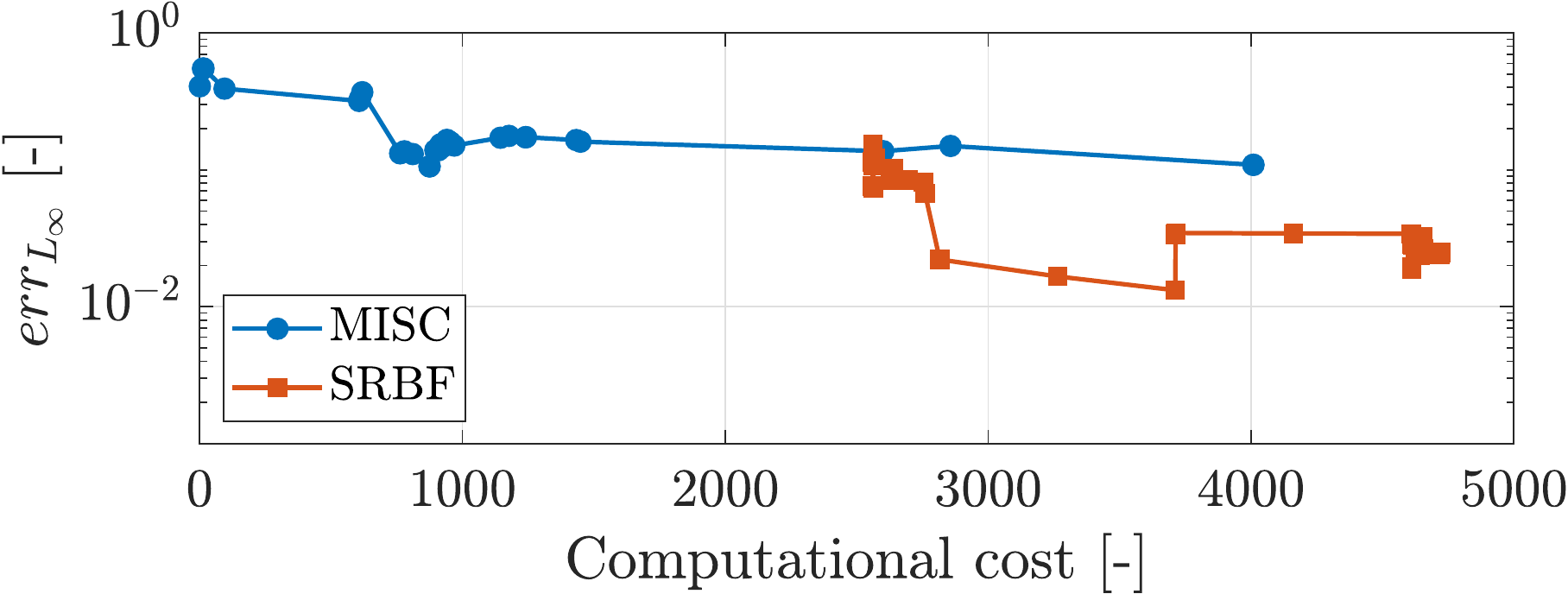}
	\caption{RoPax problem, comparison of MISC and SRBF methods: relative error of the approximation $R_T$ in $L_2$ (left) and $L_{\infty}$ norm (right) (see Eq.~\eqref{eq:err_norm}).}
	\label{fig:RoPax_norms}
\end{figure}

\begin{figure}[tp] 
	\centering
	\subfigure[Variance difference and $L_2^2$ metrics]{\includegraphics[width=0.3\textwidth]{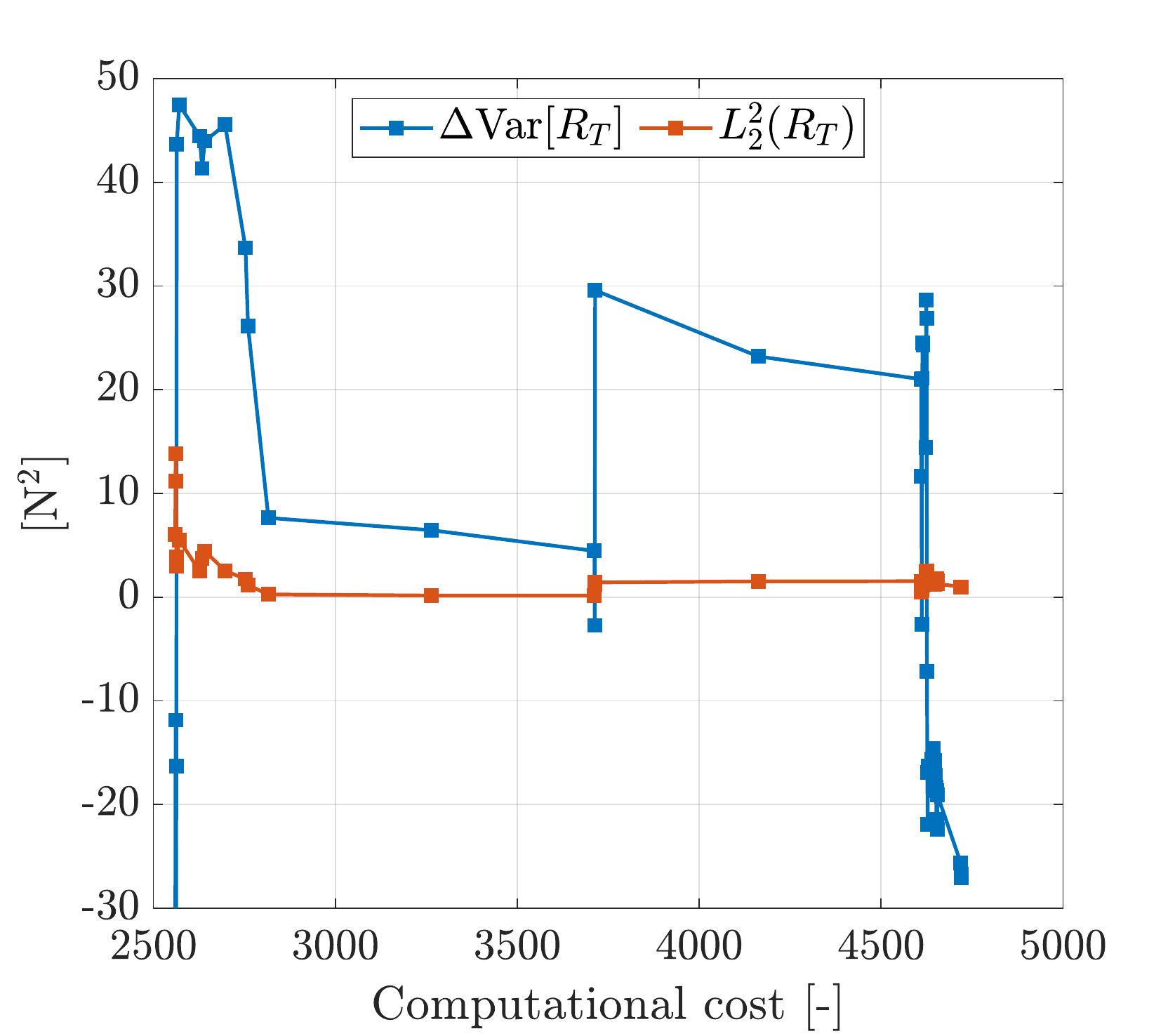}}
	\subfigure[Variance difference components]{\includegraphics[width=0.3\textwidth]{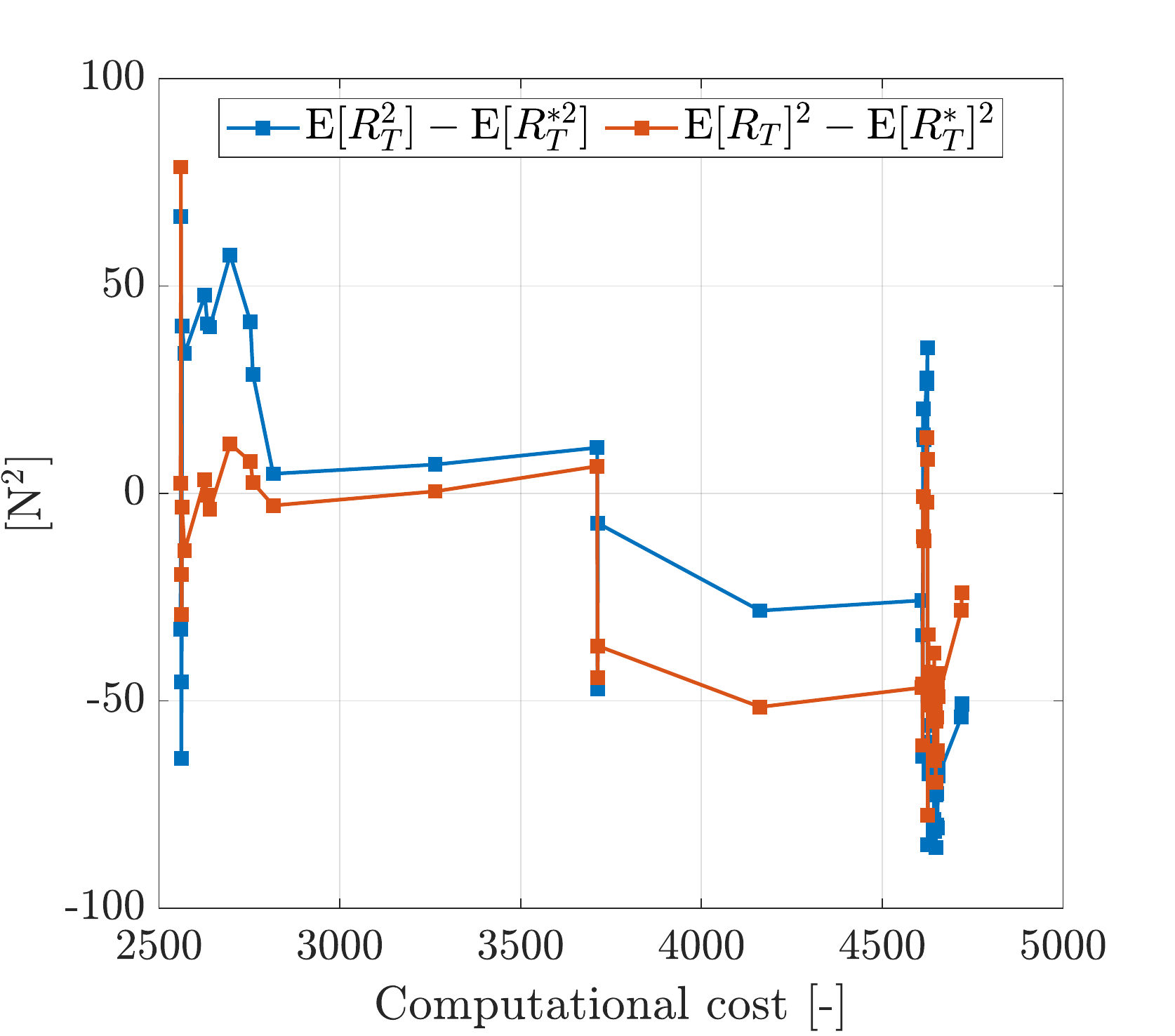}}
	\subfigure[$L_2^2$ metric components]{\includegraphics[width=0.3\textwidth]{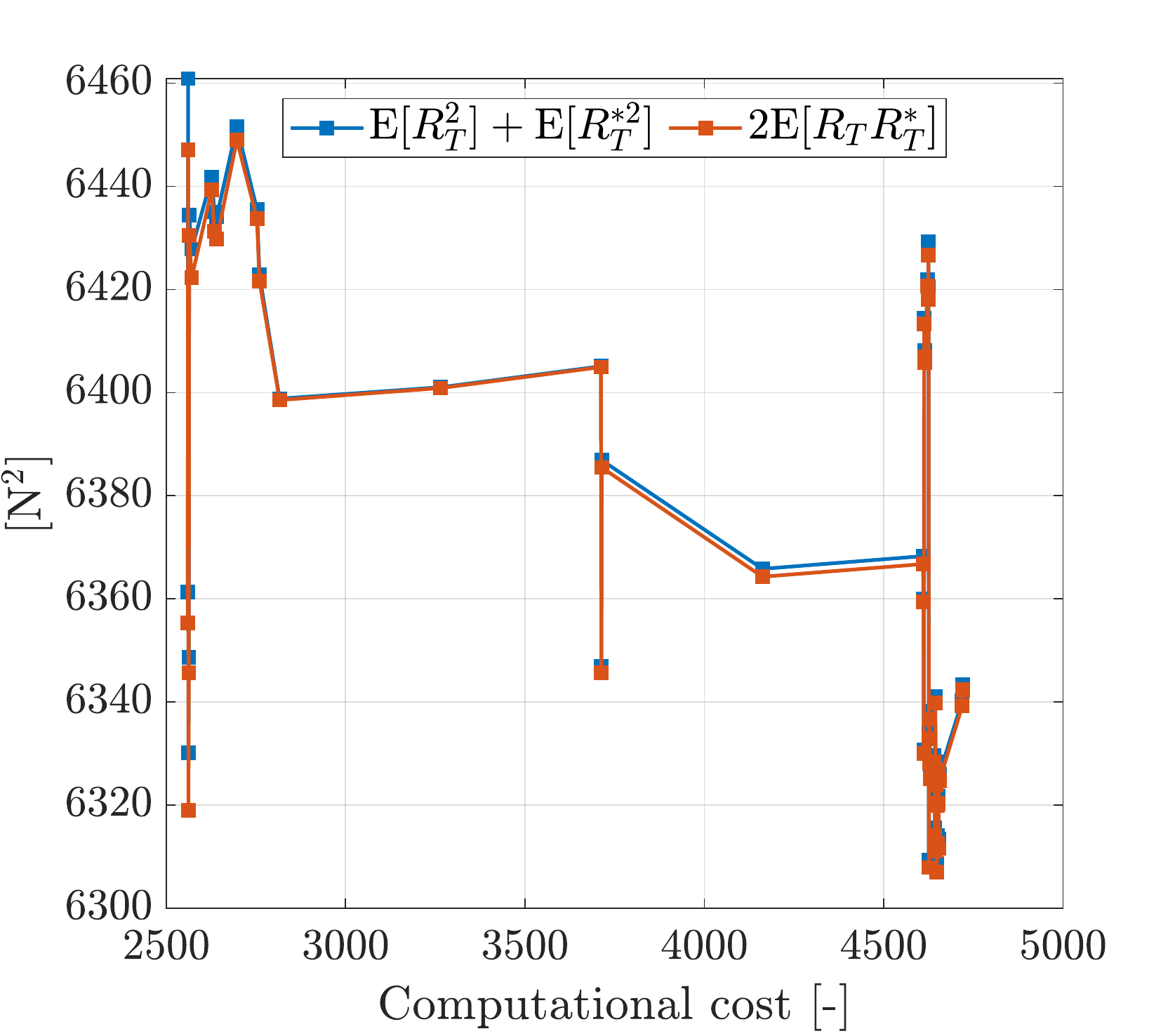}}
	\subfigure[Variance difference and $L_2^2$ metrics, detail]{\includegraphics[width=0.3\textwidth]{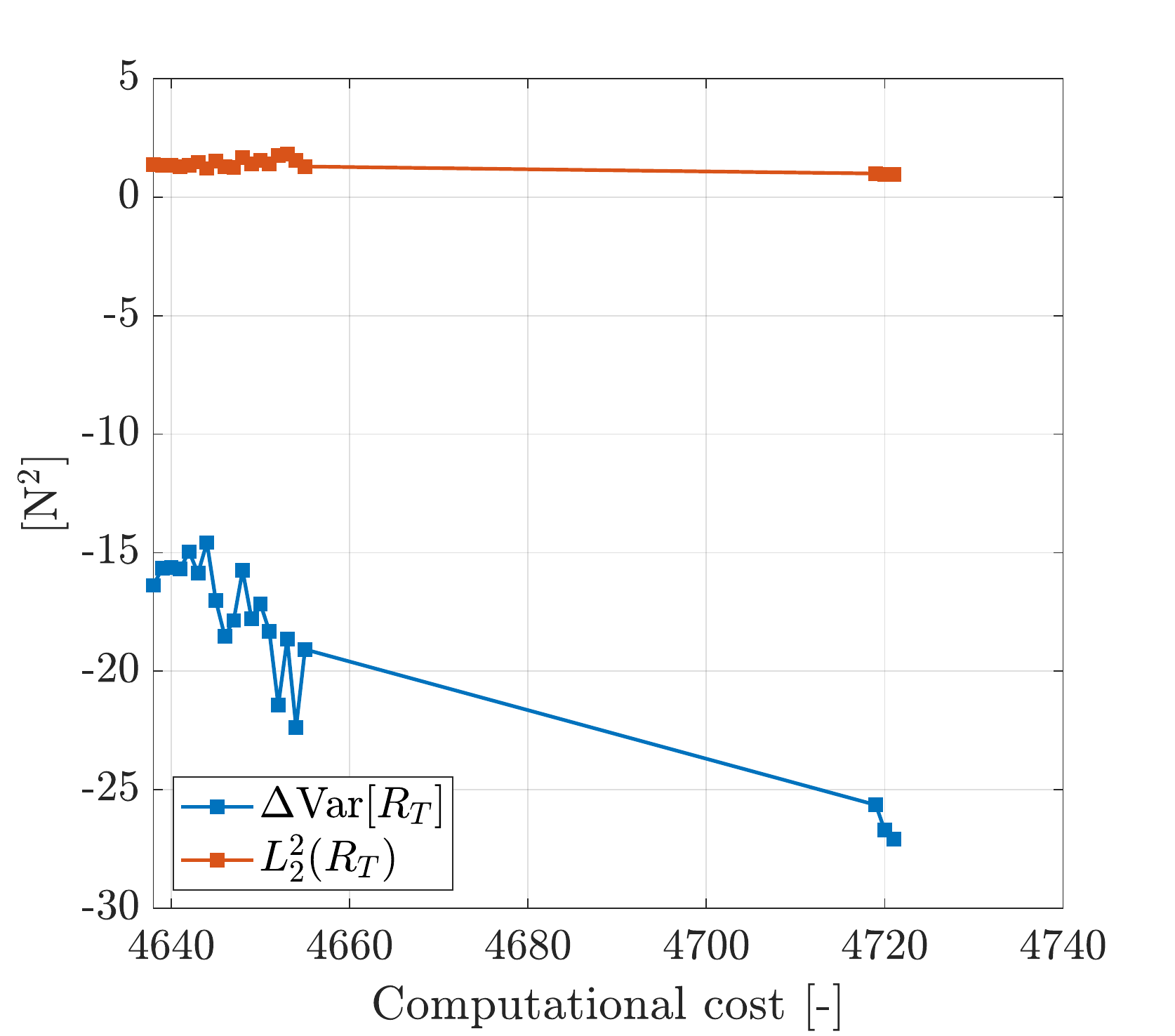}}
	\subfigure[Variance difference components, detail]{\includegraphics[width=0.3\textwidth]{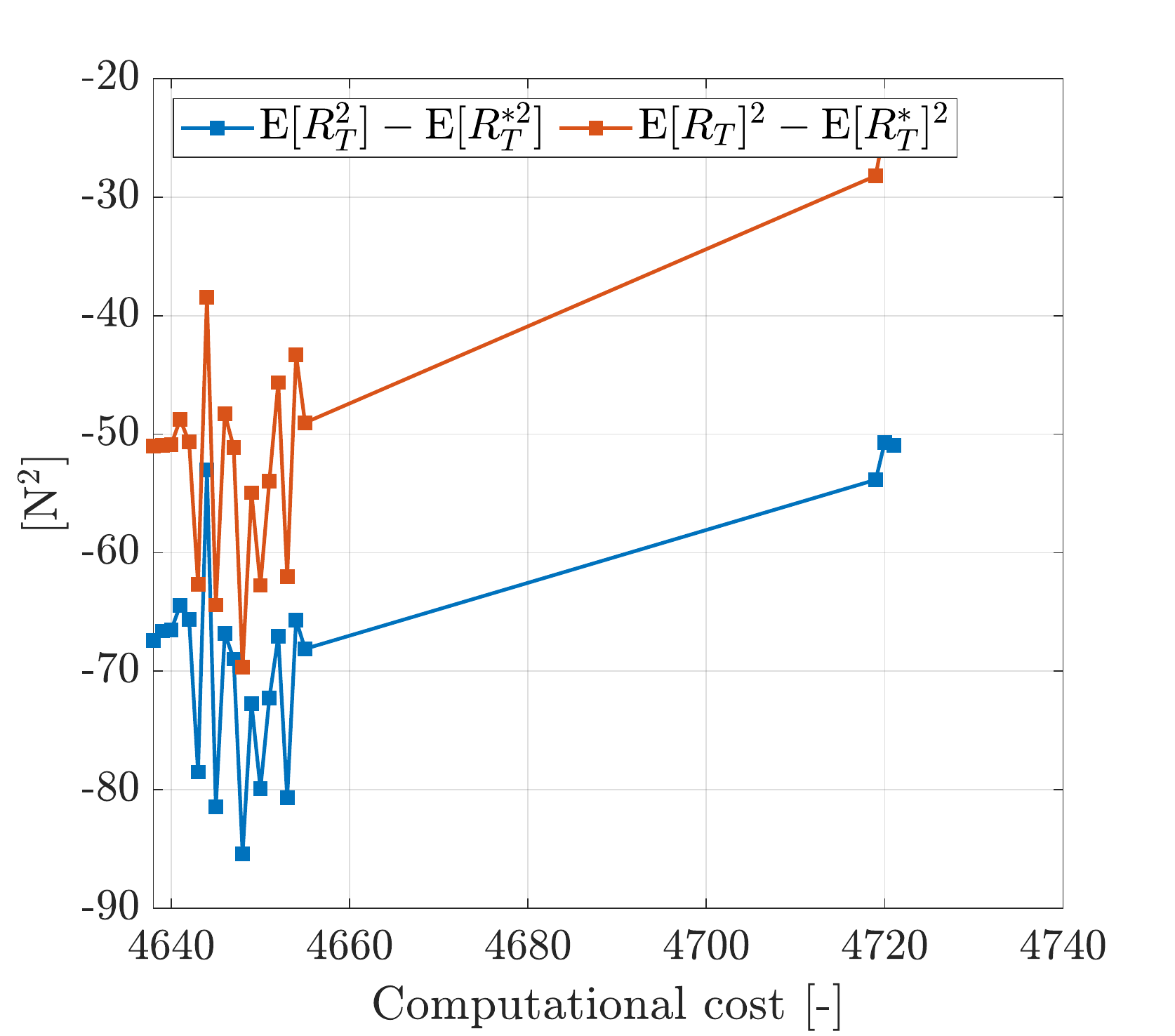}}
	\subfigure[$L_2^2$ metric components, detail]{\includegraphics[width=0.3\textwidth]{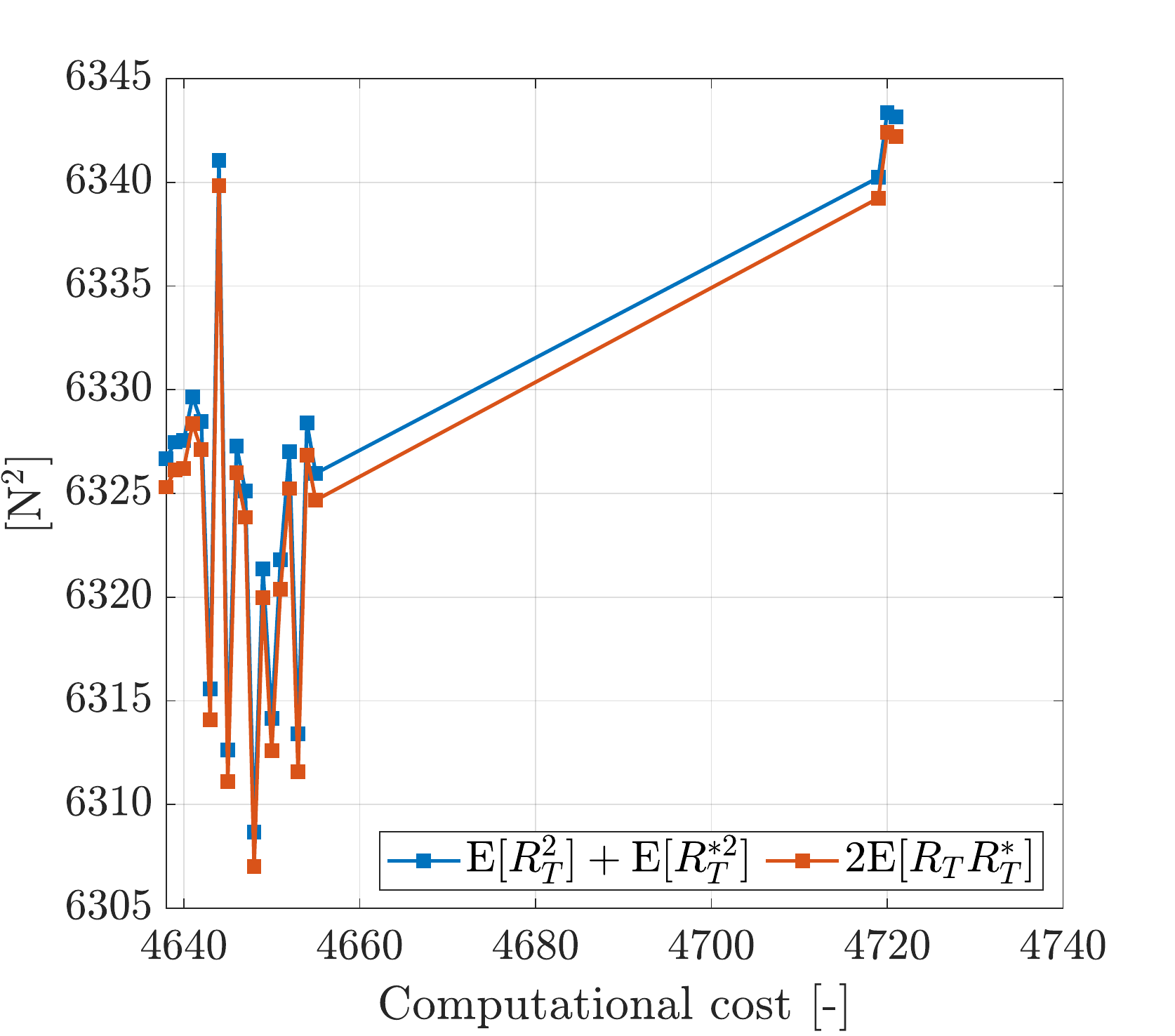}}
	\caption{RoPax problem, results for the SRBF method. Convergence of the variance difference and $L_2^2$ metrics and their components.}
	\label{fig:RoPax_srbf_varl2}
\end{figure}

Finally, in Fig.~\ref{fig:RoPax_pdf} the PDFs {obtained with both methods and with the reference surrogate, as well as} the results of the KS test statistic (cf.\ Eq.~\eqref{eq:KS_stat}), are plotted.
Both MISC and SRBF predict well the position of the mode of the PDF and its magnitude and the tails of the distribution, {with good agreement with the reference solution}. 
In the range $[40 \ 70]\text{N}$, {the PDF obtained by the MISC method is more ``wobbly'',
again due to the presence of noise in the evaluations of the solver, which corrupts the surrogate.} 
{Finally, the KS test statistic is seen to be convergent for both methods, implying}
convergence towards the reference CDF. Both convergences are not monotonic due to the influence of the noisy simulations.

\begin{figure}[tp]
	\centering
	\subfigure[PDF of $R_T$ at the final iteration (obtained by Matlab's {\it ksdensity} enforcing positive support since $R_T$ takes positive values only)]{\includegraphics[width = 0.47\textwidth]{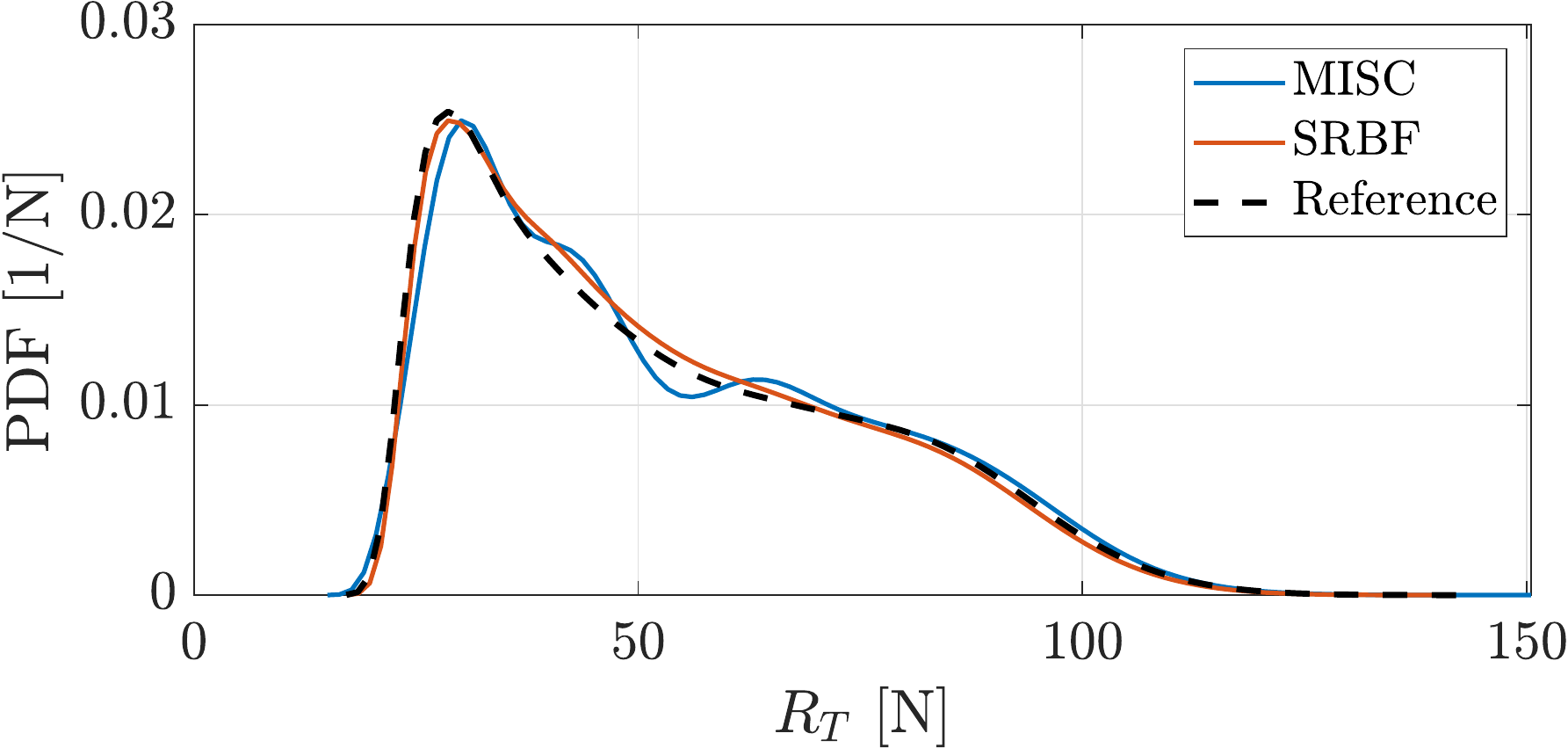}}
	\subfigure[KS test statistic (see Eq.~\eqref{eq:KS_stat})]{\includegraphics[width = 0.47\textwidth]{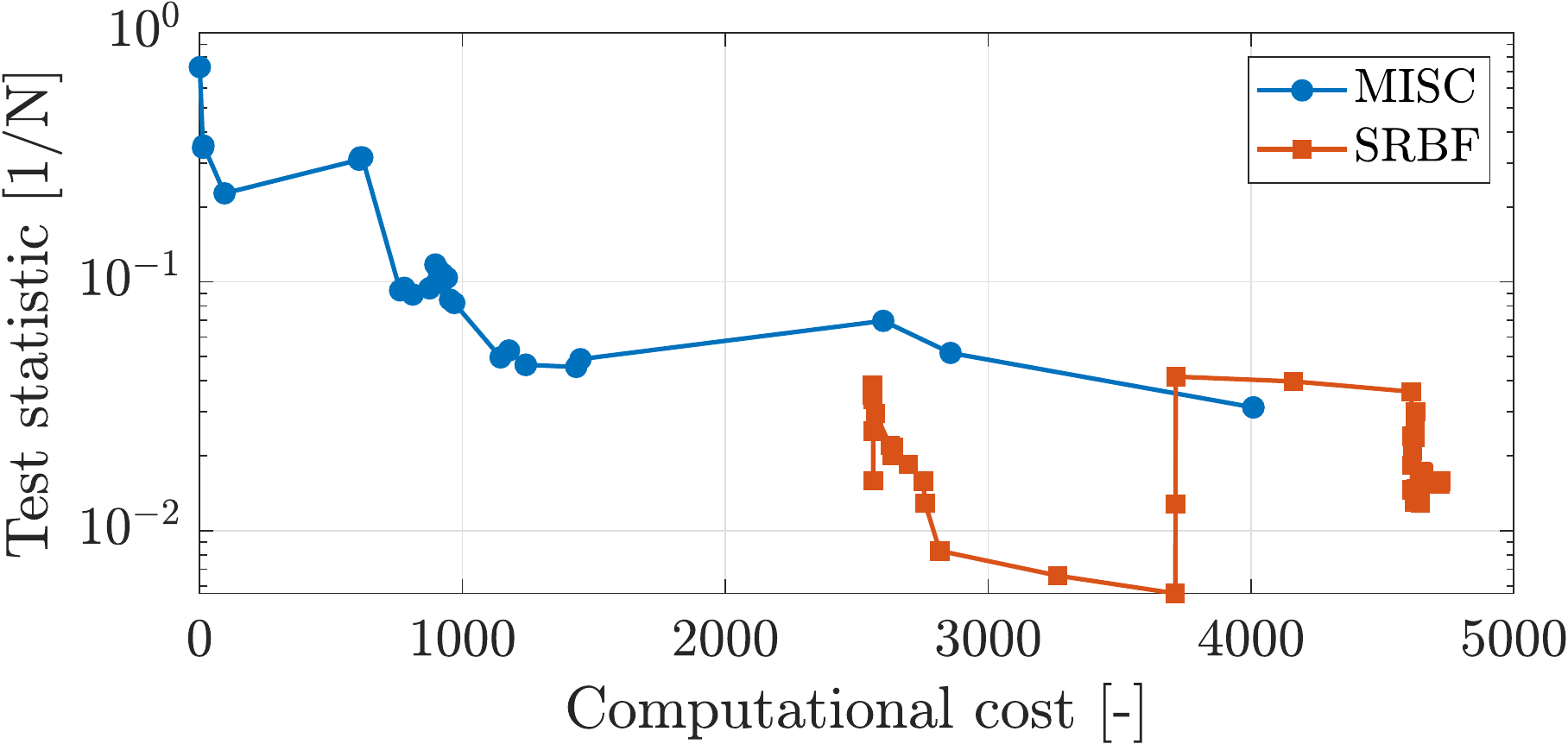}}
	\caption{RoPax problem, comparison of MISC and SRBF methods.}
	\label{fig:RoPax_pdf}
\end{figure}

\section{Conclusions and future work}\label{sect:conclusions}
{In this work, two multi-fidelity methods for forward UQ applications, MISC and SRBF, have been presented and applied to two examples to highlight their strengths and critical points. 
The first numerical test considered in this work is an analytical example, that served as benchmark for the results of the second test, which instead consists in a realistic application in naval engineering and it is more demanding for a number of reasons (noisy evaluations of the quantity of interest,  large setup time and computational costs). For the former, the different fidelities considered are Taylor expansions of increasing order of a given function,  while for the latter the fidelities are obtained by stopping the multi-grid computations in the RANS solver at different grid levels.}	
	
{In detail, we have considered the a-posteriori adaptive MISC method  already presented in \cite{jakeman2019adaptive}, with slight modifications on the profit computation,  and we have highlighted in passing that MISC is not an interpolatory method, contrary to its single-fidelity counterpart  (i.e., sparse-grids); this detail was never previously discussed (up to the authors' knowledge) in the MISC literature. 
SRBF has been used as an interpolatory surrogate model for the analytical test problem and as a regressive surrogate model \cite{wackers2020-AIAA} for the RoPax problem.}
  
{ For both tests, we have computed a number of error metrics for the quantity of interest (value of the function for the analytic test / advancement resistance for the ferry problem):  convergence of the approximation of the first four centered moments, mean squared and maximum prediction errors over a set of samples, and convergence of the CDF (as measured by the Kolmogorov--Smirnov test statistic).}

Overall, both MISC and SRBF confirmed to be viable multi-fidelity approaches to forward UQ problems. 
MISC has an edge in providing reasonable estimates of most statistical quantities at reduced computational cost, but is more sensitive to the noise in the evaluations of the  quantities of interest: {indeed, noise can strongly influence the adaptive selection process of the multi-indices, corrupt the estimates of e.g.\ higher-order moments (skewness, kurtosis), and introduce artifacts in the estimation of the PDF of the quantities of interest. With respect to the first issue, a quadrature-based adaptive criterion is expected to be more robust than a criterion based on the pointwise accuracy of the surrogate model.}
A possible strategy to mitigate the second issue, that consists in computing higher-order statistical  moments by taking Monte Carlo samples of the MISC surrogate model and then computing the moments from such set of values by sample formulas  (sample variance/skweness/kurtosis), has been proposed but it is not entirely satisfactory, since it is not clear how to choose  an appropriate number of samples (enough to be accurate, not too many to avoid resolving the scale of the noise).  This aspect deserves more investigations and is one of the subjects of future work.  Another practical issue is caused by the non-monotonic behavior of the profits,  where some indices with low profits shade useful neighbors, thus delaying the convergence of MISC.  More robust strategies to explore the set of multi-indices, that blend the profit-based selection of indices with other criteria  are also subject of future work; see e.g. \cite{chkifa:adaptive-interp,gerstner.griebel:adaptive}, where this problem  was discussed in the context of adaptive sparse-grids quadrature/interpolation.

Conversely, SRBF is less prone to the issue of noisy evaluations but on the other hand the current initialization strategy needs to sample all available fidelities, which results in a significantly larger initialization cost. Different initialization strategies are under investigation to reduce this gap: a possible approach would be, e.g., to build the initial training set using only one evaluation for each fidelity other than the lowest one, instead of $2N+1$ as in the current implementation, see \cite{pellegrini2021-MARINE}. This would allow the adaptive sampling method to freely define the best training set for the higher-fidelities. 

{Note that the discussion on computational costs in this work was set up in terms of sheer evaluation costs of the fidelities, i.e., neglecting the {CPU-time} related to the other operations required by MISC and SRBF. In details, MISC needs to keep track of profits, update the multi-index sets $\Lambda$ and $\text{Mar}(\Lambda)$, compute the combination technique coefficients in Eq. \eqref{eq:misc}, generate the tensor grids corresponding to the multi-indices in $\Lambda$ and $\text{Mar}(\Lambda)$, and evaluate the associated multivariate Lagrange polynomials. Conversely, SRBF requires to solve the linear system \eqref{eq:lsrbf} at each iteration (the dimension of the linear system grows at each iteration; furthermore, when the LOOCV procedure is performed, the linear system is solved for each tested value of the number of centers $\mathcal{K}$), as well as to solve the problem in Eq. \eqref{eq:MUAS}.} 

An aspect worth investigating for both methods is how to incorporate {available} soft information (monotonicity, multimodality, etc.)   on the physical nature of the problem in the construction of the surrogate models. In the particular RoPax problem considered here, the resistance is expected to be    monotone increasing with respect to advancement speed and draught: such property could be preserved   by employing, e.g., least-squares regressions with appropriate polynomial degrees and/or monotonic smoothing, see e.g. \cite{dykstra1982}.
  
Concerning the application to naval UQ, once fixed the current limitations of both methods,   future research will address more complex test cases, such as a larger number of uncertain parameters and more realistic conditions, to take, e.g., regular/irregular waves into account.

\section*{Acknowledgments}
CNR-INM is grateful to Dr. Woei-Min Lin, Dr. Elena McCarthy, and Dr. Salahuddin Ahmed of the Office of Naval Research and Office of Naval Research Global, for their support through NICOP grant N62909-18-1-2033. 
Dr. Riccardo Pellegrini is partially supported through CNR-INM project OPTIMAE. 
The HOLISHIP project (HOLIstic optimisation of SHIP design and operation for life cycle, www.holiship.eu) is also acknowledged, funded by the European Union’s Horizon 2020 research and innovation program under grant agreement N. 689074. 
Lorenzo Tamellini and Chiara Piazzola have been supported by the PRIN 2017 project 201752HKH8 ``Numerical Analysis for Full and Reduced Order Methods for the efficient and accurate solution of complex systems governed by Partial Differential Equations (NA-FROM-PDEs)''. 
Lorenzo Tamellini also acknowledges the support of GNCS-INdAM (Gruppo Nazionale Calcolo Scientifico - Istituto Nazionale di Alta Matematica). 
Finally, CNR-INM is grateful to CINECA for providing HPC capabilities through the ISCRA-C grant HP10CRM564 ``Enabling high-fidelity Shape Optimization by Design-space Augmented Dimensionality Reduction (ESODADR)''.

\bibliographystyle{elsarticle-num}
\section*{\refname}
\bibliography{biblio}

\end{document}